# SPECTRAL ANALYSIS OF SINAI'S WALK FOR SMALL EIGENVALUES[1]

By Anton Bovier and Alessandra Faggionato

*Weierstrass Institut für Angewandte Analysis und Stochastik and Technische Universität Berlin, and Università "La Sapienza"*

Sinai's walk can be thought of as a random walk on $\mathbb{Z}$ with random potential $V$, with $V$ weakly converging under diffusive rescaling to a two-sided Brownian motion. We consider here the generator $\mathbb{L}_N$ of Sinai's walk on $[-N,N]\cap\mathbb{Z}$ with Dirichlet conditions on $-N,N$. By means of potential theory, for each $h>0$, we show the relation between the spectral properties of $\mathbb{L}_N$ for eigenvalues of order $o(\exp(-h\sqrt{N}))$ and the distribution of the $h$-extrema of the rescaled potential $V_N(x)\equiv V(Nx)/\sqrt{N}$ defined on $[-1,1]$. Information about the $h$-extrema of $V_N$ is derived from a result of Neveu and Pitman concerning the statistics of $h$-extrema of Brownian motion. As first application of our results, we give a proof of a refined version of Sinai's localization theorem.

**1. Introduction.** Random walks in random environments are a major paradigm for the dynamics of systems in complex environments (see [23] for a recent in depth review). One of the simplest special cases is the one-dimensional nearest-neighbor random walk with i.i.d transition probabilities, $p_x, 1-p_x$, in the regime where $\mathbf{E}\ln\frac{p_x}{1-p_x}=0$ and $\mathbf{E}\ln^2\frac{p_x}{1-p_x}>0$ . In this regime Sinai [21] discovered remarkable slowing down of the diffusive time scale. Since then, the model was investigated very intensely and in great detail both in the probabilistic and the physics literature; see, for example, [6, 7, 8, 10, 11, 12, 13, 15, 16, 20]. Rather recently [8], this model was considered from the point of view of the popular concept of *ageing* which is a particular manifestation of the slow down of the dynamics characterized by a particular behavior of autocorrelation functions. It was shown that ageing results, in this model, rather directly from Sinai's localization theorem that we shall explain below. Another approach toward the characterization of slow dynamics would be through the spectral properties of the generator of the process. In a recent paper we have carried this out in full detail in the simplest possible model, Bouchaud's trap model on the complete graph [4] where we have shown, in particular, that all the standard ageing properties of the model can be derived easily from spectral data. Recently, the spectrum of the generator of Sinai's random walk was analyzed in [10, 15] using renormalization group methods. In the present paper we give a refined and fully rigorous analysis of the bottom part of the spectrum of Sinai's random walk and show that this leads to a very easy proof of a (refined) version of Sinai's localization theorem. Another application of the spectral information will show a limit law that expresses the fact that Sinai's random walk can be seen as a process that on an infinite sequence of (random) time scales *appears* to be approaching equilibrium exponentially. Let us note that Comets and Popov [7] have used control of principal eigenvalues of the generator of Sinai's walk in suitable intervals to obtain moderate deviation results.

[1]Supported in part by the DFG in the Dutch–German Bilateral Research Group "Mathematics of Random Spatial Models from Physics and Biology."

Let us note that the spectral analysis of the generator can also be considered as that of a corresponding quantum mechanical Schrödinger operator. This operator has been considered in the context of two-dimensional electrons in a particular random magnetic field and as an effective Hamiltonian of polyacetylene (see [1] for a discussion and references).

1.1. *Sinai's walk. Definitions and key facts.* Before stating our results, let us fix the notation. We define an environment as a sequence, $\omega = \{\omega_x\}_{x\in\mathbb{Z}}$ with $\omega_x \in [0,1]$. For a given environment, $\omega$, Sinai's walk $(X_n, n \geq 0)$ is a discrete time random walk on $\mathbb{Z}$ with transition probabilities

$$\begin{aligned} &\mathrm{Prob}(X_{n+1} = x+1 | X_n = x) = \omega_x, \\ &\mathrm{Prob}(X_{n+1} = x-1 | X_n = x) = 1-\omega_x. \end{aligned} \tag{1.1}$$

We use $\mathrm{P}_x^\omega$ to denote the law of the random walk $(X_n, n \geq 0)$ starting at $x \in \mathbb{Z}$.

We will consider random environments consisting of i.i.d. sequences of random variables, $\omega_x$, $x \in \mathbb{Z}$, whose law will be denoted by $\mathbf{P}$. We will make the usual ellipticity assumption that, for some $\kappa > 0$,

$$\omega_x \in [\kappa, 1-\kappa] \qquad \forall x \in \mathbb{Z}. \tag{1.2}$$

We set $\Omega \equiv [\kappa, 1-\kappa]^{\mathbb{Z}}$. To be in the situation of Sinai's walk, we assume further that

$$\mathbf{E}\left(\ln\left(\frac{\omega_x}{1-\omega_x}\right)\right) = 0, \tag{1.3}$$

where $\mathbf{E}$ denotes the expectation w.r.t. to $\mathbf{P}$, and

$$\sigma^2 \equiv \mathbf{E}\left[\ln^2\left(\frac{\omega_x}{1-\omega_x}\right)\right] > 0. \tag{1.4}$$

Let us finally define the measure $\mathbf{P}_x \equiv \mathbf{P} \otimes \mathrm{P}_x^\omega$ on $\Omega \times \mathbb{Z}^{\mathbb{N}}$ as

$$\mathbf{P}_x(F \times G) = \int_F \mathrm{P}_x^\omega(G)\mathbf{P}(d\omega) \qquad \forall F \in \mathcal{F}, G \in \mathcal{G},$$

where $\mathcal{F}$, $\mathcal{G}$ are respectively the $\sigma$-algebra of Borel subsets of $\Omega$ and $\mathbb{Z}^{\mathbb{N}}$.

In this setting, it is well known (see, e.g., [23]) that the random walk is recurrent $\mathbf{P}_0$-almost surely. Moreover, Sinai [21] proved that there exists a function, $\mathfrak{m}^{(n)}(\omega)$, depending only on the environment, such that

$$\frac{X_n}{\ln^2 n} - \mathfrak{m}^{(n)} \to 0 \qquad \text{in } \mathbf{P}_0\text{-probability}, \tag{1.5}$$

as $n \to \infty$. A refinement of the above localization result was obtained by Golosov for a slightly modified random walk [12]: for a suitable distribution function $F$,

$$\lim_{n\uparrow\infty} \mathbf{P}_0(X_n - \mathfrak{m}^{(n)} \ln^2 n \le y) = F(y) \qquad \forall y \in \mathbb{R}, \tag{1.6}$$

namely, under $\mathbf{P}_0$, the random variable $X_n - \mathfrak{m}^{(n)} \ln^2 n$ converges in law. Moreover, as shown independently in [11] and [13], the distribution of the random variable $\sigma^2\mathfrak{m}^{(n)}(\omega)$ under $\mathbf{P}$ converges weakly as $n \to \infty$ to a suitable functional $L$ of the Brownian motion with

$$\frac{d\mathrm{Prob}[L \le x]}{dx} = \frac{2}{\pi}\sum_{k=0}^{\infty} \frac{(-1)^k}{(2k+1)} \exp\left\{-\frac{(2k+1)^2\pi^2}{8}|x|\right\}.$$

Sinai's walk can be thought of as a random walk on $\mathbb{Z}$ with random potential. Namely, define the potential $V(x)$, $x \in \mathbb{Z}$, as

$$V(x) = \begin{cases} \sum_{i=1}^{x} \ln\frac{1-\omega_i}{\omega_i}, & \text{if } x \ge 1, \\ 0, & \text{if } x = 0, \\ -\sum_{i=x+1}^{0} \ln\frac{1-\omega_i}{\omega_i}, & \text{if } x \le -1. \end{cases} \tag{1.7}$$

Then, the jump probabilities can be expressed as

$$\omega_x = e^{-\nabla V(x)/2}/Z, \qquad 1-\omega_x = e^{\nabla V(x)/2}/Z, \tag{1.8}$$

where $Z$ denotes the normalizing constant and $\nabla V(x) \equiv V(x) - V(x-1)$. The behavior of the potential $V$ is well described by Donsker's invariance principle. Given $N \in \mathbb{Z}_+$, define the rescaled potential $V^{(N)} \in C(\mathbb{R})$ as

$$V^{(N)}(t) \equiv \frac{s}{\sqrt{N}} V(k) + \frac{1-s}{\sqrt{N}} V(k+1) \qquad \text{if } t = s\frac{k}{N} + (1-s)\frac{k+1}{N}, k \in \mathbb{Z}, s \in [0,1]. \tag{1.9}$$

For later applications, note that $V^{(N)}$ is a Lipschitz function with Lipschitz constant $c(\kappa)\sqrt{N}$. Due to the independence of $\{\omega_x\}_{x\in\mathbb{Z}}$ and assumptions (1.3) and (1.4), endowing the space $C(\mathbb{R})$ with the topology of uniform convergence on compact subsets, by the Donsker's invariance principle the random path $V^{(N)}$ converges in distribution to $B = (B_t, t \in \mathbb{R})$, the two-sided Brownian motion with $B_0 = 0$ and variance $\sigma^2$.

The Komlós–Major–Tusnády strong approximation theorem [14] (see also Proposition 3) gives an even stronger result: given $L > 0$, there exist positive constants $C_1, C_2, C_3$ such that, for each $N \in \mathbb{Z}_+$, there exists a coupling on an enlarged probability space between $(V^{(N)}(x), x \in [-L, L])$ and the two-sided Brownian motion $B$ with variance $\sigma$ such that

$$P^{(N)}\left(\sup_{x\in[-L,L]} |V^{(N)}(x) - B_x| > \frac{C_1 \ln N}{\sqrt{N}}\right) < \frac{C_2}{N^{C_3}}. \tag{1.10}$$

*Notational warning*: in what follows, $c(\kappa)$ will denote a generic constant depending only on $\kappa$ [see (1.2)] and it can change from expression to expression.

1.2. *Generators with Dirichlet conditions.* Our objective will be to control the spectrum of the generator of Sinai's walk with Dirichlet conditions outside a (large) interval $\{-N+1, \dots, N-1\}$. We write $\mathbb{P} \equiv \mathbb{P}(\omega)$ for the transition matrix of the random walk for a fixed environment. For $D \subset \mathbb{Z}$, we define the transition matrix with Dirichlet conditions outside $D$ as $\mathbb{P}(D) \equiv (\mathbb{P}_{x,y})_{x,y\in D}$. It is convenient to define the "generator," $\mathbb{L}$, of the discrete-time chain as $\mathbb{L} \equiv \mathbb{I} - \mathbb{P}$, as well as the corresponding Dirichlet operators $\mathbb{L}(D)$. Note that $\mathbb{L}(D)$ is the restriction to $D$ of the generator of Sinai's walk killed when it leaves $D$.

Given $u \in \mathbb{R}^D$, let us define $\tilde{u} \in \mathbb{R}^{\mathbb{Z}}$ as $\tilde{u} \equiv u\mathbb{I}_D$, then $(\mathbb{L}(D)u)(x) = (\mathbb{L}\tilde{u})(x)$ for any $x \in D$. In particular, $\lambda$ is an eigenvalue of $\mathbb{L}(D)$, shortly $\lambda \in \sigma(\mathbb{L}(D))$, iff $\exists v \in \mathbb{R}^{\mathbb{Z}}$ such that

$$\begin{cases} (\mathbb{L}-\lambda)v(x) = 0, & \text{if } x \in D, \\ v(x) = 0, & \text{if } x \notin D. \end{cases} \tag{1.11}$$

Identifying $v$ with $v_{|D}$, we say that $v$ satisfying (1.11) is an eigenvector of $\mathbb{L}(D)$ with eigenvalue $\lambda$.

Let us first describe some simple spectral results concerning $\mathbb{L}(D)$. Note that the measure $\mu$ on $\mathbb{Z}$ defined as

$$\mu(x) \equiv e^{-V(x)}/\omega_x \qquad \forall x \in \mathbb{Z},$$

satisfies

$$\mu(x)\omega_x = \mu(x+1)(1-\omega_{x+1}) = e^{-V(x)} \qquad \forall x \in \mathbb{Z}. \tag{1.12}$$

In particular, it is a reversible measure for $\mathbb{L}(D)$ for all $D \subset \mathbb{Z}$, that is, $\mathbb{L}(D)$ is a symmetric operator on $\mathbb{L}^2(D,\mu)$ having left eigenvector $\mu u$ with eigenvalue $\lambda$ whenever $u$ is a (right) eigenvector with eigenvalue $\lambda$. Moreover, denoting by $(\cdot,\cdot)$ the scalar product on $\mathbb{L}^2(\mathbb{Z},\mu)$, one easily obtains for all $f \in L^2(\mathbb{Z},\mu)$ that the *Dirichlet form* is given by the expression

$$(f,\mathbb{L}f) = \sum_{x\in\mathbb{Z}} \mu(x)\omega_x (f(x+1)-f(x))^2. \tag{1.13}$$

*Periodicity.* Note that the Markov chains we are defining are periodic. Define $\Sigma_{\rm o}$ $[\Sigma_{\rm e}]$ the subspace of $\mathbb{R}^D$ having even [odd] coordinates equal to zero. Trivially, $\mathbb{R}^D = \Sigma_{\rm o} \oplus \Sigma_{\rm e}$ and $\mathbb{P}(\Sigma_{\rm o}) \subset \Sigma_{\rm e}$, $\mathbb{P}(\Sigma_{\rm e}) \subset \Sigma_{\rm o}$. This implies the following a-priori information on the spectra, whose proof is left to the reader:

LEMMA 1. *Let $D \equiv [a,b]\cap\mathbb{Z}$. Then the matrix $\mathbb{P}(D)$ has simple eigenvalues $-1 < \lambda_1 < \lambda_2 < \cdots < \lambda_{|D|} < 1$ and $\lambda_i = -\lambda_{|D|-i+1}$ for all $i : 1 \le i \le |D|$. Moreover, if $\mathbb{P}\psi = \lambda\psi$, where $\psi = \psi_{\rm o} + \psi_{\rm e}$ with $\psi_{\rm o} \in \Sigma_{\rm o}$, $\psi_{\rm e} \in \Sigma_{\rm e}$, then $\mathbb{P}\psi' = -\lambda\psi'$ where $\psi' = \psi_{\rm o} - \psi_{\rm e}$.*

1.3. *$h$-extrema and saddles.* The small eigenvalues of the generators will be labeled by the deep minima of the potential. This will require some further notation. Given a continuous path $\gamma \in C([-1,1])$, we say that $x \in [-1,1]$ is a *$h$-minimum* (for $\gamma$) if there exist $a,b \in [-1,1]$ with

$$\begin{aligned} &a < x < b, \\ &\gamma(a) \ge \gamma(x) + h, \qquad \gamma(b) \ge \gamma(x) + h \quad \text{and} \quad \gamma(x) = \min_{[a,b]} \gamma. \end{aligned} \tag{1.14}$$

We say that $x \in [-1,1]$ is a *$h$-maximum* (for $\gamma$) if one of the following three complementary conditions is satisfied:

(i) $x$ is a $h$-minimum for $-\gamma$,

(ii) $\exists b \in (x,1]$ such that $\gamma(x) - \gamma(b) \ge h$, $\gamma(x) = \max_{[-1,b]}\gamma$ and $\min_{[-1,x]}\gamma > \gamma(x) - h$,

(iii) $\exists a \in [-1,x)$ such that $\gamma(x) - \gamma(a) \ge h$, $\gamma(x) = \max_{[a,1]}\gamma$ and $\min_{[x,1]}\gamma > \gamma(x) - h$.

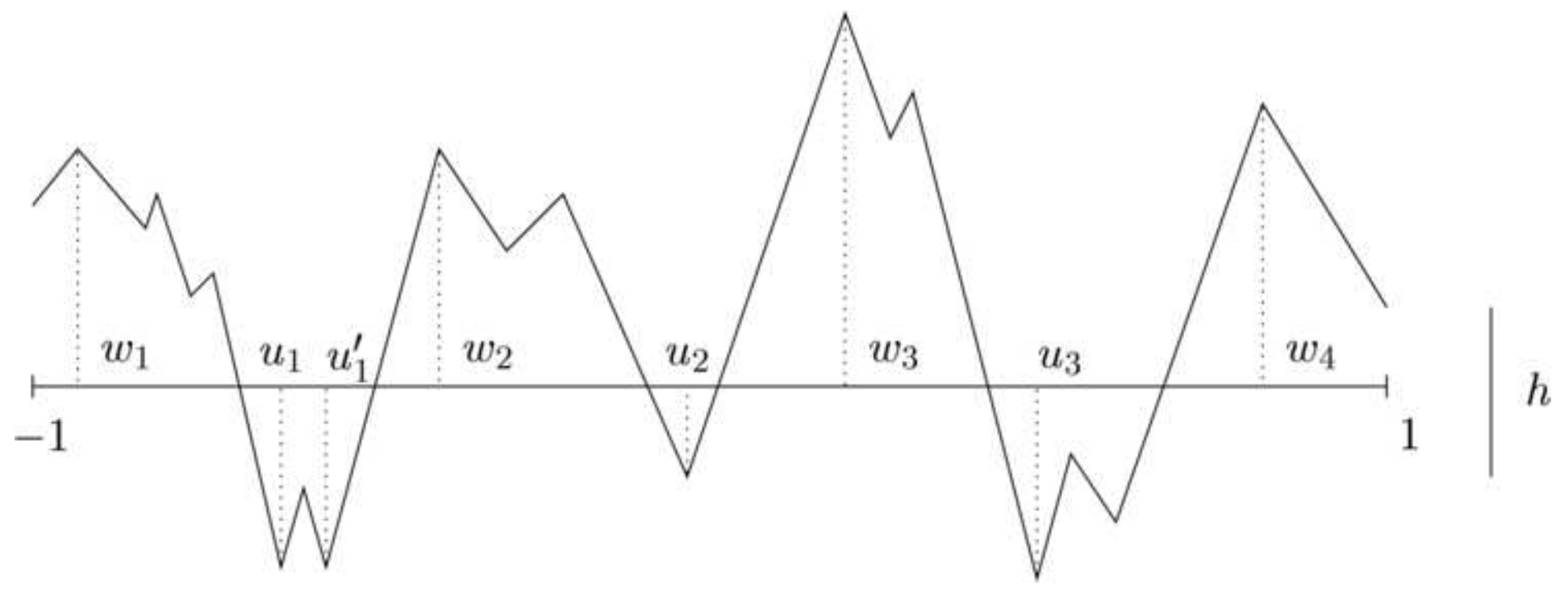

FIG. 1. *$h$-extrema of a path in $C([-1,1])$.*

See, for example, Figure 1 where points $u_1, u'_1, u_2, u_3$ are $h$-minima, while $w_1, w_2, w_3, w_4$ are $h$ maxima.

When considering $\gamma \in C(\mathbb{R})$, we say that $x \in \mathbb{R}$ is a $h$-minimum (for $\gamma$) if there exist $a, b \in \mathbb{R}$ satisfying (1.14) and we say that $x \in \mathbb{R}$ is a $h$-maximum (for $\gamma$) if $x$ is a $h$-minimum for $-\gamma$.

In what follows we take $\gamma \in C(I)$ with $I = [-1,1]$ or $I = \mathbb{R}$. A point $x \in I$ is called a *$h$-extremum* if it is a $h$-minimum or a $h$-maximum. We write $\mathcal{M}_h^-(\gamma)$, $\mathcal{M}_h^+(\gamma)$ and $\mathcal{E}_h(\gamma)$, respectively for the sets of $h$-minima, $h$-maxima and $h$-extrema of $\gamma$.

Given $x, x' \in \mathcal{M}_h^\pm(\gamma)$, we say that they are equivalent, $x \sim x'$, if

$$\max_{z \in [x \wedge x', x \vee x']} |\gamma(z) - \gamma(x)| < h.$$

Note that $\gamma(x) = \gamma(x')$ whenever $x \sim x'$ and that $z \sim x$ if $z \in \mathcal{M}_h^\pm(\gamma), x \sim x'$ and $z \in [x \wedge x', x \vee x']$. One can easily prove that each equivalence class is a closed subset of $I$ and for each compact subset $K \subset I$, $K$ intersects a finite number of equivalence classes. We will denote by $M_h^\pm(\gamma)$ the subset of $\mathcal{M}_h^\pm(\gamma)$ obtained by taking for each equivalence class in $\mathcal{M}_h^\pm(\gamma)/\sim$ the smallest element (if it exists). Note that if $I = [-1,1]$, then $|M_h^\pm(\gamma)| < \infty$. Finally, the piece of $\gamma$ between consecutive $h$-maxima in $M_h^+(\gamma)$ will be called *$h$-valley*.

One can easily prove the following lemma (see Figure 1):

LEMMA 2. *Given $\gamma \in C([-1,1])$, if $M_h^-(\gamma) = \{u_1, \dots, u_q\}$ with $q \geq 1$ and $u_1 < u_2 < \cdots < u_q$, then $M_h^+(\gamma) = \{w_1, w_2, \dots, w_{q+1}\}$ with*

$$-1 \leq w_1 < u_1 < w_2 < u_2 < \cdots < w_q < u_q < w_{q+1} \leq 1.$$

*Moreover, for all $i \in \{1, \dots, q\}$ and $j \in \{2, \dots, q\}$,*

$$\gamma(u_i) = \min_{[w_i, w_{i+1}]} \gamma, \gamma(w_1) = \max_{[-1,u_1]} \gamma, \gamma(w_j) = \max_{[u_{j-1},u_j]} \gamma, \gamma(w_{q+1}) = \max_{[u_q,1]} \gamma.$$

Given $\gamma \in C(I)$ and disjoint finite sets $A, B \subset I$, we define $Z(A,B)$ as the set of *saddle* points between $A$ and $B$:

$$Z(A,B) \equiv \Big\{ z \in I : \exists a \in A, b \in B \text{ with } a \wedge b \le z \le a \vee b \\ \text{and } \gamma(z) = \min_{a\in A, b \in B} \max_{a \wedge b \le x \le a \vee b} \gamma(x) \Big\}.$$

Moreover, we set

$$z^*(A,B) \equiv \min(Z(A,B))$$

[the definition is well posed since $Z(A,B)$ is compact]. Note that $z^*(A,B) \in M_h^+(\gamma)$ whenever $A, B \subset M_h^-(\gamma)$.

Finally, given $h, \delta > 0$, we define the family of *good paths* in $C([-1,1])$, $\mathcal{A}_{h,\delta}$, as the set of paths $\gamma$ satisfying the following conditions:

1. $M_h^-(\gamma) \neq \varnothing$,
2. 

$$\gamma(z^*(x, M_h^-(\gamma) \cup \{-1,1\} \setminus \{x\})) - \gamma(x) \ge h + \delta \qquad \forall x \in M_h^-(\gamma), \tag{1.15}$$

3. for a suitable labeling $M_h^-(\gamma) = \{x_1, x_2, \dots, x_q\}$,

$$\gamma(z^*(x_k, S_{h,k-1})) - \gamma(x_k) \ge \max_{q \ge j > k} \{\gamma(z^*(x_j, S_{h,k-1})) - \gamma(x_j)\} + \delta \\ \forall k : 1 \le k \le q-1, \tag{1.16}$$

where

$$\begin{cases} S_{h,k} = \{-1,1\}, & \text{if } k = 0, \\ S_{h,k} = \{x_1, x_2, \dots, x_k\} \cup \{-1,1\}, & \text{if } 1 \le k \le q. \end{cases}$$

Condition (1.16) is a *nondegeneracy* condition. It can be read as follows (see Figure 2): $(x_1, \gamma(x_1))$ is the most trapped starting point in $\gamma$ for a walker desiring to reach one of the points $(-1, \gamma(-1))$ and $(1, \gamma(1))$. Then $(x_2, \gamma(x_2))$ is the most trapped starting point in $\gamma$ for a walker desiring to reach one of the points $(-1, \gamma(-1))$, $(1, \gamma(1))$ and $(x_1, \gamma(x_1))$, and so on.

We note that if $\gamma \in \mathcal{A}_{h,\delta}$, then the above labeling $M_h^-(\gamma) = \{x_1, x_2, \dots, x_q\}$ is unique. In what follows, when assuming $\gamma \in A_{h,\delta}$, we will always think of $x_1, \dots, x_q$ as this labeling of $M_h^-(\gamma)$. In particular, $q = |M_h^-(\gamma)|$. It is also convenient to set (see Figure 2)

$$d_k(\gamma) = \gamma(z^*(x_k, S_{h,k-1})) - \gamma(x_k) \qquad \forall 1 \le k \le q. \tag{1.17}$$

Then it is easy to prove that condition (1.16) is equivalent to the following one:

$$d_k(\gamma) \ge d_{k+1}(\gamma) + \delta \qquad \forall k : 1 \le k \le q-1. \tag{1.18}$$

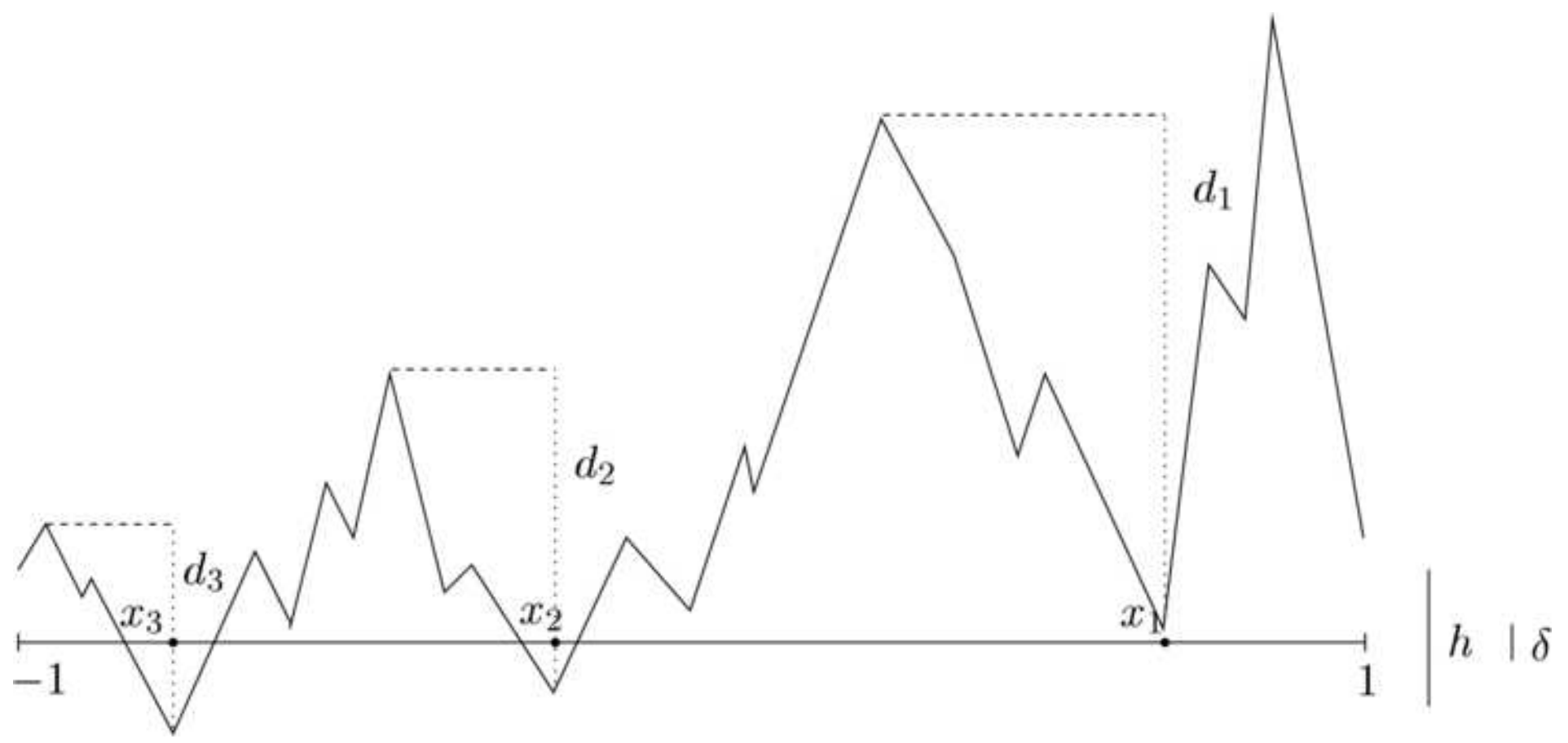

FIG. 2. *Path in* $\mathcal{A}_{h,\delta}$.

1.4. *Main results.* We can finally state our main results concerning the spectral analysis of the operators $\mathbb{L}(\{-N+1,\ldots,N-2,N-1\})$, $N\geq 1$. Note that since the rescaled potential $V^{(N)}$ defined in (1.9) converges weakly to the two-sided Brownian motion $B$, it is more natural to work on the rescaled lattice $\mathbb{Z}/N$. In particular, we introduce the infinite matrix $\mathcal{L}^{(N)}$ with entries

$$\mathcal{L}^{(N)}_{x,y}\equiv \mathbb{L}_{Nx,Ny} \qquad \forall x,y\in\mathbb{Z}/N,$$

and, for each $D\subset\mathbb{Z}/N$, we denote by $\mathcal{L}^{(N)}(D)$ the restriction of $\mathcal{L}^{(N)}$ to $D\times D$, that is,

$$\mathcal{L}^{(N)}(D)\equiv(\mathcal{L}^{(N)}_{x,y})_{x,y\in D}. \tag{1.19}$$

Defining

$$I_N\equiv(-1,1)\cap(\mathbb{Z}/N) \qquad \forall N\geq 1, \tag{1.20}$$

$u$ is an eigenvector of $\mathbb{L}(\{-N+1,\ldots,N-2,N-1\})$ with eigenvalue $\lambda$ iff $u(N\cdot)$ is an eigenvector of $\mathcal{L}^{(N)}(I_N)$ with eigenvalue $\lambda$. Note that the operator $\mathcal{L}^{(N)}(I_N)$ is symmetric on $\mathbb{L}^2(I_N,\mu_N)$, where

$$\mu_N(x)\equiv\mu(Nx) \qquad \forall x\in\mathbb{Z}/N. \tag{1.21}$$

In what follows, we denote by $V_N$ the restriction of the rescaled potential $V^{(N)}$ [see (1.9)] to $[-1,1]$, that is,

$$V_N:[-1,1]\to\mathbb{R}, \qquad V_N(t)\equiv V^{(N)}(t) \qquad \forall t\in[-1,1]. \tag{1.22}$$

Given a good path $\gamma \in \mathcal{A}_{h,\delta}$ such that $M_h^-(\gamma) = \{x_1, x_2, \dots, x_q\}$ is the special labeling satisfying condition (1.16), we use the following notation:

$$\begin{cases} M_h^-(\gamma) = \{x_1, x_2, \dots, x_q\}, & \\ M_{h,k}^-(\gamma) \equiv \varnothing, & \text{if } k = 0, \\ M_{h,k}^-(\gamma) \equiv \{x_1, \dots x_k\}, & \text{if } 1 \le k \le q, \\ S_{h,k}(\gamma) \equiv M_{h,k}^-(\gamma) \cup \{-1, 1\}, & \text{if } 0 \le k \le q, \\ S_{h,k}^*(\gamma) \equiv S_{h,k}(\gamma) \cup (\mathbb{Z}/N \setminus (-1,1)), & \text{if } 0 \le k \le q, \\ d_k \equiv \gamma(z^*(x_k, S_{h,k-1})) - \gamma(x_k), & \text{if } 1 \le k \le q. \end{cases} \tag{1.23}$$

We omit $\gamma$ from the above notation when the path is understood.

Finally, we write $P_{x,N}^\omega$ for the law of the rescaled Sinai's random walk $(X_n/N, n \ge 0)$ starting in $x \in \mathbb{Z}/N$, with environment $\omega$, and we set

$$\tau_A = \min\{n \ge 1 : X_n/N \in A\}, \qquad A \subset \mathbb{Z}/N.$$

THEOREM 1. *Given $Q, h, \delta > 0$ if $V_N \in \mathcal{A}_{h,\delta}$, the number of $h$-minima $q \equiv |M_h^-| \le Q$ and the rescaling parameter $N \ge N(\delta, Q)$, then the following holds:*

*The rescaled generator $\mathcal{L}^{(N)}(I_N)$ with Dirichlet conditions outside $(-1,1)$ has exactly $q$ eigenvalues smaller than the principal eigenvalue $\lambda_N^*$ of the operator with Dirichlet conditions outside $(-1,1)$ and on the set of $h$-minima $M_h^-$, that is,*

$$\begin{cases} \sigma(\mathcal{L}^{(N)}(I_N)) \cap [0, \lambda_N^*) = \{\lambda_1^{(N)} < \lambda_2^{(N)} < \dots < \lambda_q^{(N)}\}, \\ \lambda_N^* \equiv \min \sigma(\mathcal{L}^{(N)}(I_N \setminus M_h^-)). \end{cases} \tag{1.24}$$

*Moreover, the threshold $\lambda_N^*$ satisfies*

$$\lambda_N^* \ge N^{-2} e^{-h\sqrt{N}}. \tag{1.25}$$

*Setting*

$$h_k(x) = \begin{cases} P_{N,x}^\omega(\tau_{x_k} < \tau_{S_{h,k-1}}), & \text{if } x \notin M_h^-, \\ 1, & \text{if } x = x_k, \\ 0, & \text{if } x \in S_{h,k-1}, \end{cases} \tag{1.26}$$

*and denoting $\|\cdot\|_2$ the norm in $L^2(I_N, \mu_N)$, the first $q$ eigenvalues $\lambda_k^{(N)}$ admit the probabilistic approximation*

$$\lambda_k^{(N)} = \frac{\mu_N(x_k) P_{N,x_k}^\omega(\tau_{S_{h,k-1}} < \tau_{x_k})}{\|h_k\|_2^2} (1 + O(e^{(-\delta/10)\sqrt{N}})) \tag{1.27}$$

*and satisfy*

$$c(\kappa) N^{-2} e^{-\sqrt{N} d_k} \le \lambda_k^{(N)} \le c'(\kappa) e^{-\sqrt{N} d_k}. \tag{1.28}$$

*Moreover, for* $1 \leq k \leq q$*, there exists a normalized eigenvector* $\psi_k^{(N)}$ *with eigenvalue* $\lambda_k^{(N)}$ *such that*

$$\left\| \psi_k^{(N)} - \frac{h_k}{\|h_k\|_2} \right\|_2 \leq e^{(-\delta/10)\sqrt{N}}. \tag{1.29}$$

Proof. The theorem follows easily from Lemma 7, Proposition 4 and Theorem 6. □

Due to (1.18) and (1.28), under the conditions of Theorem 1, the first $q$ eigenvalues of $\mathcal{L}^{(N)}(I_N)$ split as follows:

$$\lambda_k^{(N)} \leq c(\kappa) N^2 e^{-\delta\sqrt{N}} \lambda_{k+1}^{(N)} \qquad \forall k = 1, \ldots, q-1. \tag{1.30}$$

Finally, we observe that the hypothesis of Theorem 1 is satisfied with probability tending to one, as $N$ tends to infinity.

Theorem 2. *For any* $\alpha > 0$*, there exist* $h > 0, \delta > 0$*, and* $Q < \infty$*, such that*

$$\liminf_{N\uparrow\infty} \mathbf{P}(\mathcal{A}_{h,\delta} \cap \{|M_h^-(V_N)| \leq Q\}) \geq 1 - \alpha. \tag{1.31}$$

This theorem, and its extended version given by Theorem 5 of Section 2, follows from a result of Neveu and Pitman (Proposition 1) about the $h$-extrema of Brownian motion and on the KMT approximation theorem [(1.10) Proposition 3]. The proof of Theorem 5 is given in Section 2.

Remark 1. The above theorems reproduce (and partly refine) results obtained in [15] via a (nonrigorous) renormalization group (RG). In Appendix A we show that the labeling of $h$-minima satisfying condition (1.16) is equivalent to the labeling obtained in [15] via the RG.

Remark 2. Theorem 1 is very similar in nature to the results in [3] and [5] on metastable Markov chains, respectively reversible diffusions in smooth potentials and our proofs will follow the strategy outlined in these papers. The purpose of Theorem 1 is to provide a precise relation between spectral properties of the generator and geometric properties of the random potential $V_N$, in particular, its $h$-extrema.

The hypothesis of Theorem 1 provide, the analogue of the nondegeneracy conditions required, for example, in [5]. The validity of these hypothesis, as asserted by Theorem 2, as well as information on the statistical properties of the eigenvalues can be derived from the statistical properties on the $h$–extrema of $V_N$. Due to the KMT approximation theorem, the rescaled potential $V_N$ can be thought of as a $L^\infty$-perturbation of the Brownian motion.

Hence, as described in Theorem 5, the $h$-extrema of $V_N$ are well approximated by the $h$-extrema of the Brownian motion, whose statistics are provided by the Neveu–Pitman results. A detailed analysis is given is Section 2.

REMARK 3. Our spectral analysis is based on potential theory, briefly discussed in Section 3. Both the formulas (1.27) and (1.29) can be restated in terms of capacity and equilibrium potential. In fact, the equilibrium potential $h_{x_k,S^*_{h,k-1}}$ associated to $x_k, S^*_{h,k-1}$ coincides with $h_k$, while the capacity $\mathrm{cap}(x_k, S^*_{h,k-1})$ between $x_k$ and $S_{h,k-1}$ satisfies

$$\mathrm{cap}(x_k, S^*_{h,k-1}) = \mu_N(x_k) P^{\omega}_{N,x_k}(\tau_{S_{h,k-1}} < \tau_{x_k}).$$

Hence, (1.27) corresponds to the formula

$$\lambda_k^{(N)} = \frac{\mathrm{cap}(x_k, S^*_{h,k-1})}{\|h_{x_k,S^*_{h,k-1}}\|_2^2}(1 + O(e^{(-\delta)10\sqrt{N}})).$$

A more detailed description of the eigenvector $\psi_k^{(N)}$ via potential theory is given in Theorem 6. As explained in Section 3, as we are in dimension one, both the equilibrium potential and the capacity admit simple expressions that, together with the results of Section 2, allow to get from Theorems 1 and 6 rather precise quantitative estimates on the eigenvalues and eigenfunctions.

We briefly describe the strategy leading to the spectral analysis of $\mathcal{L}^{(N)}(I_N)$ for small eigenvalues, referring to Sections 5 and 6 for more details. We suppose $V_N \in \mathcal{A}_{h,\delta}$ and define

$$\mathcal{L}_k^{(N)} \equiv \mathcal{L}^{(N)}(I_N \setminus M^-_{h,k}), \qquad 0 \le k \le q-1, \tag{1.32}$$

that is, $\mathcal{L}_k^{(N)}$ is the rescaled generator with Dirichlet conditions on $M^-_{h,k} = \{x_1, x_2, \ldots, x_k\}$ and outside $(-1,1)$ (see Figure 3). We call $\bar{\lambda}_k^{(N)}$ the principal eigenvalue of $\mathcal{L}_k^{(N)}$:

$$\bar{\lambda}_k^{(N)} \equiv \min \sigma(\mathcal{L}_k^{(N)}), \qquad 0 \le k \le q-1. \tag{1.33}$$

Note that

$$\mathcal{L}^{(N)}(I_N) = \mathcal{L}_0^{(N)}, \qquad \lambda_1^{(N)} = \bar{\lambda}_0^{(N)}. \tag{1.34}$$

The analysis of the above principal eigenvalues $\bar{\lambda}_k^{(N)}$ and the associated principal eigenvectors is given in Section 5 (see Propositions 4 and 5) and is based on potential theory. In Section 6, by a perturbation argument, we prove for $0 \le k < q$ that the eigenvalue $\lambda_{k+1}^{(N)}$ and the eigenvector $\psi_{k+1}^{(N)}$ are well approximated respectively by the principal eigenvalue $\bar{\lambda}_k^{(N)}$ and the principal eigenvector of $\mathcal{L}_k^{(N)}$, thus, leading to Theorem 1 and Theorem 6.

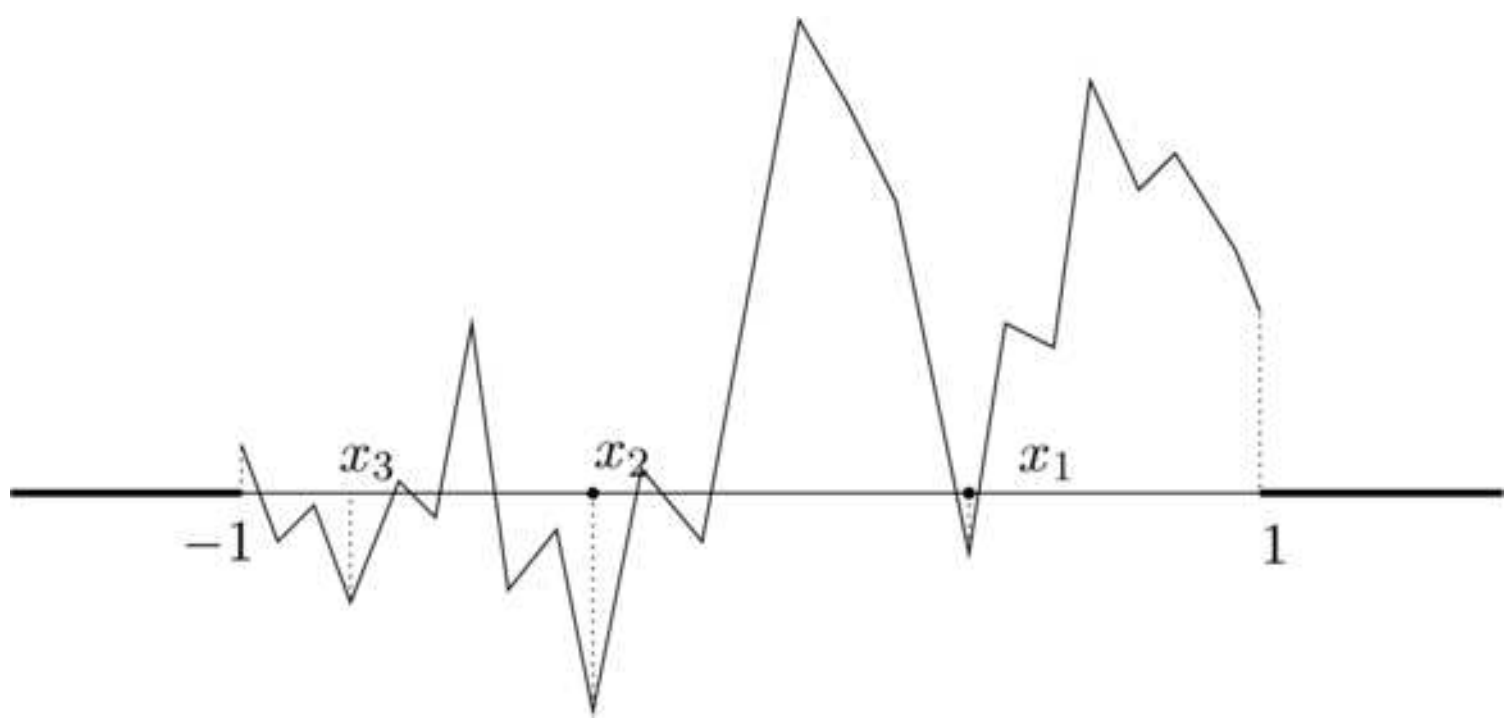

FIG. 3. *$\mathcal{L}_2^{(N)}$ has Dirichlet conditions on the marked regions.*

As application of our spectral analysis, we will give a spectral proof of a refined version of Sinai's theorem:

THEOREM 3. *Recall the definition of $V^{(1)}$ given in (1.9). For each $n \in \mathbb{Z}_+$ and $\omega \in \Omega$, let $m^{(n)}(\omega) \in M^-_{\ln n}(V^{(1)})$ be the $\ln n$-minimum corresponding to the bottom of the $\ln n$-valley covering the origin and set $\mathfrak{m}^{(n)} \equiv m^{(n)}(\omega)/\ln^2 n$. Fix $\alpha > 0$ and a positive function $\rho$ on $(0,\infty)$ such that*

$$\lim_{x\downarrow 0} x^2/\rho(x) = 0. \tag{1.35}$$

*Then, for each $n$, there exists a Borel subset $\Omega_n \subset \Omega$ with $\mathbf{P}(\Omega_n) \geq 1-\alpha$ and*

$$\lim_{n\uparrow\infty} \inf_{\omega\in\Omega_n} P_0^\omega\left(\left|\frac{X_n}{\ln^2 n} - \mathfrak{m}^{(n)}(\omega)\right| \leq \delta_n\right) = 1, \qquad \delta_n \equiv \rho\left(\frac{\ln(\ln n)}{\ln n}\right). \tag{1.36}$$

REMARK 4. We point out that, due to Golosov's localization result (1.6), the limit in (1.36) must hold with $\delta_n \equiv \rho(1/\ln n)$. Hence, Theorem 3 is not optimal. The gap is not due to the spectral method and could be filled as follows. The proof of Theorem 3 (given in Section 7) needs some knowledge on the local behavior of $V^{(\ln n)}$ around $\mathfrak{m}^{(n)}$ [equivalently, on the local behavior of $V^{(1)}$ around $m^{(n)}$] in order to prove (7.30). As already discussed in Remark 2, we study the geometric properties of the rescaled potential by comparing it with the Brownian motion via the KMT Approximation Theorem. While this method provides good information about the global statistics of the $h$-extrema of $V^{(\ln n)}$ in a given box, it gives a very rough picture of the local behavior of $V^{(\ln n)}$ around $\mathfrak{m}^{(n)}$. This lack of information is paid by the factor $\ln\ln n$ in the definition of $\delta_n$ in (1.36). In order to avoid it, a direct analysis of $V^{(\ln n)}$ near to $\mathfrak{m}^{(n)}$ as in [12] is necessary, and also sufficient.

From the spectral information we can derive easily another characterization of the long term dynamics. Let $A_n \subset \mathbb{Z}$ be the box covered by the $\ln n$-valley of $V^{(1)}$ containing the origin. Set $D_n = (m^{(n)} - \delta_n \ln^2 n, m^{(n)} - \delta_n \ln^2 n)$, where $m^{(n)}, \delta_n$ are as in Theorem 3. Then construct a sequence of boxes $A_{n_k}$ as follows: start with $A_{n_0}$, $n_0$ large. Then increase $n$ to $n_1$ such that, for the first time, $m^{(n_1)} \neq m^{(n_0)}$, and so on. Finally, define $\lambda_k$ as the second (in increasing order) eigenvalue of $\mathbb{L}(A_{n_k})$. Then the following holds:

THEOREM 4. *Fix $\alpha > 0$. Then for each $k$, there exists a subset $\Omega_{n_k} \subset \Omega$ with $\mathbf{P}(\Omega_{n_k}) \geq 1 - \alpha$ and*

$$\lim_{k\uparrow\infty} \sup_{\omega\in\Omega_{n_k}} |P_0^\omega(X_{t/\lambda_k} \in D_{n_k}) - (1 - e^{-t})| = 0. \tag{1.37}$$

Theorem 4 throws a somewhat nonageing like view on Sinai's model. It says that there is an infinity of diverging (and well separated) time-scales on which the process looks as if it would approach equilibrium exponentially.

To see ageing effect, one needs to go into a different regime of time scales. In fact, Dembo, Guionnet and Zeitouni [8] (see also [23]) have shown that

$$\lim_{n\uparrow\infty} \mathbf{P}_0(X_n \sim X_{n^h}) = h^{-2}(\tfrac{5}{3} + \tfrac{2}{3}e^{-(h-1)}),$$

that is, ageing occurs on an exponential time scale. Note that this result follows easily from Theorem 3 and the right-hand side is just the probability that $\mathfrak{m}^n = \mathfrak{m}^{n^h}$, as observed in [8].

We divide the remainder of this paper as follows. In Section 2 we recall a theorem of Neveu and Pitman [17] about the statistics of $h$-extrema of Brownian motion and use it to derive the required statistical properties of the random potentials. In particular, we prove Theorem 2 and its extension Theorem 5. In Section 3 we recall some elementary background from potential theory for later use. In Section 4 we compute hitting times and conditional hitting times of our process. In Section 5 we compute principle eigenvalues and eigenvectors for $\mathcal{L}_k^{(N)}$ defined in (1.32). In Section 6 we prove Theorem 6. In Section 7 we use the spectral results to prove Theorem 3 and Theorem 4.

**2. *h*-extrema of Brownian motion and random walks.** The following result about the statistics of $h$-extrema for Brownian motion is due to Neveu and Pitman [17]. We state it here for the Brownian motion $B = (B_t, t \in \mathbb{R})$ with variance $\sigma^2$.

PROPOSITION 1 (Neveu and Pitman [17]). *The set of $h$-extrema $\mathcal{E}_h(B)$ for the Brownian motion $B=(B_t, t\in\mathbb{R})$ is a stationary renewal process. Setting $\mathcal{E}_h(B)=\{S_n^{(h)}\}_{n\in\mathbb{Z}}$, with*

$$\cdots < S_{-1}^{(h)} < S_0^{(h)} \le 0 < S_1^{(h)} < S_2^{(h)} < \cdots,$$

*then the trajectories between $h$-extrema (called $h$-*slopes*)*

$$(B_{S_n^{(h)}+t} - B_{S_n^{(h)}} : 0 \le t \le S_{n+1}^{(h)} - S_n^{(h)}) \tag{2.1}$$

*are independent and, for $n \neq 0$, identically distributed, up to changes of sign. In particular, the variables*

$$|B_{S_{n+1}^{(h)}} - B_{S_n^{(h)}}| - h, \qquad n\in\mathbb{Z},$$

*are independent and exponentially distributed with mean $h$, whereas the variables $S_{n+1}^{(h)} - S_n^{(h)}$, $n\neq 0$, are i.i.d, with Laplace transform*

$$\mathbf{E}_B(\exp\{-\lambda(S_{n+1}^{(h)} - S_n^{(h)})\}) = 1/\cosh\left(\frac{h\sqrt{2\lambda}}{\sigma}\right) \tag{2.2}$$

*and mean $h^2/\sigma^2$.*

[Note that in [17] the r.h.s. of (2.2) is written with $\sqrt{2\lambda}$ replaced by $\sqrt{2}\lambda$. As explained in [6], Section 2, the correct form if given by (2.2).]

Note that $\mathcal{M}_h^-(\gamma)=M_h^-(\gamma)$ for $\mathbf{P}_B$ almost all $\gamma$ and that, since $(B_t, t\in\mathbb{R}) \overset{\text{law}}{=} (B_{ta^2}/a, t\in\mathbb{R})$ for all $a>0$,

$$(S_n^{(h)}, n\in\mathbb{Z}) \overset{\text{law}}{=} (a^2 S_n^{(h/a)}, n\in\mathbb{Z}) \qquad \forall a>0. \tag{2.3}$$

As in [6], in order to describe the law of the trajectory (2.1) for $n\neq 0$, it is convenient to introduce the Polish space, $\mathbf{G}$, of continuous paths, $\gamma:[0,\ell(\gamma)]\to\mathbb{R}$, defined on some interval $[0,\ell(\gamma)]$, equipped with the metric

$$d(\gamma,\gamma') \equiv |\ell(\gamma)-\ell(\gamma')| + \max_{t\in[0,1]} |\gamma(t\ell(\gamma)) - \gamma'(t\ell(\gamma))|.$$

In the sequel, we will consider random paths as $\mathbf{G}$-valued random variables.

Starting from the Brownian motion $B$ we define (see also Figure 4)

$$\begin{aligned} S_t &\equiv \max\{B_s : 0\le s\le t\},\\ \tau &\equiv \min\{t>0 : S_t = B_t + h\},\\ \beta &\equiv S_\tau,\\ \alpha &\equiv \max\{t : 0\le t\le\tau \text{ and } B_t=\beta\}. \end{aligned}$$

Note that $\mathbf{P}_B$ a.s. there exists a unique $s\in[0,\tau]$ such that $B_s=\beta$. As proved in [17], the random paths $(B_t : 0\le t\le\alpha)$ and $(\beta - B_t : 0\le t\le \tau-\alpha)$ are independent.

Moreover, in [17] the following result is proved.

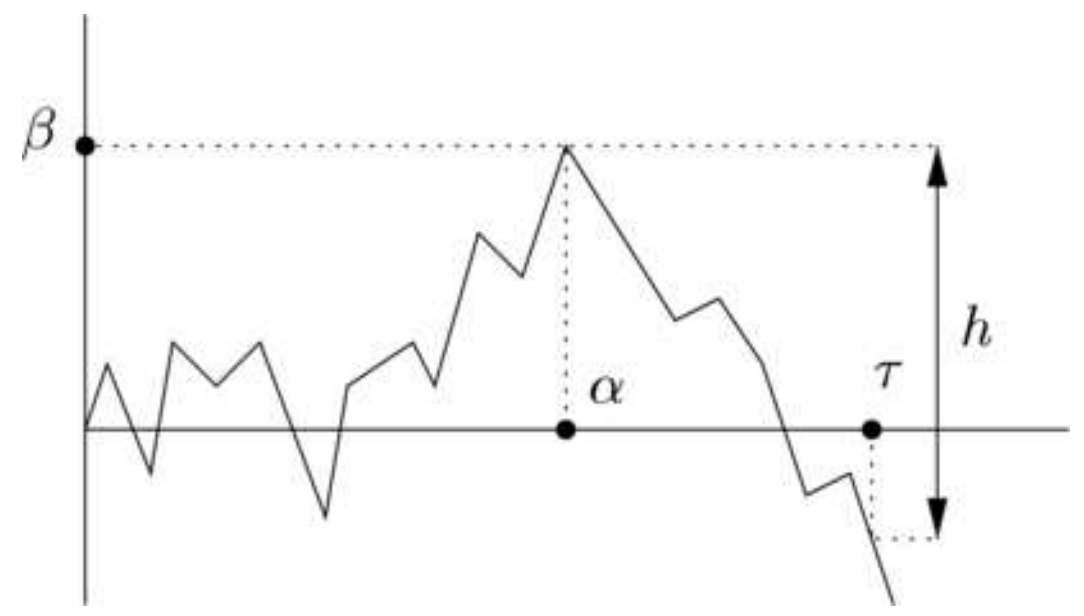

FIG. 4. *Definition of $\beta, \alpha, \tau$.*

PROPOSITION 2 (Neveu and Pitman [17]). *For $n \neq 0$, the random path*

$$(|B_{S_n^{(h)}+t} - B_{S_n^{(h)}}| : 0 \leq t \leq S_{n+1}^{(h)} - S_n^{(h)})$$

*has on* $\mathbf{G}$ *the same law of the random path $\gamma_B : [0, \tau] \to \mathbb{R}_+$ defined as*

$$\gamma_B(t) \equiv \begin{cases} \beta - B_{\alpha+t}, & \text{if } t \in [0, \tau - \alpha], \\ h + B_{t-(\tau-\alpha)}, & \text{if } t \in [\tau - \alpha, \tau). \end{cases}$$

The above proposition, Corollary 4.4 in [18], Chapter XII and the reflection invariance of Brownian motion easily imply the following result concerning the behavior of the Brownian motion near to an $h$-extremum.

COROLLARY 1. *Given $n \in \mathbb{Z}$, let*

$$T_{n,+}^{(h)} \equiv \min\{t \in (0, S_{n+1}^{(h)} - S_n^{(h)}) : |B_{S_n^{(h)}+t} - B_{S_n^{(h)}}| = h\}, \tag{2.4}$$

$$T_{n,-}^{(h)} \equiv \max\{t \in (0, S_n^{(h)} - S_{n-1}^{(h)}) : |B_{S_n^{(h)}} - B_{S_n^{(h)}-t}| = h\}. \tag{2.5}$$

*Moreover, let $Z = BES^3(0)$, namely, $Z = (Z_t, \geq 0)$ is a Bessel process of dimension 3 starting at the origin, independent of the Brownian motion $B$. Let $T_h$ be the hitting time*

$$T_h = \min\{t > 0 : \sigma Z_t = h\}.$$

*Then the random paths*

$$(|B_{S_n^{(h)}+t} - B_{S_n^{(h)}}|,\ 0 \leq t \leq T_{n,+}^{(h)}) \qquad \text{with } n \neq 0,$$

$$(|B_{S_n^{(h)}} - B_{S_n^{(h)}-t}|, 0 \leq t \leq T_{n,-}^{(h)}) \qquad \text{with } n \neq 1, \text{ and } (\sigma Z_t, 0 \leq t \leq T_h)$$

*have the same law on* $\mathbf{G}$.

LEMMA 3. *Let $Z = BES^3(0)$, then*

$$P\Big(\inf_{s \geq t} Z_s < \varepsilon\Big) < \sqrt{2}\varepsilon/\sqrt{\pi t} \qquad \forall \varepsilon, t > 0.$$

PROOF. Let us define $J_t \equiv \inf_{s\geq t} Z_s$. By Pitman theorem (see Theorem 3.5 in [18], Chapter VI), $J_t$ has the same law of $S_t$. In particular, by the reflection principle for Brownian motion,

$$P(J_t<\varepsilon)=P(S_t<\varepsilon)=P(|B_t|<\varepsilon)=P(|B_1|<\varepsilon/\sqrt{t})<\sqrt{2}\varepsilon/\sqrt{\pi t}. \quad \square$$

REMARK 5. Using renewal theory, one can describe the law on $\mathbf{G}$ of the random path

$$(|B_{S_0^{(h)}+t}-B_{S_0^{(h)}}|: 0\leq t\leq S_1^{(h)}-S_0^{(h)}).$$

In [6] it is shown that, conditioning on the length $S_{n+1}^{(h)}-S_n^{(h)}$, the law of the path $(|B_{S_n^{(h)}+t}-B_{S_n^{(h)}}|: 0\leq t\leq S_{n+1}^{(h)}-S_n^{(h)})$ does not depend on $n$, for all $n\in\mathbb{Z}$. Moreover, the random variable $S_{n+1}^{(h)}-S_n^{(h)}$ has respectively probability density $(\sigma/h)^2 f(x(\sigma/h)^2)\,dx$ and $x(\sigma/h)^4 f(x(\sigma/h)^2)\,dx$ if $n\neq 0$ and $n=0$, where

$$f(x)=\mathbb{I}_{x>0}\frac{\pi}{2}\sum_{k\in\mathbb{Z}}(-1)^k\left(k+\frac{1}{2}\right)\exp\left\{-\frac{\pi^2}{2}\left(k+\frac{1}{2}\right)^2 x\right\}. \tag{2.6}$$

The above result will not be used in what follows, while we will need some information about the distribution of $S_1^{(h)}$. This can be obtained from renewal theory as follows. Calling $F$ and $G$ respectively the distribution functions of $S_{n+1}^{(h)}-S_n^{(h)}$ for $n\in\mathbb{Z}\setminus\{0\}$ and $S_1^{(h)}$, formula (4.7) in [9], Chapter 3, reads

$$G(t)=\int_0^t(1-F(y))\,dy\Big/\int_0^t y\,dF(y).$$

Due to the above identity and integration by parts, one obtains

$$\mathbf{P}_B(S_1^{(h)}\leq t)=\frac{tP(X^{(h)}>t)+E(X^{(h)};X^{(h)}\leq t)}{E(X^{(h)})}, \tag{2.7}$$

where $X^{(h)}$ is a random variable with Laplace transform

$$E(\exp\{-\lambda X^{(h)}\})=1/\cosh(h\sqrt{2\lambda}/\sigma). \tag{2.8}$$

We conclude this section with a technical lemma whose proof is given in Appendix B. Given a path $\gamma\in C(\mathbb{R})$, we define $\gamma^*\equiv\{\gamma_t\}_{|t|\leq 1}\in C([-1,1])$. Note that $\mathbf{P}_B$ a.s. $M_h^-(\gamma^*)\subset M_h^-(\gamma)\cap[-1,1]$ and $|M_h^-(\gamma)\cap[-1,1]|-|M_h^-(\gamma^*)|\leq 2$.

LEMMA 4. *Let $h, H, \beta, \delta, \varepsilon$ be positive constants with $h < H$. Recall the definition of $\mathcal{A}_{h,\delta}$ given in Section 1.3 and define the events $\mathcal{B}_{h,\delta}, \mathcal{C}_{h,\delta}$, $\mathcal{D}_{h,\beta,\varepsilon}$ as*

$$\mathcal{B}_{h,\delta} \equiv \{\gamma : \exists n \in \mathbb{Z} \text{ s.t. } |\gamma(S_{n+1}^{(h)}) - \gamma(S_n^{(h)})| < h + \delta \text{ and } S_n^{(h)}, S_{n+1}^{(h)} \in [-1,1]\},$$

$$\mathcal{C}_{h,\delta} \equiv \{\gamma : \exists (x,y) \neq (x',y') \in M_h^-(\gamma) \times M_h^+(\gamma) \cap [-1,1]^2, \quad \text{s.t. } ||\gamma(x) - \gamma(y)| - |\gamma(x') - \gamma(y')|| < \delta\}$$

*and*

$$\mathcal{D}_{h,\beta,\varepsilon} \equiv \left\{\exists n \in \mathbb{Z} : S_n^{(h)} \in [-1,1] \text{ and } \inf_{t \in [-T_{n,-}^{(h)}, T_{n,+}^{(h)}] \setminus [-\beta,\beta]} |B_{S_n^{(h)}+t} - B_{S_n^{(h)}}| < \varepsilon\right\}.$$

*Then,*

$$\mathbf{P}_B(\gamma : |\mathcal{E}_h(\gamma) \cap [-1,1]| \geq 4) \geq 1 - c(\alpha,\sigma)h^{\alpha} \qquad \forall \alpha > 0, \tag{2.9}$$

$$\mathbf{P}_B(\gamma : |\mathcal{E}_h(\gamma) \cap [-1,1]| \geq n) \leq e\left(1 + \frac{h^2}{2\sigma^2}\right)^{-n}, \tag{2.10}$$

$$\mathbf{P}_B(\mathcal{B}_{h,\delta}) \leq c(H,\sigma)(1 - e^{-\delta/h})h^{-4}, \tag{2.11}$$

$$\mathbf{P}_B(\mathcal{C}_{h,\delta}) \leq c(H,\sigma)\delta h^{-11}, \tag{2.12}$$

$$\mathbf{P}_B(\mathcal{D}_{h,\beta,\varepsilon}) \leq c(\sigma)\left(1 + \frac{\varepsilon^{1/4}}{\beta^{1/8}}\right)\frac{\varepsilon^{1/2}}{\beta^{1/4}}, \tag{2.13}$$

*where $c(\sigma)$, $c(\alpha,\sigma)$, $c(H,\sigma)$ are suitable positive constants depending respectively on $\sigma$, $\alpha,\sigma$ and $H,\sigma$. In particular, for each $\alpha > 0$, there exist positive constants $h, \delta, Q$ such that*

$$\mathbf{P}_B(\gamma : \gamma^* \in \mathcal{A}_{h,\delta} \text{ and } |\mathcal{M}_h^-(\gamma^*)| \leq Q) \geq 1 - \alpha, \tag{2.14}$$

*where $\gamma^* \equiv \{\gamma_t\}_{|t| \leq 1}$.*

2.1. *$L^\infty$-perturbations of Brownian motion.* So far, we have collected properties of Brownian motion. We want to use the KMT strong approximation result to deduce analogous results for the rescaled random walks $V_N$ defined in (1.22).

PROPOSITION 3. *For suitable positive constants $C_1, C_2, C_3$ depending on $\sigma$, given $N \in \mathbb{Z}_+$, there exists a coupling on an enlarged probability space between $V^{(N)}$ and the two-sided Brownian motion $B$ with variance $\sigma$ such that*

$$P^{(N)}\left(\sup_{x \in [-1,1]} |V_N(x) - B_x| > \frac{C_1 \ln N}{\sqrt{N}}\right) < \frac{C_2}{N^{C_3}}. \tag{2.15}$$

Proof. This follows easily from the Komlós–Major–Tusnády strong approximation theorem [14] and some elementary regularity estimates controlling the variation of Brownian motion between lattice points of $\mathbb{Z}/N$. $\square$

The following lemma describes the effect of a $L^\infty$-perturbation on the location of $h$-extrema.

Lemma 5. *Let $h, \varepsilon > 0$ and let $\gamma, \gamma' \in C([-1,1])$ such that*

$$M^-_{h+\varepsilon}(\gamma) = M^-_{h-\varepsilon}(\gamma), \qquad M^+_{h+\varepsilon}(\gamma) = M^+_{h-\varepsilon}(\gamma), \qquad \|\gamma - \gamma'\|_\infty \le \frac{\varepsilon}{4}.$$

*Let $M^-_h(\gamma) = \{u_1, u_2, \dots, u_q\}$, $M^+_h(\gamma) = \{w_1, w_2, \dots, w_q, w_{q+1}\}$, where*

$$-1 \le w_1 < u_1 < w_2 < \cdots < u_q < w_{q+1} \le 1,$$

*and set*

$$u'_i \equiv \min\Bigl\{z \in [w_i, w_{i+1}] : \gamma'(z) = \min_{[w_i, w_{i+1}]} \gamma'\Bigr\} \qquad \forall i = 1, \dots, q,$$

$$w'_i \equiv \min\Bigl\{z \in [u_{i-1}, u_i] : \gamma'(z) = \max_{[u_{i-1}, u_i]} \gamma'\Bigr\} \qquad \forall i = 1, \dots, q+1,$$

*where $u_0 \equiv -1, u_{q+1} \equiv 1$. Then*

$$M^-_h(\gamma') = \{u'_1, u'_2, \dots, u'_q\}, \qquad M^+_h(\gamma') = \{w'_1, w'_2, \dots, w'_{q+1}\}.$$

*Moreover,*

$$\begin{gathered} u'_i \in \{x \in [w_i, w_{i+1}] : \gamma(x) \le \gamma(u_i) + \varepsilon/2\}, \\ |\gamma'(u'_i) - \gamma(u_i)| \le \frac{\varepsilon}{4}, \forall i = 1, \dots, q, \end{gathered} \tag{2.16}$$

$$\begin{gathered} w'_i \in \{x \in [u_{i-1}, u_i] : \gamma(x) \ge \gamma(w_i) - \varepsilon/2\}, \\ |\gamma'(w'_i) - \gamma(w_i)| \le \frac{\varepsilon}{4}, \forall i = 1, \dots, q+1. \end{gathered} \tag{2.17}$$

Proof. We leave the simple case $q = 0$ to the reader and assume here that $q \ge 1$. Due to Lemma 2,

$$\gamma'(w_{i+1}) - \gamma'(u'_i) \ge \gamma'(w_{i+1}) - \gamma'(u_i) \ge h + \varepsilon/2 \qquad \forall i = 1, \dots, q$$

and, similarly, $\gamma'(w_i) - \gamma'(u'_i) \ge h + \varepsilon/2$. This, together with the definition of $u'_i$, implies that $u'_i \in M^-_h(\gamma')$. Let us suppose that $|M^-_h(\gamma')| > |M^-_h(\gamma)|$. Then at least one of the following cases must hold:

(C1) $\exists u \in [-1, w_1]$ such that $u \in M^-_h(\gamma')$,
(C2) $\exists u \in [w_{q+1}, 1]$ such that $u \in M^-_h(\gamma')$,

(C3) for some $i \in \{1, \dots, q\}$, there exists $u \in [w_i, w_{i+1}] \setminus \{u'_i\}$ such that $u \in M_h^-(\gamma')$.

We treat the case (C3). Let us suppose, for example, that $u < u'_i$. Then $\exists y : u < y < u'_i$ such that $\gamma'(y) - \gamma'(u) > h$ and $\gamma'(y) - \gamma'(u'_i) > h$. But this would imply that $\gamma(y) - \gamma(u) > h - \varepsilon/2$ and $\gamma(y) - \gamma(u'_i) > h - \varepsilon/2$ and therefore that $M_{h-\varepsilon}^-(\gamma) \cap [w_i, w_{i+1}]$ has at least two elements in contradiction with the hypothesis that $M_{h-\varepsilon}^-(\gamma) = M_h^-(\gamma)$. Hence, (C3) cannot hold. By similar arguments, one can prove that both (C1) and (C2) cannot be valid. This completes the proof that $M_h^-(\gamma') = \{u'_1, u'_2, \dots, u'_q\}$. The proof that $M_h^+(\gamma') = \{w'_1, w'_2, \dots, w'_{q+1}\}$ is similar. In order to prove the first assertion of (2.16), we need to show that $\gamma(u'_i) \le \gamma(u_i) + \varepsilon/2$. To this end, it is enough to observe that

$$\gamma(u'_i) \le \gamma'(u'_i) + \varepsilon/4 \le \gamma'(u_i) + \varepsilon/4 \le \gamma(u_i) + \varepsilon/2,$$

where the second inequality comes from the definition of $u'_i$. Note that similarly one gets $\gamma'(u'_i) \ge \gamma(u_i) - \varepsilon/4$, thus completing the proof of (2.16). The proof of (2.17) is similar. $\square$

THEOREM 5. *Let $h, \delta > 0$ and let $P^{(N)}$, $C_1, C_2, C_3$ be as in Proposition 3. Set*

$$\varepsilon = \varepsilon(N) \equiv 4 \frac{C_1 \ln N}{\sqrt{N}}$$

*and fix a function $\beta : \mathbb{Z}_+ \to (0, \infty)$ such that $\lim_{N \uparrow \infty} \varepsilon(N)/\sqrt{\beta(N)} = 0$.*

*Let $B^* \equiv (B_t, t \in [-1, 1])$. On the enlarged probability space with probability measure $P^{(N)}$, let $\mathcal{G}_{h,\delta}$ be the event that the following conditions are fulfilled:*

(i)

$$\sup_{x \in [-1,1]} |V_N(x) - B_x| \le \frac{\varepsilon}{4} = \frac{C_1 \ln N}{\sqrt{N}},$$

(ii)

$$|M_h^-(V_N)| = |M_h^-(B^*)|, \qquad |M_h^+(V_N)| = |M_h^+(B^*)|,$$

(iii)

$$\begin{cases} |V_N(u'_i) - B^*(u_i)| \le \varepsilon/4 \text{ and } |u'_i - u_i| \le \beta(N), & \forall 1 \le i \le q, \\ |V_N(w'_i) - B^*(w_i)| \le \varepsilon/4 \text{ and } |w'_i - w_i| \le \beta(N), & \forall 2 \le i \le q, \\ |V_N(w'_i) - B^*(w_i)| \le \varepsilon/4 & \\ \qquad \text{for } i = 1, q+1, w'_1 \in [-1, u_1), w'_{q+1} \in (u_q, 1], & \end{cases}$$

*where*

$$\begin{cases} M_h^-(V_N) = \{u_1' < u_2' < \cdots < u_q'\}, \\ M_h^+(V_N) = \{w_1' < w_2' < \cdots < w_q' < w_{q+1}'\}, \\ M_h^-(B^*) = \{u_1 < u_2 < \cdots < u_q\}, \\ M_h^+(B^*) = \{w_1 < w_2 < \cdots < w_q < w_{q+1}\}. \end{cases}$$

*Then*

$$\lim_{N\uparrow\infty} P^{(N)}(\mathcal{G}_{h,\delta}) = 1. \tag{2.18}$$

*In particular, Theorem 2 holds.*

PROOF. Due to Lemma 5, (i) together with the condition

$$M_{h-\varepsilon}^-(B^*) = M_{h+\varepsilon}^-(B^*), \qquad M_{h-\varepsilon}^+(B^*) = M_{h+\varepsilon}^+(B^*) \tag{2.19}$$

implies (ii). Note that for all $h' > 0$ $M_{h'}^-(B^*) \subset M_{h'}^-(B) \cap [-1,1]$ $\mathbf{P}_B$-a.s., while the smallest and the largest elements of $M_{h'}^+(B^*)$ could be no $h'$-maxima of $B$. Due to (2.11) with $h, \delta$ replaced respectively by $h-\varepsilon$, $2\varepsilon$, and considering the behavior of $B^*$ at $w_1, w_{q+1}$, one gets that $\lim_{N\uparrow\infty} \mathbf{P}_B[(2.19) is\ fulfilled] = 1$. Let us suppose that the realization of $B$ does not belong to the event $\mathcal{D}_{h,\beta,\varepsilon}$ defined in Lemma 4 and that it satisfies (2.19). Then,

$$x \in [w_i, w_{i+1}] \quad \text{and} \quad B^*(x) \le B^*(u_i) + \varepsilon \Rightarrow x \in [w_i, w_{i+1}] \cap [u_i - \beta, u_i + \beta] \qquad \forall 1 \le i \le q,$$

$$x \in [u_{i-1}, u_i] \text{ and } B^*(x) \ge B^*(w_i) - \varepsilon \Rightarrow x \in [u_{i-1}, u_i] \cap [w_i - \beta, w_i + \beta] \qquad \forall 2 \le i \le q.$$

The above observations together with Lemmas 4 and 5 imply (2.18). We finally note that (1.31) follows easily from (2.14) and (2.18). □

**3. Potential theory.** In this section we recall some elementary facts of potential theory in our setting that we will need later. To this aim, we define $\mathrm{P}_{N,x}^\omega$, $\mathrm{E}_{N,x}^\omega$ respectively as the probability measure and the expectation associated to the rescaled Sinai's random walk $(X_n/N, n \ge 0)$ starting in $x \in \mathbb{Z}/N$ and with environment $\omega$.

3.1. *Equilibrium potential and capacity.* Given disjoint subsets $A, B \subset \mathbb{Z}/N$ with $|(A \cup B)^c| < \infty$ and given $\lambda \in \mathbb{C}$, the $\lambda$-*equilibrium potential* $h_{A,B}^\lambda$ is defined as the function on $\mathbb{Z}/N$ satisfying the following system:

$$\begin{cases} (\mathcal{L}^{(N)} - \lambda) h_{A,B}^\lambda(x) = 0, & \text{if } x \notin A \cup B, \\ h_{A,B}^\lambda(x) = 1, & \text{if } x \in A, \\ h_{A,B}^\lambda(x) = 0, & \text{if } x \in B. \end{cases} \tag{3.1}$$

The definition is well posed whenever the above system has a unique solution, that is, whenever $\lambda$ is not an eigenvalue of the matrix $\mathcal{L}^{(N)}((A\cup B)^c)$. In fact, it is simple to check that this last condition implies the uniqueness of the solution, while the existence is discussed below. $h^0_{A,B}$ has the probabilistic interpretation

$$h^0_{A,B}(x) = \mathrm{P}^\omega_{N,x}(\tau_A < \tau_B) \qquad \forall x \notin A\cup B, \tag{3.2}$$

where $\tau_A$ denotes the first hitting time of the set $A$, that is,

$$\tau_A \equiv \min\{n \geq 1 : X_n/N \in A\}.$$

If $\lambda \notin \sigma(\mathcal{L}^{(N)}((A\cup B)^c))$, then denoting by $g$ the restriction of $h^0_{A,B}(x)$ to $(A\cup B)^c$,

$$h^\lambda_{A,B}(x) = h^0_{A,B}(x) + \lambda(\mathcal{L}^{(N)}((A\cup B)^c) - \lambda)^{-1} g \qquad \forall x \notin A\cup B.$$

Due to the above identity, $h^\lambda_{A,B}(x)$ is holomorphic on $\mathbb{C}\setminus\sigma(\mathcal{L}^{(N)}((A\cup B)^c))$. Moreover, the following probabilistic interpretation holds:

$$h^\lambda_{A,B}(x) = \mathrm{E}^\omega_{N,x}(e^{-\ln(1-\lambda)\tau_A}\mathbb{I}_{\tau_A<\tau_B}) \qquad \forall x \notin A\cup B, \tag{3.3}$$

if

$$\lambda < \min\{\sigma(\mathbb{L}((A\cup B)^c))\}. \tag{3.4}$$

To simplify the notation, we set $h_{A,B} \equiv h^0_{A,B}$.

Given $A, B$ as above, we define the capacity, $\mathrm{cap}(A,B)$, between $A$ and $B$ as

$$\mathrm{cap}(A,B) \equiv \sum_{x\in A} \mu_N(x)(\mathcal{L}^{(N)}h_{A,B})(x) = -\sum_{x\in B}\mu_N(x)(\mathcal{L}^{(N)}h_{A,B})(x). \tag{3.5}$$

We note that $\mathrm{cap}(A,B) = \mathrm{cap}(B,A)$, since $h_{A,B} = 1 - h_{B,A}$, and that (3.1) and (3.2) imply

$$\mathrm{P}^\omega_{N,x}(\tau_A < \tau_B) = ((\mathbb{I} - \mathcal{L}^{(N)})h_{A,B})(x) \qquad \forall x \in \mathbb{Z}.$$

Due to the above identity,

$$\mathcal{L}^{(N)}h_{A,B}(x) = \begin{cases} \mathrm{P}^\omega_{N,x}(\tau_B < \tau_A), & \text{if } x \in A,\\ -\mathrm{P}^\omega_{N,x}(\tau_A < \tau_B), & \text{if } x \in B,\end{cases}$$

which gives a probabilistic interpretation of the capacity as

$$\mathrm{cap}(A,B) = \sum_{x\in A}\mu_N(x)\mathrm{P}^\omega_{N,x}(\tau_B < \tau_A).$$

In particular, if $a \notin B$,

$$\mathrm{cap}(a,B) = \mu_N(a)(\mathcal{L}^{(N)}h_{a,B})(a) = \mu_N(a)\mathrm{P}^\omega_{N,a}(\tau_B < \tau_a). \tag{3.6}$$

Another useful representation of the capacity is the well-known identity

$$\mathrm{cap}(A,B)=\sum_{x\in\mathbb{Z}/N}\mu_N(x)h_{A,B}(x)(\mathcal{L}^{(N)}h_{A,B})(x). \tag{3.7}$$

A simple renewal argument (see, e.g., [2]) gives a remarkably useful estimate of the equilibrium potential in terms of capacities,

$$h_{A,B}(x)\le\frac{\mathrm{cap}(x,A)}{\mathrm{cap}(x,B)}\qquad\forall x\notin A\cup B. \tag{3.8}$$

We consider now some particular cases of subsets $A,B$ that will be useful later. Let $\sup A=:a<b\equiv\inf B$. Then,

$$h_{A,B}(x)=\begin{cases}1, & \text{if } x\le a,\\ 0, & \text{if } x\ge b,\\ h_{a,b}(x), & \text{if } a\le x\le b,\end{cases} \tag{3.9}$$

where [with $\mu_N$ and $V^{(N)}$ defined in Section 1]

$$\begin{aligned}h_{a,b}(x)&\equiv\frac{\sum_{y=x}^{b-1/N}1/(\mu_N(y)\omega_y)}{\sum_{y=a}^{b-1/N}1/(\mu_N(y)\omega_y)}\\&=\frac{\sum_{y=x}^{b-1/N}e^{\sqrt{N}V^{(N)}(y)}}{\sum_{y=a}^{b-1/N}e^{\sqrt{N}V^{(N)}(y)}},\qquad a\le x\le b.\end{aligned} \tag{3.10}$$

For later applications, we define

$$h_{b,a}(x)\equiv1-h_{a,b}(x)=\frac{\sum_{y=a}^{x-1/N}e^{\sqrt{N}V^{(N)}(y)}}{\sum_{y=a}^{b-1/N}e^{\sqrt{N}V^{(N)}(y)}},\qquad a\le x\le b.$$

Due to the above identities,

$$\mathrm{cap}(A,B)=\mu_N(a)(\mathcal{L}^{(N)}h_{A,B})(a)=\mu_N(a)\omega_{aN}(1-h_{a,b}(a+1/N)),$$

thus implying that $\mathrm{cap}(A,B)=\mathrm{cap}(a,b)$, where

$$\mathrm{cap}(a,b)\equiv\frac{1}{\sum_{y=a}^{b-1/N}e^{\sqrt{N}V^{(N)}(y)}}. \tag{3.11}$$

Let us now suppose that $A=\{a\}$ and $B\cap(-\infty,a)\neq\varnothing$, $B\cap(a,\infty)\neq\varnothing$. By setting

$$b_1\equiv\max\{B\cap(-\infty,a)\},\qquad b_2\equiv\min\{B\cap(a,\infty)\},$$

one can check that

$$h_{a,B}(x)=\begin{cases}0, & \text{if } x\le b_1 \text{ or } x\ge b_2\\ h_{a,b_1}(x), & \text{if } b_1\le x\le a,\\ h_{a,b_2}(x), & \text{if } a\le x\le b_2,\end{cases}$$

thus implying that

$$\mathrm{cap}(a,B)=\mathrm{cap}(b_1,a)+\mathrm{cap}(a,b_2). \tag{3.12}$$

*Notational warning.* Given $a \in \mathbb{Z}/N$ and a finite subset $B \subset \mathbb{Z}/N$ such that $a \notin B$, $B \cap (-\infty, a) \neq \varnothing$ and $B \cap (a, \infty) \neq \varnothing$, then we set

$$\mathrm{cap}(a, B) \equiv \mathrm{cap}(a, B \cup (-\infty, \min B) \cup (\max B, \infty)).$$

We point out some simple estimates that will be useful in what follows. It is convenient to introduce the following notation: given positive sequences $a_N(\bar{\alpha}), b_N(\bar{\alpha})$, $N \in \mathbb{Z}_+$ (depending on some parameters $\bar{\alpha}$, including the environment $\omega$), we write

$$a_N(\bar{\alpha}) \sim [c_1(N), c_2(N), b_N(\bar{\alpha})]$$

if

$$c_1(N) b_N(\bar{\alpha}) \leq a_N(\bar{\alpha}) \leq c_2(N) b_N(\bar{\alpha}) \qquad \forall N \in \mathbb{Z}_+, \forall \bar{\alpha}.$$

Then, if $a < x < b$ belong to $\mathbb{Z}/N$,

$$h_{a,b}(x) \sim \biggl[\frac{c(\kappa)}{(b-a)N}, c'(\kappa)(b-a)N, \exp\{\sqrt{N}[V^{(N)}(z^*(x,b)) - V^{(N)}(z^*(a,b))]\}\biggr], \tag{3.13}$$

$$h_{b,a}(x) \sim \biggl[\frac{c(\kappa)}{(b-a)N}, c'(\kappa)(b-a)N, \exp\{\sqrt{N}[V^{(N)}(z^*(a,x)) - V^{(N)}(z^*(a,b))]\}\biggr], \tag{3.14}$$

$$\mathrm{cap}(a,b) \sim \biggl[\frac{c(\kappa)}{(b-a)N}, c'(\kappa), \exp\{-\sqrt{N} V^{(N)}(z^*(a,b))\}\biggr]. \tag{3.15}$$

We explain (3.13) [(3.14) and (3.15) can be justified in a similar way]. Due to (3.10), one easily gets

$$h_{a,b}(x) \sim \biggl[\frac{1}{(b-a)N}, (b-a)N, \exp\biggl\{\sqrt{N} \max_{[x,b-1/N]} V^{(N)} - \sqrt{N} \max_{[a,b-1/N]} V^{(N)}\biggr\}\biggr].$$

Due to condition (1.2), $V^{(N)}$ is a Lipschitz function with Lipschitz constant $c(\kappa)\sqrt{N}$. Therefore, the above equation implies (3.13).

3.2. *Dirichlet Green's function.* Given a finite subset $D$ in $\mathbb{Z}/N$, we define the Dirichlet Green's function $G_D$ as the $|D| \times |D|$-matrix

$$G_D \equiv (\mathcal{L}^{(N)}(D))^{-1}$$

[recall that $0 \notin \sigma(\mathcal{L}^{(N)}(D))$]. In particular, the Dirichlet problem

$$\begin{cases} \mathcal{L}^{(N)} f(z) = g(z), & \text{if } z \in D, \\ f(z) = 0, & \text{if } z \notin D, \end{cases}$$

has a unique solution, given by

$$f(z)=\sum_{y\in D} G_D(z,y)g(y) \qquad \forall z\in D. \tag{3.16}$$

It will be crucial in what follows to have an expression of $G_D$ in terms of equilibrium potentials and capacities. To this aim, we observe that, given $x\in D$, $h_{x,D^c}$ satisfies the Dirichlet problem

$$\begin{cases} \mathcal{L}^{(N)}h_{x,D^c}(y)=0, & \text{if } y\in D\setminus\{x\},\\ \mathcal{L}^{(N)}h_{x,D^c}(y)=\dfrac{\operatorname{cap}(x,D^c)}{\mu_N(x)}, & \text{if } y=x,\\ h_{x,D^c}(y)=0, & \text{if } y\in D^c \end{cases} \tag{3.17}$$

[the second identity follows from (3.6)]. Therefore, by (3.16),

$$h_{x,D^c}(z)=G_D(z,x)\frac{\operatorname{cap}(x,D^c)}{\mu_N(x)} \qquad \forall x,z\in D.$$

Since, by reversibility, $\mu_N(z)G_D(z,x)=\mu_N(x)G_D(x,z)$, the above identity is equivalent to

$$G_D(x,z)=\frac{h_{x,D^c}(z)\mu_N(z)}{\operatorname{cap}(x,D^c)} \qquad \forall x,z\in D. \tag{3.18}$$

**4. Hitting times.** By standard arguments ([9], Chapter 3), (1.2) implies that $\mathrm{E}^{\omega}_{N,x}(\tau_A)<\infty$ if $A\subset\mathbb{Z}/N$ and $|A^c|<\infty$. Due to this observation and since $\tau_a\mathbb{I}_{\tau_a<\tau_b}\le\tau_{\{a,b\}}$, we get, for $a<x<b$ in $\mathbb{Z}/N$,

$$\mathrm{E}^{\omega}_{N,x}(\tau_{\{a,b\}})<\infty, \qquad \mathrm{E}^{\omega}_{N,x}(\tau_a\mathbb{I}_{\tau_a<\tau_b})<\infty.$$

One can express the above expectations in terms of capacities and equilibrium potentials. In fact, the functions $w_1$, $w_2$ defined on $\mathbb{Z}/N$ as

$$w_1(x)\equiv\begin{cases}\mathrm{E}^{\omega}_{N,x}(\tau_{\{a,b\}}), & \text{if } a<x<b,\\ 0, & \text{if } x\notin(a,b),\end{cases}$$

$$w_2(x)\equiv\begin{cases}\mathrm{E}^{\omega}_{N,x}(\tau_a\mathbb{I}_{\tau_a<\tau_b}), & \text{if } a<x<b,\\ 0, & \text{if } x\notin(a,b),\end{cases}$$

satisfy the Dirichlet problems

$$\begin{cases}\mathcal{L}^{(N)}w_1(x)=1, & \text{if } a<x<b,\\ w_1(x)=0, & \text{if } x\notin(a,b),\end{cases}$$
$$\begin{cases}\mathcal{L}^{(N)}w_2(x)=h_{a,b}(x), & \text{if } a<x<b,\\ w_2(x)=0, & \text{if } x\notin(a,b).\end{cases} \tag{4.1}$$

Therefore, due to (3.16) and (3.18),

$$\mathrm{E}^{\omega}_{N,x}(\tau_{\{a,b\}}) = \sum_{y\in(a,b)\cap\mathbb{Z}/N} \frac{\mu_N(y)h_{x,\{a,b\}}(y)}{\mathrm{cap}(x,\{a,b\})}, \tag{4.2}$$

$$\mathrm{E}^{\omega}_{N,x}(\tau_a\mathbb{I}_{\tau_a<\tau_b}) = \sum_{y\in(a,b)\cap\mathbb{Z}/N} \frac{\mu_N(y)h_{x,\{a,b\}}(y)h_{a,b}(y)}{\mathrm{cap}(x,\{a,b\})} \tag{4.3}$$

for all $a<x<b$ in $\mathbb{Z}/N$, where $h_{x,\{a,b\}}(y)\equiv h_{x,\mathbb{Z}/N\setminus(a,b)}(y)$. Note that $h_{x,\{a,b\}}(y)=h_{x,a}(y)\mathbb{I}_{y\le x}+h_{x,b}(y)\mathbb{I}_{y>x}$.

Lemma 6. *Given $a<b$ in $\mathbb{Z}/N\cap[-1,1]$,*

$$\begin{aligned} &\max_{x\in(a,b)\cap\mathbb{Z}/N} \mathrm{E}^{\omega}_{N,x}(\tau_{\{a,b\}}) \\ &\quad\sim\Big[\frac{c(\kappa)}{N}, c'(\kappa)N^2, \\ &\qquad \exp\Big\{\sqrt{N}\max_{y\in(a,b)\cap\mathbb{Z}/N}[V_N(z^*(y,\{a,b\}))-V_N(y)]\Big\}\Big], \end{aligned}\tag{4.4}$$

$$\begin{aligned} &\max_{x\in(a,b)\cap\mathbb{Z}/N} \mathrm{E}^{\omega}_{N,x}(\tau_a|\tau_a<\tau_b) \\ &\quad\sim\Big[\frac{c(\kappa)}{N^2}, c'(\kappa)N^3, \\ &\qquad \exp\Big\{\sqrt{N}\max_{y\in(a,b)\cap\mathbb{Z}/N}[V_N(z^*(y,\{a,b\}))-V_N(y)]\Big\}\Big]. \end{aligned}\tag{4.5}$$

Proof. Since $h_{x,\{a,b\}}(x)=1$ and, for $y\in(a,b)\cap\mathbb{Z}/N$,

$$h_{x,\{a,b\}}(y)=\mathrm{P}^{\omega}_{N,y}(\tau_x<\tau_{\{a,b\}})=\begin{cases} h_{x,a}(y), & \text{if } a<y<x,\\ h_{x,b}(y), & \text{if } x<y<b,\end{cases}$$

by means of the results of the previous section, we obtain

$$\frac{\mu_N(y)h_{x,\{a,b\}}(y)}{\mathrm{cap}(x,\{a,b\})}\sim\Big[\frac{c(\kappa)}{N}, c'(\kappa)N^2, \exp\{\sqrt{N}W_{x,y}\}\Big],$$

$$\frac{\mu_N(y)h_{x,\{a,b\}}(y)h_{a,b}(y)}{\mathrm{cap}(x,\{a,b\})h_{a,b}(x)}\sim\Big[\frac{c(\kappa)}{N^2}, c'(\kappa)N^3, \exp\{\sqrt{N}\tilde{W}_{x,y}\}\Big],$$

where

$$W_{x,y}\equiv\begin{cases} V_N(z^*(x,\{a,b\})) & \\ \qquad +V_N(z^*(a,y))-V_N(z^*(a,x))-V_N(y), & \text{if } a<y\le x,\\ V_N(z^*(x,\{a,b\}))-V_N(x), & \text{if } y=x,\\ V_N(z^*(x,\{a,b\})) & \\ \qquad +V_N(z^*(y,b))-V_N(z^*(x,b))-V_N(y), & \text{if } x<y<b,\end{cases}$$

$\tilde{W}_{x,y} \equiv W_{x,y} + V_N(z^*(y,b)) - V_N(z^*(x,b))$.

We claim that

$$W_{x,y} \wedge \tilde{W}_{x,y} \leq V_N(z^*(y,\{a,b\})) - V_N(y). \tag{4.6}$$

This can be easily checked by straightforward computations as follows. By setting

$$\begin{cases} M_1 \equiv \max_{[a,y]} V_N, \qquad M_2 \equiv \max_{[y,x]} V_N, \qquad M_3 \equiv \max_{[x,b]} V_N, \\ \hfill \text{if } a < y \leq x, \\ M_1 \equiv \max_{[y,b]} V_N, \qquad M_2 \equiv \max_{[x,y]} V_N, \qquad M_3 \equiv \max_{[a,x]} V_N, \\ \hfill \text{if } x < y < b, \end{cases} \tag{4.7}$$

we can write

$$\begin{aligned} W_{x,y} + V_N(y) &= (M_1 \vee M_2) \wedge M_3 + M_1 - M_1 \vee M_2, \\ V_N(z^*(y,\{a,b\})) &= M_1 \wedge (M_2 \vee M_3). \end{aligned} \tag{4.8}$$

Moreover, the following inequalities hold:

$$\begin{cases} \tilde{W}_{x,y} \leq W_{x,y}, & \text{if } x \leq y, \\ \tilde{W}_{x,y} = W_{x,y} + M_2 \vee M_3 - M_3, & \text{if } x > y. \end{cases}$$

At this point, (4.6) can be checked by considering the six possible orderings of $M_1, M_2, M_3$:

Having proved (4.6), (4.4) and (4.5) can be easily derived from (4.2) and (4.3) together with the identity $h_{a,b}(x) = \mathrm{P}^{\omega}_{N,x}(\tau_a < \tau_b)$ and the observation that $W_{y,y} = \tilde{W}_{y,y} = V_N(z^*(y,\{a,b\})) - V_N(y)$. □

We conclude this section by recalling a generalization of (4.3). Given two disjoint subsets $A, B \subset \mathbb{Z}/N$ with $|(A \cup B)^c| < \infty$, the function $w$ on $\mathbb{Z}/N$ defined as

$$w(x) \equiv \begin{cases} \mathrm{E}^{\omega}_{N,x}(\tau_A \mathbb{I}_{\tau_A < \tau_B}), & \text{if } x \notin (A \cup B), \\ 0, & \text{if } x \in A \cup B, \end{cases}$$

is a finite function satisfying the Dirichlet problem

$$\begin{cases} \mathcal{L}^{(N)} w(x) = h_{A,B}(x), & \text{if } x \notin A \cup B, \\ w(x) = 0, & \text{if } x \in A \cup B. \end{cases}$$

In particular, due to (3.16) and (3.18), we get

$$\begin{aligned} \mathrm{E}^{\omega}_{N,x}(\tau_A \mathbb{I}_{\tau_A < \tau_B}) &= \sum_{y \notin A \cup B} G_{(A \cup B)^c}(x,y) h_{A,B}(y) \\ &= \sum_{y \notin A \cup B} \frac{\mu_N(y) h_{x,A \cup B} h_{A,B}(y)}{\mathrm{cap}(x, A \cup B)}. \end{aligned} \tag{4.9}$$

**5. Principal eigenvalues and eigenvectors.** In this section we fix $h,\delta>0$ and $V_N\in\mathcal{A}_{h,\delta}$ and we usually omit the index $h$ and the reference to the path $V_N$ from the standard notation. In particular, we write $M^-\equiv M_h^-(V_N)=\{x_1,x_2,\dots,x_q\}$, where $x_1,x_2,\dots,x_q$ is the labeling satisfying condition (1.16). Moreover, we set, for $0\le k\le q$,

$$\begin{cases} M_k^-\equiv\{x_1,\dots,x_k\},\\ S_k\equiv M_k^-\cup\{-1,1\},\\ S_k^*\equiv M_k^-\cup(\mathbb{Z}/N\setminus(-1,1)).\end{cases}$$

Our target here is to study the principal eigenvalue $\bar\lambda_k^{(N)}$ and the related eigenvector of the generator on $\mathbb{Z}/N$ with Dirichlet conditions on $S_k^*$. To this aim, recall the definitions (1.32) and (1.33):

$$\mathcal{L}_k^{(N)}\equiv\mathcal{L}^{(N)}((S_k^*)^c),\qquad \bar\lambda_k^{(N)}\equiv\min\sigma(\mathcal{L}_k^{(N)}).$$

Moreover, observe that $\mathcal{L}_0^{(N)}=\mathcal{L}^{(N)}(I_N)$ and that, due to Corollary 3 in Appendix C, $\bar\lambda_0^{(N)}<\bar\lambda_1^{(N)}<\cdots<\bar\lambda_q^{(N)}$.

To get an upper bound on $\bar\lambda_k^{(N)}$, we recall its variational characterization:

$$\bar\lambda_k^{(N)}=\inf_{\substack{f\in\mathbb{R}^{\mathbb{Z}/N}\\ f\equiv0\text{ on }S_k^*,f\not\equiv0}}\frac{(f,\mathcal{L}^{(N)}f)}{\|f\|_2^2},\tag{5.1}$$

where $(\cdot,\cdot)$ and $\|\cdot\|_2$ denote respectively the scalar product and the norm in $L^2(\mathbb{Z}/N,\mu_N)$.

A lower bound can be obtained using a Donsker–Varadhan like argument as explained in [3], Lemma 4.2:

$$\bar\lambda_k^{(N)}\ge\frac{1}{\sup_{x\notin S_k^*}\mathrm{E}_{N,x}^{\omega}(\tau_{S_k^*})}.\tag{5.2}$$

LEMMA 7. *If $V_N\in\mathcal{A}_{h,\delta}$, then*

$$c(\kappa)N^{-2}e^{-\sqrt{N}d_{k+1}}\le\bar\lambda_k^{(N)}\le c'(\kappa)e^{-\sqrt{N}d_{k+1}}\qquad\forall k:0\le k\le q-1,\tag{5.3}$$

$$\bar\lambda_q^{(N)}\ge c(\kappa)N^{-2}e^{-h\sqrt{N}}.\tag{5.4}$$

*In particular,*

$$\bar\lambda_k^{(N)}\le c(\kappa)N^2e^{-\delta\sqrt{N}}\bar\lambda_{k+1}^{(N)}\qquad\forall k:0\le k\le q-1.\tag{5.5}$$

PROOF. We first derive an upper bound for $\bar\lambda_k^{(N)}$, for $0\le k\le q-1$. Let

$$a\equiv z^*([-1,x_{k+1})\cap S_k,x_{k+1}),\qquad b\equiv z^*((x_{k+1},1]\cap S_k,x_{k+1}),$$

where the saddle points are w.r.t. $V_N$, and set $f \equiv \mathbb{I}_{(a,b)\cap\mathbb{Z}/N}$. Then

$$\|f\|_2^2 \geq \mu_N(x_{k+1}) \geq c(\kappa)e^{-\sqrt{N}V_N(x_{k+1})},$$

while, by (1.13),

$$(f, \mathcal{L}^{(N)} f) = e^{-\sqrt{N}V_N(a)} + e^{-\sqrt{N}V_N(b-1/N)} \leq e^{-\sqrt{N}V_N(a)} + c(\kappa)e^{-\sqrt{N}V_N(b)}.$$

Since $f \equiv 0$ on $S_k^*$, by (5.1), we get

$$\bar{\lambda}_k^{(N)} \leq c(\kappa)e^{-\sqrt{N}d_{k+1}}.$$

To bound $\bar{\lambda}_k^{(N)}$ from below, we derive from (5.2) and (4.4) that

$$\bar{\lambda}_k^{(N)} \geq c(\kappa)N^{-2}\exp\Big\{-\sqrt{N}\max_{x\notin S_k^*}[V_N(z^*(x,S_k)) - V_N(x)]\Big\}. \tag{5.6}$$

Due to Lemma 2, the maximum in the above expression is achieved for some $x \in M^- \setminus S_k = \{x_{k+1}, \dots, x_q\}$. Then, due to (1.16), the maximum has to be achieved at $x = x_{k+1}$, thus concluding the proof of (5.3). To prove (5.4), we observe that (5.6) remains true for $k = q$. This, together with the definition of $h$-extrema, implies (5.4). Finally, (5.5) with $0 \leq k \leq q-2$ follows from (1.16), (1.18) and (5.3), while, for $k = q-1$, it follows from (5.3), (5.4) and (1.15). □

PROPOSITION 4. *Given $V_N \in \mathcal{A}_{h,\delta}$ and $k \in \{1, 2, \dots, q\}$, $\bar{\lambda}_{k-1}^{(N)}$ is a simple eigenvalue of $\mathcal{L}_{k-1}^{(N)}$ with eigenfunction $h^{\lambda}_{x_k,S^*_{k-1}}$, where $\lambda \equiv \bar{\lambda}_{k-1}^{(N)}$, and*

$$\begin{aligned} &\frac{\operatorname{cap}(x_k, S_{k-1}^*)}{\|h_{x_k,S_{k-1}^*}\|_2^2}(1 - c(\kappa)N^2 e^{-\delta\sqrt{N}}) \\ &\qquad \leq \bar{\lambda}_{k-1}^{(N)} \leq \frac{\operatorname{cap}(x_k, S_{k-1}^*)}{\|h_{x_k,S_{k-1}^*}\|_2^2}(1 + c(\kappa)N^2 e^{-\delta\sqrt{N}}). \end{aligned} \tag{5.7}$$

PROOF. Note that $\mathcal{L}_k^{(N)}$ is obtained from $\mathcal{L}_{k-1}^{(N)}$ by adding Dirichlet conditions on $x_k$, since $S_k^* = S_{k-1}^* \cup \{x_k\}$. In particular, $\lambda < \bar{\lambda}_k^{(N)}$ is an eigenvalue of $\mathcal{L}_{k-1}^{(N)}$ with eigenfunction $f \in \mathbb{R}^{\mathbb{Z}/N}$ iff $\exists \phi \in \mathbb{R}$ such that

$$\begin{cases} (\mathcal{L}^{(N)} - \lambda)f(y) = 0, & \text{if } y \notin S_k^*, \\ f(y) = \phi, & \text{if } y = x_k, \\ f(y) = 0, & \text{if } y \in S_{k-1}^*, \end{cases} \tag{5.8}$$

and

$$(\mathcal{L}^{(N)} - \lambda)f(x_k) = 0. \tag{5.9}$$

Note that (5.8) is equivalent to the identity $f = \phi h^{\lambda}_{x_k, S^*_{k-1}}$. Since $f$ is an eigenfunction, $\phi \neq 0$ and, therefore, $\phi$ can be taken equal to 1. Since $\bar{\lambda}^{(N)}_{k-1} < \bar{\lambda}^{(N)}_{k}$ due to Corollary 3, this implies the first statement of the proposition. Let us write

$$h^{\lambda} \equiv h^{\lambda}_{x_k, S^*_{k-1}}, \qquad h \equiv h_{x_k, S^*_{k-1}}, \qquad \delta h^{\lambda} \equiv h^{\lambda} - h, \tag{5.10}$$

and show that $h^{\lambda}$ can be approximated by $h$ for $\lambda < \bar{\lambda}^{(N)}_{k}$. In fact, the function $\delta h^{\lambda}$ satisfies the system

$$\begin{cases} \delta h^{\lambda}(y) = 0, & \text{if } y \in S^*_k, \\ (\mathcal{L}^{(N)} - \lambda)\delta h^{\lambda}(y) = \lambda h(y), & \text{if } y \notin S^*_k, \end{cases} \tag{5.11}$$

thus implying $\delta h^{\lambda} = \lambda(\mathcal{L}^{(N)}_k - \lambda)^{-1} h$ and, therefore,

$$\begin{aligned} \|\delta h^{\lambda}\|_2 &= \|\delta h^{\lambda}\|_{L^2((S^*_k)^c, \mu_N)} \\ &\leq \frac{\lambda}{\bar{\lambda}^{(N)}_k - \lambda} \|h\|_{L^2((S^*_k)^c, \mu_N)} \leq \frac{\lambda}{\bar{\lambda}^{(N)}_k - \lambda} \|h\|_2. \end{aligned} \tag{5.12}$$

Due to (5.9), (5.11) and the identity

$$(\mathcal{L}^{(N)} h)(x) = A \delta_{x, x_k} \qquad \forall x \notin (S^*_{k-1})^c, \text{ where } A \equiv \frac{\text{cap}(x_k, S^*_{k-1})}{\mu_N(x_k)}$$

[which follows from (3.6)], we obtain

$$(\mathcal{L}^{(N)} - \lambda)\delta h^{\lambda} = \lambda h - A\delta_{x, x_k} \qquad \text{on } (S^*_{k-1})^c. \tag{5.13}$$

By taking the scalar product in $\mathbb{L}^2(\mu_N)$ with $h$, which is zero on $S^*_{k-1}$, we get

$$(h, \mathcal{L}^{(N)} \delta h^{\lambda}) - \lambda(h, \delta h^{\lambda}) = \lambda(h, h) - A(h, \delta_{x, x_k}). \tag{5.14}$$

Due to reversibility, the first addendum is zero. Moreover, since $A(h, \delta_{x, x_k}) = \text{cap}(x_k, S^*_{k-1})$, (5.12) and (5.14) imply

$$\left| \lambda - \frac{\text{cap}(x_k, S^*_{k-1})}{\|h\|_2^2} \right| \leq \frac{\lambda^2}{\bar{\lambda}^{(N)}_k - \lambda}. \tag{5.15}$$

The assertion follows now by considering the case $\lambda \equiv \bar{\lambda}^{(N)}_{k-1}$ and using (5.5).
□

We conclude the section with a description of the principal eigenfunction of $\mathcal{L}^{(N)}_{k-1}$.

PROPOSITION 5. *If $V_N \in \mathcal{A}_{h,\delta}$ and $k \in \{1,2,\dots,q\}$, then the function $h^{\lambda}_{x_k,S^*_{k-1}}$, $\lambda = \bar{\lambda}^{(N)}_{k-1}$, satisfies*

$$h_{x_k,S^*_{k-1}}(y) \le h^{\lambda}_{x_k,S^*_{k-1}}(y) \le h_{x_k,S^*_{k-1}}(1 + c(\kappa)N^5 e^{-\delta\sqrt{N}}). \tag{5.16}$$

PROOF. For simplicity of notation, let $\lambda \equiv \bar{\lambda}^{(N)}_{k-1}$ and let $h^{\lambda}, h, \delta h^{\lambda}$ be defined as in (5.10).

Since $\bar{\lambda}^{(N)}_{k-1} < \bar{\lambda}^{(N)}_{k}$ and the latter is the principal eigenvalue of $\mathcal{L}^{(N)}_k$, $h^{\lambda}$ admits the probabilistic interpretation (3.3). Comparing it with (3.2), one gets the inequality on the left in (5.16).

To prove the inequality on the right, we observe that $\delta h^{\lambda}$ satisfies (5.11), and therefore,

$$\begin{cases} \mathcal{L}^{(N)} \delta h^{\lambda}(y) = \lambda h^{\lambda}(y), & \text{if } y \notin S^*_k, \\ \delta h^{\lambda}(y) = 0, & \text{if } y \in S^*_k. \end{cases}$$

Due to the above Dirichlet problem and (3.16), we obtain

$$\frac{h^{\lambda}(y)}{h(y)} = 1 + \frac{\lambda}{h(y)} \sum_{z \notin S^*_k} G_{(S^*_k)^c}(y,z) \frac{h^{\lambda}(z)}{h(z)} h(z) \qquad \forall y \notin S^*_k.$$

The above identity implies

$$\frac{h^{\lambda}(y)}{h(y)} \le 1 + \frac{\lambda M}{h(y)} \sum_{z \notin S^*_k} G_{(S^*_k)^c}(y,z) h(z),$$

where $M \equiv \max_{z \notin S^*_k} \frac{h^{\lambda}(z)}{h(z)}$. From the above inequality, (3.2) and (4.9), we derive

$$M \le 1 + \lambda M \max_{y \notin S^*_k} \mathrm{E}^{\omega}_{N,y}(\tau_{x_k} | \tau_{x_k} < \tau_{S^*_{k-1}}). \tag{5.17}$$

Due to Lemma 6,

$$\begin{aligned} &\max_{y \notin S^*_k} \mathrm{E}^{\omega}_{N,y}(\tau_{x_k} | \tau_{x_k} < \tau_{S^*_{k-1}}) \\ &\qquad \le c(\kappa) N^3 \exp\Big\{ \sqrt{N} \max_{y \notin S^*_k} [V_N(z^*(y,S_k)) - V_N(y)] \Big\}. \end{aligned} \tag{5.18}$$

If $1 \le k < q$, then Lemma 2 and (1.16) imply that the above maximum is achieved for $y = x_{k+1}$. (5.3) then implies

$$\mathrm{E}^{\omega}_{N,y}(\tau_{x_k} | \tau_{x_k} < \tau_{S^*_{k-1}}) \le \frac{c(\kappa) N^3}{\bar{\lambda}^{(N)}_k}, \qquad \text{if } 1 \le k < q.$$

Therefore, due to (5.5) and (5.17), we get $M \le 1 + c(\kappa) M N^5 e^{-\delta\sqrt{N}}$, which implies (5.16). If $k = q$, then the r.h.s. of (5.18) can be bounded by $e^{h\sqrt{N}}$. Due to (5.3) and condition (1.15), one get that $M \le 1 + c(\kappa) M N^3 e^{-\delta\sqrt{N}}$, which implies (5.16). □

**6. The set $\boldsymbol{\sigma(\mathcal{L}^{(N)}(I_N)) \cap (0, \bar\lambda_k^{(N)})}$.** As in Section 5, we fix $h, \delta > 0$, $V_N \in \mathcal{A}_{h,\delta}$ and we usually omit the index $h$ and the reference to the path $V_N$ from the standard notation.

Given $1 \le k \le q$, $\lambda < \bar\lambda_k^{(N)}$ is an eigenvalue of $\mathcal{L}^{(N)}(I_N) = \mathcal{L}_0^{(N)}$ with eigenvector $f^\lambda \in \mathbb{R}^{\mathbb{Z}/N}$ if and only if, for suitable constants $\phi^\lambda(y)$ with $y \in M_k^-$,

$$\begin{cases} (\mathcal{L}^{(N)} - \lambda) f^\lambda(y) = 0, & \text{if } y \notin S_k^*, \\ f^\lambda(y) = \phi^\lambda(y), & \text{if } y \in M_k^-, \\ f^\lambda(y) = 0, & \text{if } y \in S_k^* \setminus M_k^-, \end{cases} \tag{6.1}$$

[note that $S_k^* \setminus M_k^- = \mathbb{Z}/N \setminus (-1,1)$] and

$$(\mathcal{L}^{(N)} - \lambda) f^\lambda(y) = 0 \qquad \forall y \in M_k^-. \tag{6.2}$$

System (6.1) is equivalent to the identity

$$f^\lambda(y) = \sum_{x \in M_k^-} \phi_x^\lambda h^\lambda_{x, S_k^* \setminus \{x\}}(y) \qquad \forall y \in \mathbb{Z}/N. \tag{6.3}$$

It is convenient to introduce a shortened notation by defining

$$h_x^\lambda \equiv h^\lambda_{x, S_k^* \setminus \{x\}}, \qquad h_x \equiv h_{x, S_k^* \setminus \{x\}}$$

(note that $h_x^\lambda$ depends on $k$). Assuming (6.3), condition (6.2) is equivalent to

$$\sum_{x \in M_k^-} \phi_x^\lambda ((\mathcal{L}^{(N)} - \lambda) h_x^\lambda)(y) = 0 \qquad \forall y \in M_k^-. \tag{6.4}$$

Let us denote by $\mathcal{E}_k(\lambda)$ the $k \times k$-matrix

$$(\mathcal{E}_k(\lambda))_{x,z} = ((\mathcal{L}^{(N)} - \lambda) h_z^\lambda)(x) \qquad \forall x, z \in M_k^-, \tag{6.5}$$

and by $\hat{\mathcal{E}}_k(\lambda)$ the $k \times k$-matrix

$$(\hat{\mathcal{E}}_k(\lambda))_{x,z} = \frac{1}{\mu_N(x)} \frac{(\mathcal{E}_k(\lambda))_{x,z}}{\|h_x\|_2 \|h_z\|_2}. \tag{6.6}$$

Note that both $\mathcal{E}_k(\lambda)$ and $\hat{\mathcal{E}}_k(\lambda)$ are well defined and holomorphic on $\mathbb{C} \setminus \sigma(\mathcal{L}_k^{(N)})$.

Then the above observations imply:

LEMMA 8. *$\lambda < \bar\lambda_k^{(N)}$ is an eigenvalue of $\mathcal{L}_0^{(N)}$ iff* $\det(\mathcal{E}_k(\lambda)) = 0$. *In this case, $f^\lambda : \mathbb{Z}/N \to \mathbb{R}$ is an eigenvector of $\mathcal{L}_0^{(N)}$ with eigenvalue $\lambda$ iff $f^\lambda = \sum_{x \in M_k^-} \phi_x^\lambda h_x^\lambda$ for some eigenvector $\phi^\lambda : M_k^- \to \mathbb{R}$ of $\mathcal{E}_k(\lambda)$ with eigenvalue* 0. *Moreover,* $\det(\mathcal{E}_k(\lambda)) = 0$ *iff* $\det(\hat{\mathcal{E}}_k(\lambda)) = 0$, *and $\mathcal{E}_k(\lambda)\phi = 0$ iff $\hat{\mathcal{E}}_k(\lambda)\hat\phi = 0$, where $\hat\phi_x = \phi_x \|h_x\|_2$ for all $x \in M_k^-$.*

Note that we can write, for any $x, z \in M_k^-$,

$$(\mathcal{E}_k(\lambda))_{x,z} = \mu_N(x)((h_x, \mathcal{L}^{(N)} h_z) - \lambda(h_x, h_z) - \lambda(h_x, \delta h_z^\lambda)), \tag{6.7}$$

where $\delta h_z^\lambda(y) \equiv h_z^\lambda(y) - h_z(y)$, respectively

$$(\hat{\mathcal{E}}_k(\lambda))_{x,y} = \mathcal{K}_{x,z}^{(k)} - \lambda \mathbb{I}_{x,z} - \lambda A_{x,z}^{(k)} - \lambda B_{x,z}^{(k)} \qquad \forall x, z \in M_k^-, \tag{6.8}$$

where

$$\mathcal{K}_{x,x}^{(k)} \equiv \frac{(h_x, \mathcal{L}^{(N)} h_z)}{\|h_x\|_2 \|h_z\|_2}, \qquad A_{x,z}^{(k)} \equiv \frac{(h_x, h_z)}{\|h_x\|_2 \|h_z\|_2}(1 - \delta_{x,z}),$$

$$B_{x,z}^{(k)} \equiv \frac{(h_x, \delta h_z^\lambda)}{\|h_x\|_2 \|h_z\|_2}.$$

The above $k \times k$-matrix $\mathcal{K}^{(k)}$ is called the *normalized capacity matrix* $\mathcal{K}^{(k)}$. Due to (3.7), $\mathcal{K}^{(k)}$ is a symmetric matrix with diagonal elements given by

$$\mathcal{K}_{x,x}^{(k)} = \|h_x\|_2^{-2} \operatorname{cap}(x, S_k^* \setminus \{x\}). \tag{6.9}$$

Note that due to (5.7), if $V_N \in A_{h,\delta}$, then

$$\left| \frac{\mathcal{K}_{x_k,x_k}^{(k)}}{\bar\lambda_{k-1}^{(N)}} - 1 \right| \le c N^2 e^{-\delta\sqrt{N}}. \tag{6.10}$$

Moreover, ordering the entries of $\mathcal{K}^{(k)}$ along the increasing order of $x_1, x_2, \ldots, x_k$, the matrix $\mathcal{K}^{(k)}$ is Jacobian. Namely, let $\{x_1, x_2, \ldots, x_k\} = \{u_1, u_2, \ldots, u_k\}$ with $u_1 < u_2 < \cdots < u_k$ and set $A_{i,j} \equiv \mathcal{K}_{u_i,u_j}^{(k)}$. Then $A = (A_{i,j})_{1 \le i,j \le d}$ is symmetric and $A_{i,j} = 0$ if $|i - j| > 1$. This follows easily from the following observation: setting $u_0 \equiv -1$, $u_{k+1} \equiv 1$, then the support of $h_{u_i}$, $\mathcal{L}^{(N)} h_{u_i}$ is respectively given by $(u_{i-1}, u_{i+1})$, $[u_{i-1}, u_{i+1}]$ for all $i = 1, \ldots, k$.

The following proposition will allow us to think of $\hat{\mathcal{E}}_k(\lambda)$, with $\lambda \in \mathbb{C} \setminus \sigma(\mathcal{L}_k^{(N)})$ small, as obtained by a small perturbation from the matrix

$$(K_{x_k,x_k}^{(k)} \delta_{x_k,x} \delta_{x_k,z} - \lambda \delta_{x,z})_{x,z}, \qquad \text{with } x, z \in M_k^-.$$

PROPOSITION 6. *If $V_N \in \mathcal{A}_{h,\delta}$ and $1 \le k \le q$, then*

$$K^{(k)}_{x_j,x_j} \le CN^2 e^{-\delta\sqrt{N}} K^{(k)}_{x_k,x_k} \qquad \forall 1 \le j < k, \tag{6.11}$$

$$K^{(k)}_{x_i,x_j} \le CN^2 e^{(-\delta/2)\sqrt{N}} K^{(k)}_{x_k,x_k} \qquad \forall 1 \le i,j \le k, (i,j) \ne (k,k), \tag{6.12}$$

$$A^{(k)}_{x_i,x_j} \le CN^2 e^{(-\delta/4)\sqrt{N}} \qquad \forall 1 \le i,j \le k, i \ne j \tag{6.13}$$

$$|B^{(k)}_{x_i,x_j}| \le \frac{|\lambda|}{\operatorname{dist}(\lambda, \sigma(\mathcal{L}^{(N)}_k))} \qquad \forall \lambda \in \mathbb{C} \setminus \sigma(\mathcal{L}^{(N)}_k), \forall 1 \le i,j \le k. \tag{6.14}$$

PROOF. Proposition 6 is analogous to the corresponding statements in [5] and we refer for the details to that paper.

The main ingredient of the proof are the nondegeneracy conditions. In fact, equation (6.11) follows from (6.10) and (5.5). Equation (6.12) follows from (6.11) using the Schwarz inequality, $|(h_{x_i}, \mathcal{L}^{(N)} h_{x_j})| \le (h_{x_i}, \mathcal{L}^{(N)} h_{x_i})^{1/2} (h_{x_j}, \mathcal{L}^{(N)} h_{x_j})^{1/2}$, that is, $|\mathcal{K}^{(k)}_{x_i,x_j}| \le \sqrt{\mathcal{K}^{(k)}_{x_i,x_i} \mathcal{K}^{(k)}_{x_j,x_j}}$.

Equation (6.13) is just a statement that the functions $h_x$ and $h_y$ are almost orthogonal. This is completely analogous to Lemma 4.5 of [5].

To prove equation (6.14), note that $\delta h^\lambda_{x_j} = h^\lambda_{x_j} - h_{x_j}$ satisfies the Dirichlet problem

$$\begin{cases} (\mathcal{L}^{(N)} - \lambda)\delta h^\lambda_{x_j}(y) = \lambda h_{x_j}(y), & \text{if } y \notin S^*_k, \\ \delta h^\lambda(y) = 0, & \text{if } y \in S^*_k. \end{cases}$$

Thus, $\delta h^\lambda_{x_j} = \lambda(\mathcal{L}^{(N)}_k - \lambda)^{-1} h_{x_j}$, implying

$$\|\delta h^\lambda_{x_j}\|_2 \le \frac{|\lambda|}{\operatorname{dist}(\lambda, \sigma(\mathcal{L}^{(N)}_k))} \|h_{x_j}\|_2.$$

(6.14) now follows from the Schwarz inequality. □

We can now prove the main result of this section:

THEOREM 6. *If $V_N \in \mathcal{A}_{h,\delta}$, $q \equiv |M^-_h| \le Q$, $N \ge N(\delta, Q)$, then the following hold:*

$$\sigma(\mathcal{L}^{(N)}_0) \cap [0, \bar\lambda^{(N)}_q) = \{\bar\lambda^{(N)}_0 = \lambda^{(N)}_1 < \lambda^{(N)}_2 < \cdots < \lambda^{(N)}_q\} \tag{6.15}$$

*and*

$$\left| \frac{\lambda^{(N)}_k}{\bar\lambda^{(N)}_{k-1}} - 1 \right| \le e^{(-\delta/10)\sqrt{N}} \qquad \forall k = 1, 2, \ldots, q. \tag{6.16}$$

*Moreover,* $\lambda_k^{(N)}$ *is a simple eigenvalue with normalized eigenfunction* $\psi_k^{(N)}$*:*

$$\psi_k^{(N)} = a_k^{(k)} \frac{h^{\lambda}_{x_k, S^*_{k-1}}}{\|h^{\lambda}_{x_k, S^*_{k-1}}\|_2} + \sum_{j=1}^{k-1} a_j^{(k)} \frac{h^{\lambda}_{x_j, S^*_k \setminus \{x_j\}}}{\|h^{\lambda}_{x_j, S^*_k \setminus x_j}\|_2}, \qquad \lambda \equiv \lambda_k^{(N)}, \tag{6.17}$$

*where* $a_j^{(k)}$*,* $1 \le j \le k$*, are constants satisfying*

$$1 - e^{(-\delta/10)\sqrt{N}} \le a_k^{(k)} \le 1, \qquad |a_j^{(k)}| \le e^{(-\delta/10)\sqrt{N}} \qquad \forall 1 \le j \le k-1. \tag{6.18}$$

*In particular,*

$$\left\| \psi_k^{(N)} - \frac{h_{x_k, S^*_{k-1}}}{\|h_{x_k, S^*_{k-1}}\|_2} \right\|_2 \le e^{(-\delta/10)\sqrt{N}}. \tag{6.19}$$

PROOF. To prove that the set $\sigma(\mathcal{L}_0^{(N)}) \cap [0, \bar\lambda_q^{(N)})$ has cardinality at least $q$, we apply Lagrange's theorem [22] stating the following: let $\varphi$ be a holomorphic function defined on a open set $D \subset \mathbb{C}$ containing a point $a$. If there exists a contour $\gamma$ around $a$ and inside $D$ such that $|\varphi(z)| < |z - a|$ for any $z$ in the support of $\gamma$, then the equation

$$a - z + \varphi(z) = 0 \tag{6.20}$$

has a unique solution in the interior of $\gamma$.

Fix $1 \le k \le q$ and recall the definition of $\bar\lambda_k^{(N)}$ given in (1.33). Since $\mathcal{L}_0^{(N)}$ has only positive eigenvalues and due to Lemma 8, $\lambda < \bar\lambda_k^{(N)}$ is an eigenvalue of $\mathcal{L}_0^{(N)}$ if and only if

$$\det(\hat{\mathcal{E}}_k(\lambda)) = 0. \tag{6.21}$$

Let us define

$$D_k \equiv \{\lambda \in \mathbb{C} : |\Im(\lambda)| < \Re(\lambda), e^{(-\delta/8)\sqrt{N}} \bar\lambda_{k-1}^{(N)} < \Re(\lambda) \le e^{(-\delta/4)\sqrt{N}} \bar\lambda_k^{(N)}\}.$$

Note that, due to (5.5), $D_k$ is nonempty if $N \ge N(\delta)$. Moreover, for $N \ge N(\delta)$,

$$\frac{|\lambda|}{\operatorname{dist}(\lambda, \sigma(\mathcal{L}_k^{(N)}))} \le \sqrt{2} \frac{\Re(\lambda)}{\bar\lambda_k^{(N)} - \Re(\lambda)} \le 2 e^{(-\delta/4)\sqrt{N}} \qquad \forall \lambda \in D_k.$$

Due to (6.8), (6.10), Proposition 6 and the above estimate, for all $\lambda \in D_k$, we can write

$$\hat{\mathcal{E}}_k(\lambda) = V^{(k)}(\lambda) + W^{(k)}(\lambda), \tag{6.22}$$

where, for all $x, y \in M_k^-$ and for $N \geq N(\delta)$,

$$V_{x,y}^{(k)} = K_{x_k,x_k}^{(k)} \delta_{x,x_k} \delta_{y,x_k} - \lambda \delta_{x,y}, \tag{6.23}$$

$$K_{x_k,x_k}^{(k)} \leq \bar{\lambda}_{k-1}^{(N)} (1 + e^{(-\delta/2)\sqrt{N}}), \tag{6.24}$$

$$|W_{x,y}^{(k)}(\lambda)| \leq cN^2 e^{(-\delta/4)\sqrt{N}} (|\lambda| + \bar{\lambda}_{k-1}^{(N)}) \leq 2cN^2 e^{(-\delta/8)\sqrt{N}} |\lambda|. \tag{6.25}$$

Note that the last inequality in (6.25) follows from the definition of $D_k$.

In what follows we suppose $N \geq N(\delta)$ such that $2cN^2 e^{(-\delta/8)\sqrt{N}} < e^{(-\delta/9)\sqrt{N}}$, thus implying that $|W_{x,y}^{(k)}(\lambda)| < e^{(-\delta/9)\sqrt{N}} |\lambda|$.

Let us write

$$\begin{aligned} &\det(\hat{\mathcal{E}}_k(\lambda)) \\ &\quad = \sum_\tau (-1)^{\operatorname{sgn}(\tau)} (\hat{\mathcal{E}}_k(\lambda))_{x_1,\tau(x_1)} (\hat{\mathcal{E}}_k(\lambda))_{x_2,\tau(x_2)} \cdots (\hat{\mathcal{E}}_k(\lambda))_{x_k,\tau(x_k)}, \end{aligned} \tag{6.26}$$

where $\tau$ varies among the permutations of $M_k^-$ and $\operatorname{sgn}(\tau)$ denotes its sign. Let us consider the addendum in the r.h.s. associated to $\tau$ equal to the identity, that is,

$$(K_{x_k,x_k}^{(k)} - \lambda + W_{x_k,x_k}^{(k)}(\lambda)) \prod_{j=1}^{k-1} (-\lambda + W_{x_j,x_j}^{(k)}(\lambda)).$$

It can be written as $K_{x_k,x_k}^{(k)}(-\lambda)^{k-1} + (-\lambda)^k + \tilde{\phi}(\lambda)$, where $\tilde{\phi}(\lambda)$ is a holomorphic function on $D_k$ with $|\tilde{\phi}(\lambda)| \leq c(k)(|\lambda| + \bar{\lambda}_{k-1}^{(N)})|\lambda|^{k-1} e^{-(\delta/9)\sqrt{N}}$.

Note that if the permutation $\tau$ is different from the identity and if $\lambda \in D_k$, then

$$\begin{aligned} |(\hat{\mathcal{E}}_k(\lambda))_{x_j,\tau(x_j)}| &\leq |\lambda|(1 + e^{(-\delta/9)\sqrt{N}}) && \forall 1 \leq j < k, \\ |(\hat{\mathcal{E}}_k(\lambda))_{x_{j_0},\tau(x_{j_0})}| &\leq |\lambda| e^{(-\delta/9)\sqrt{N}}, && \text{for some } 1 \leq j_0 \leq k, \\ |(\hat{\mathcal{E}}_k(\lambda))_{x_k,\tau(x_k)}| &\leq (|\lambda| + \bar{\lambda}_{k-1}^{(N)})(1 + e^{(-\delta/9)\sqrt{N}}) \end{aligned}$$

[in the last estimate we have used (6.24)]. The above observations imply that, for $\lambda \in D_k$,

$$\det(\hat{\mathcal{E}}_k(\lambda))/(-\lambda)^{k-1} = K_{x_k,x_k}^{(k)} - \lambda + \phi(\lambda), \tag{6.27}$$

where $\phi(\lambda)$ is a holomorphic function with

$$|\phi(\lambda)| \leq c'(k) e^{(-\delta/9)\sqrt{N}} (|\lambda| + \bar{\lambda}_{k-1}^{(N)}). \tag{6.28}$$

Let $\gamma$ be the circle in $\mathbb{C}$ around $K_{x_k,x_k}^{(k)}$ of radius $r = 6c'(k)e^{-(\delta/9)\sqrt{N}} \bar{\lambda}_{k-1}^{(N)}$. Due to this choice, (6.10) and (5.5), if $N \geq N(Q, \delta)$ and $\lambda \in \operatorname{supp}(\gamma)$, then $\lambda \in D_k$ and the r.h.s. of (6.28) is strictly bounded from above by $r = |K^{(k)}_{x_k,x_k} - \lambda|$. Therefore, by Lagrange's theorem, there is one and only one eigenvalue $\lambda_k^{(N)}$ of $\mathcal{L}_0^{(N)}$ inside $\gamma$ [thus implying that $\lambda_k^{(N)} \in D_k$].

Since all the sets $D_k$ are disjoint,

$$|\sigma(\mathcal{L}_0^{(N)}) \cap [0, \bar{\lambda}_k^{(N)})| \geq q,$$

while, due to Proposition 9, the l.h.s. is not larger than $q$. That completes the proof of (6.15).

Let $\underline{a} = (a_x)_{x \in M_k^-}$ be a (right) eigenvector of $\hat{\mathcal{E}}_k(\lambda_k^{(N)})$ with eigenvalue 0. We can suppose that $\underline{a}$ is normalized, that is, $\sum_{j=1}^k |a_{x_j}|^2 = 1$, and $a_{x_k} \geq 0$. For $1 \leq i < k$, the identity $(\hat{\mathcal{E}}_k(\lambda_k^{(N)})\underline{a})_{x_i} = 0$ reads

$$a_{x_i} = \sum_{\substack{1 \leq j \leq k \\ j \neq i}} \frac{W^{(k)}_{x_i,x_j}(\lambda_k^{(N)})}{\lambda_k^{(N)}} a_{x_j}.$$

The above expression and the normalization assumption imply

$$|a_{x_i}| \leq \sqrt{k} e^{(-\delta/9)\sqrt{N}}. \tag{6.29}$$

Since $a_{x_k}^2 = 1 - \sum_{i=1}^{k-1} |a_{x_i}|^2$, we get

$$1 \geq a_{x_k} \geq 1 - k^{3/2} e^{(-\delta/9)\sqrt{N}}. \tag{6.30}$$

Estimates (6.29) and (6.30) together with Lemma 8 imply (6.17) and (6.18).

To prove (6.16), let $\underline{a}$ be defined as above. Then the identity $(\hat{\mathcal{E}}_k(\lambda_k^{(N)})\underline{a})_{x_k} = 0$ reads

$$\left(1 - \frac{K^{(k)}_{x_k,x_k}}{\lambda_k^{(N)}}\right) = \sum_{j=1}^{k-1} \frac{W^{(k)}_{x_k,x_j}(\lambda_k^{(N)}) a_{x_j}}{\lambda_k^{(N)} a_{x_k}}.$$

By the Schwarz inequality, due to (6.29) and (6.30),

$$\left|1 - \frac{K^{(k)}_{x_k,x_k}}{\lambda_k^{(N)}}\right| \leq cq e^{(-\delta/10)\sqrt{N}}.$$

The above estimate together with (6.10) implies (6.16).

The last estimate (6.19) follows by straightforward computations from (6.17), (6.29), (6.30) and (5.16). □

**7. Subdiffusive behavior (proofs of Theorem 3 and Theorem 4).** We begin with the proof of Theorem 3.

Given a path $\gamma \in C(\mathbb{R})$, denote by $(m_1(\gamma), m(\gamma), m_2(\gamma))$ the consecutive 1-extrema (disregarding equivalent points) of the 1-valley of $\gamma$ covering the origin (if existing), namely,

$$\begin{aligned} m_1(\gamma) &\equiv \max\{x : x < 0 \text{ and } x \in M_1^+(\gamma)\},\\ m_2(\gamma) &\equiv \min\{x : x \geq 0 \text{ and } x \in M_1^+(\gamma)\},\\ \{m(\gamma)\} &\equiv (m_1(\gamma), m_2(\gamma)) \cap M_1^-(\gamma). \end{aligned}$$

In particular, the $\ln n$-extrema of the $\ln n$-valley of $V^{(1)}$ covering the origin can be written as

$$\begin{aligned} a^{(n)}(\omega) &\equiv \max\{x : x < 0 \text{ and } x \in M_{\ln n}^+(V^{(1)})\} = m_1(V^{(\ln^2 n)}(\omega)) \ln^2 n,\\ b^{(n)}(\omega) &\equiv \min\{x : x \geq 0 \text{ and } x \in M_{\ln n}^+(V^{(1)})\} = m_2(V^{(\ln^2 n)}(\omega)) \ln^2 n,\\ m^{(n)}(\omega) &\equiv m(V^{(\ln^2 n)}(\omega)) \ln^2 n = \mathfrak{m}^{(n)}(\omega) \ln^2 n. \end{aligned}$$

Note that the above quantities are defined $\mathbf{P}$ a.s. since $\limsup_{x\to\pm\infty} V(x) = \infty$ and $\liminf_{x\to\pm\infty} V(x) = -\infty$ $\mathbf{P}$ a.s.

Given $0 < \beta, \delta, \delta' < 1$, we denote by $\mathcal{B}_{\beta,\delta,\delta'}$ the set of paths $\gamma \in C(\mathbb{R})$ such that the following properties hold [for $m_1 = m_1(\gamma)$, $m_2 = m_2(\gamma)$, $m = m(\gamma)$]:

$$-m_1, m_2 \leq 1/\delta', \tag{7.1}$$

$$M_1^-(\gamma) \cap [-1/\delta', m) \neq \varnothing, \qquad M_1^-(\gamma) \cap (m, 1/\delta') \neq \varnothing, \tag{7.2}$$

$$M_{1-\delta}^-(\gamma) \cap [-1/\delta', 1/\delta'] = M_{1+\delta}^-(\gamma) \cap [-1/\delta', 1/\delta'], \tag{7.3}$$

$$\gamma(m_1) \wedge \gamma(m_2) \geq \max_{[0\wedge m, 0\vee m]} V + \delta, \tag{7.4}$$

$$\gamma(m) \geq -1/\delta, \tag{7.5}$$

$$\gamma(m_1) \wedge \gamma(m_2) \geq \max_{|x-m|\leq\beta} \gamma + \delta. \tag{7.6}$$

Due to the properties of Brownian motion and by means of the results of Section 2, one can show that there exist $\beta, \delta, \delta', n_0$ such that the set

$$\bar{\Omega}_n \equiv \{\omega \in \Omega : V^{(\ln^2 n)}(\omega) \in \mathcal{B}_{\beta,\delta,\delta'}\}$$

has probability $\mathbf{P}(\bar{\Omega}_n) \geq 1 - \alpha/2$ if $n \geq n_0$.

Fix $\beta, \delta, \delta'$ as above and set $N \equiv \ln^2 n$. Let $P^{(N)}$ and $C_1, C_2, C_3$ be as in Proposition 3 with the interval $[-1,1]$ replaced by $[-1/\delta', 1/\delta']$. Set

$$\varepsilon_n \equiv \frac{C_1 \ln N}{\sqrt{N}}, \qquad \delta_n \equiv \rho(\varepsilon_n),$$

where $\rho$ is defined as in Theorem 3. Note that $\delta_n$ differs from its definition in Theorem 3 by the factor $C_1$. It is simple to check that this is not restrictive.

Consider on the enlarged probability space with measure $P^{(N)}$ the event $\mathcal{C}^{(N)}_{\delta,\delta'}$ that the following conditions are satisfied, where $\bar{B} = (B(x) : x \in [-1/\delta', 1/\delta'])$ and $C_4$ is a fixed positive constant with $C_4 > 4$ and $1 + 2C_1 - C_1 C_4 < 0$:

$$\sup_{|x|\le 1/\delta'} |V^{(N)}(x) - B(x)| \le \varepsilon_n, \tag{7.7}$$

$$M^-_{1+C_4\varepsilon_n}(\bar{B}) = M^-_{1-C_4\varepsilon_n}(\bar{B}), \qquad M^+_{1+C_4\varepsilon_n}(\bar{B}) = M^+_{1-C_4\varepsilon_n}(\bar{B}), \tag{7.8}$$

$$\inf_{\delta_n/2 \le s \le T^{(1)}_{k,+}} |B_{S^{(h)}_k + s} - B_{S^{(h)}_k}| \ge C_4 \varepsilon_n \quad \text{if } k \in \mathbb{Z} \text{ and } S^{(h)}_k \in [-1/\delta', 1/\delta'], \tag{7.9}$$

$$\inf_{\delta_n/2 \le s \le T^{(1)}_{k,-}} |B_{S^{(h)}_k} - B_{S^{(h)}_k - s}| \ge C_4 \varepsilon_n \quad \text{if } k \in \mathbb{Z} \text{ and } S^{(h)}_k \in [-1/\delta', 1/\delta'], \tag{7.10}$$

$$\inf_{|s| \le \delta_n/2} |B_{S^{(h)}_k + s} - B_{S^{(h)}_k}| \le 1/2 \quad \text{if } k \in \mathbb{Z} \text{ and } S^{(h)}_k \in [-1/\delta', 1/\delta']. \tag{7.11}$$

Recall that $\{S^{(h)}_k\}_{k\in\mathbb{Z}}$ is the set of $h$-extrema of $B$, while the random times $T^{(h)}_{k,\pm}$ have been defined in Corollary 1. We stress that the above conditions (7.7)–(7.11) will be used only for proving (7.30) below, which could be done more efficiently by a direct analysis of the potential $V$ around the minimum $m$ as in [12] (see Remark 4).

By means of the results of Section 2, one can check that

$$\lim_{N\uparrow\infty} P^{(N)}(\mathcal{C}^{(N)}_{\delta,\delta'}) = 1.$$

We point out that in order to estimate the probability of the events (7.9) and (7.10) one has to use Lemma 3 [see also the proof of (2.13) in order to treat the case $n = 0$] together with the property that $\lim_{n\uparrow\infty} \varepsilon_n/\sqrt{\delta_n} = 0$.

Let $\Omega_n$ be the event in the enlarged probability space given by

$$\Omega_n \equiv \mathcal{C}^{(N)}_{\delta,\delta'} \cap \{\bar{V}^{(N)} \in \mathcal{B}_{\beta,\delta,\delta'}\}. \tag{7.12}$$

Then $P^{(N)}(\Omega_n) \ge 1 - \alpha$ if $n$ is large enough. In what follows, we will assume that the event $\Omega_n$ is realized.

Let us set

$$\begin{aligned} A_n &\equiv (a^{(n)}, b^{(n)}) \cap \mathbb{Z}, \\ D_n &\equiv ((\mathfrak{m}^{(n)} - \delta_n)\ln^2 n, (\mathfrak{m}^{(n)} + \delta_n)\ln^2) \cap A_n. \end{aligned} \tag{7.13}$$

Recall that $\mathbb{L}(A_n)$ is defined as $\mathbb{L}(A_n) = (\mathbb{L}_{x,y})_{x,y\in A_n}$. We write $\mathbb{P}(A_n)$ for the restriction of the jump probability matrix to $A_n \times A_n$, namely, $\mathbb{P}(A_n) = \mathbb{I} - \mathbb{L}(A_n)$. Then

$$\begin{aligned} &\mathrm{P}_0^\omega\left(\left|\frac{X_n}{\ln^2 n} - \mathfrak{m}^{(n)}\right| \le \delta_n\right) \\ &\quad = \mathrm{P}_0^\omega(X_n \in D_n) \ge \mathrm{P}_0^\omega(X_n \in D_n, X_k \in A_n \forall 0 \le k \le n) \\ &\quad = \sum_{y\in D_n} (\mathbb{P}(A_n)^n)_{0,y} = \frac{1}{\mu(0)}(1_0, \mathbb{P}(A_n)^n 1_{D_n}), \end{aligned} \tag{7.14}$$

where, in general, $1_Y$ denotes the characteristic function of the set $Y$ and $(\cdot,\cdot)$ denotes the scalar product in $\mathbb{L}^2(A_n,\mu)$ (the related norm will be denoted by $\|\cdot\|$).

By the same arguments of Section 5.5, due to (7.1) and (7.3), we obtain that the principal eigenvalue $\lambda_1^{(n)}$ of $\mathbb{L}(A_n)$ is a simple eigenvalue satisfying

$$c'(\ln n)^{-4} n^{-1-\delta} \le \lambda_1^{(n)} \le c n^{-1-\delta}. \tag{7.15}$$

Moreover, defining the function $h^\lambda$ on $A_n$ as

$$h^\lambda(y) \equiv h^\lambda_{m^{(n)},A_n^c}(y) \qquad \forall y \in A_n,$$

and setting $h \equiv h^0$, a principal eigenvector of $\mathbb{L}(A_n)$ is given by $h^{\lambda_1^{(n)}}$ and (see Proposition 5)

$$h(x) \le h^{\lambda_1^{(n)}}(x) \le h(x)(1 + p(\ln n) n^{-\delta}) \qquad \forall x \in A_n,$$

where, here and in what follows, $p$ denotes a generic polynomial having positive coefficients. In particular, the eigenvector $\psi_1^{(n)}$ obtained by normalizing $h^{\lambda_1^{(n)}}$, that is, $\psi_1^{(n)} \equiv h^{\lambda_1^{(n)}}/\|h^{\lambda_1^{(n)}}\|$, satisfies

$$\begin{aligned} &\left|\psi_1^{(n)}(x) - \frac{h(x)}{\|h\|}\right| \le \frac{h(x)}{\|h\|} p(\ln n) n^{-\delta} \qquad \forall x \in A_n, \\ &\left\|\psi_1^{(n)} - \frac{h}{\|h\|}\right\| \le p(\ln n) n^{-\delta}. \end{aligned} \tag{7.16}$$

We denote by $\lambda_2^{(n)} < \lambda_3^{(n)} < \cdots < \lambda_{|A_n|}^{(n)}$ the remaining (simple) eigenvalues of $\mathbb{L}(A_n)$ and by $\psi_2^{(n)}, \psi_3^{(n)}, \dots, \psi_{|A_n|}^{(n)}$ the related normalized eigenvectors. Due to (7.1), (7.3) and Theorem 1, $\lambda_2^{(n)}$ can be bounded from below as

$$\lambda_2^{(n)} \ge c(\ln n)^{-4} n^{-1+\delta}. \tag{7.17}$$

Since $\mathbb{P}(A_n)$ has simple eigenvalues given by

$$1 - \lambda_1^{(n)} > 1 - \lambda_2^{(n)} > \cdots > 1 - \lambda_{|A_n|}^{(n)},$$

with related eigenvectors $\psi_1^{(n)}, \psi_2^{(n)}, \dots, \psi_{|A_n|}^{(n)}$, we can write

$$\frac{1}{\mu(0)}(1_0, \mathbb{P}(A_n)^n 1_{D_n}) = \sum_{j=1}^{|A_n|}(1-\lambda_j^{(n)})^n(\psi_j^{(n)}, 1_{D_n})\psi_j^{(n)}(0). \tag{7.18}$$

Let $\Pi$ be the orthogonal projection of $L^2(A_n, \mu)$ along the subspace generated by $\psi_k^{(n)}$ with $2 \le k \le |A_n|-1$. Since, by Lemma 1,

$$\sup_{1<j<|A_n|} |1-\lambda_j^{(n)}| = 1-\lambda_2^{(n)},$$

we obtain the bound

$$\begin{aligned} &\left|\sum_{j=2}^{|A_n|-1}(1-\lambda_j^{(n)})^n(\psi_j^{(n)}, 1_{D_n})\psi_j^{(n)}(0)\right| \\ &\quad = \left|\frac{1}{\mu(0)}(1_0, \mathbb{P}(A_n)^n \Pi 1_{D_n})\right| \\ &\quad \le c(\kappa)(1-\lambda_2^{(n)})^n \|1_{D_n}\|. \end{aligned} \tag{7.19}$$

We claim that

$$\lim_{n\uparrow\infty} \sup_{\Omega_n} |(1-\lambda_1^{(n)})^n(\psi_1^{(n)}, 1_{D_n})\psi_1^{(n)}(0) - 1| = 0, \tag{7.20}$$

$$\lim_{n\uparrow\infty} \sup_{\Omega_n} |(1-\lambda_{|A_n|}^{(n)})^n(\psi_{|A_n|}^{(n)}, 1_{D_n})\psi_{|A_n|}^{(n)}(0)| = 0, \tag{7.21}$$

$$\lim_{n\uparrow\infty} \sup_{\Omega_n} (1-\lambda_2^{(n)})^n \|1_{D_n}\| = 0. \tag{7.22}$$

Note that the above estimates together with (7.14) and (7.18) imply Theorem 3.

Let us prove (7.20). Due to (7.15),

$$\lim_{n\uparrow\infty} \sup_{\omega\in\Omega_n} |(1-\lambda_1^{(n)})^n - 1| = 0, \tag{7.23}$$

while, due to (7.16),

$$\left|(\psi_1^{(n)}, 1_{D_n})\psi_1^{(n)}(0) - \left(\frac{h}{\|h\|}, 1_{D_n}\right)\frac{h(0)}{\|h\|}\right| \le p(\ln n) n^{-\delta}\frac{\|1_{D_n}\|}{\|h\|}h(0).$$

Applying Lemma 9 completes the proof.

To prove (7.22), we observe that, due to (7.1) and (7.5),

$$\|1_{D_n}\|^2 \le c \ln^2 n \cdot \exp\{-V(m^{(n)})\} \le c' \ln^2 n \cdot n^{1/\delta'}. \tag{7.24}$$

The above estimate together with (7.17) implies (7.22).

Finally, note that due to (7.20) and (7.22), it is clear that

$$\limsup_{n\uparrow\infty}\sup_{\Omega_n}(1-\lambda^{(n)}_{|A_n|})^n(\psi^{(n)}_{|A_n|},1_{D_n})\psi^{(n)}_{|A_n|}(0)=0.$$

But, on the other hand, since $1-\lambda_{|A_n|}<0$, and all quantities vary slowly with $n$,

$$(1-\lambda^{(n)}_{|A_n|})^n(\psi^{(n)}_{|A_n|},1_{D_n})\sim-(1-\lambda^{(n+1)}_{|A_{n+1}|})^{n+1}(\psi^{(n+1)}_{|A_{n+1}|},1_{D_{n+1}}),$$

implying that

$$\limsup_{n\uparrow\infty}\sup_{\Omega_n}(1-\lambda^{(n)}_{|A_n|})^n(\psi^{(n)}_{|A_n|},1_{D_n})=-\liminf_{n\uparrow\infty}\inf_{\Omega_n}(1-\lambda^{(n)}_{|A_n|})^n(\psi^{(n)}_{|A_n|},1_{D_n}),$$

which yields (7.21).

LEMMA 9.

$$\lim_{n\uparrow\infty}\sup_{\Omega_n}|h(0)-1|=0, \tag{7.25}$$

$$\lim_{n\uparrow\infty}\sup_{\Omega_n}\left|\frac{(h,1_{D_n})}{\|h\|^2}-1\right|=0, \tag{7.26}$$

$$\lim_{n\uparrow\infty}\sup_{\Omega_n}\left|\frac{\|1_{D_n}\|}{\|h\|}-1\right|=0. \tag{7.27}$$

PROOF. Let us suppose that the event $\Omega_n$ is verified. In order to prove (7.25), suppose, for example, that $a^{(n)}\le 0<m^{(n)}$, thus implying that $1-h(0)=h_{a^{(n)},m^{(n)}}(0)$. Therefore, due to (3.10) and assumptions (7.1) and (7.4),

$$1-h(0)\le c\ln^2 n\cdot\exp\Big\{\max_{[0,m^{(n)}-1]}V-\max_{[a^{(n)},m^{(n)}-1]}V\Big\}\le c\ln^2 n\cdot n^{-\delta},$$

thus implying (7.25).

We prove now (7.26). The proof of (7.27) is similar and we will omit it. Let us first bound $1-h(x)$ for $x\in D_n$. Suppose, for example, that $x<m^{(n)}$, thus implying $1-h(x)=h_{a^{(n)},m^{(n)}}(x)$. Due to (3.10), (7.1) and (7.6),

$$1-h(x)\le c\ln^2 n\cdot\exp\Big\{\max_{[x,m^{(n)}-1]}V-\max_{[a^{(n)},m^{(n)}-1]}V\Big\}\le c\ln^2 n\cdot n^{-\delta}. \tag{7.28}$$

In particular,

$$\left|\sum_{x\in D_n}\mu(x)h^2(x)-\sum_{x\in D_n}\mu(x)h(x)\right|\le c\ln^2 n\cdot n^{-\delta}\sum_{x\in D_n}\mu(x)h(x). \tag{7.29}$$

Let us write

$$\frac{(h, 1_{D_n})}{\|h\|^2} = \frac{W_1(n)}{W_2(n)} \cdot \frac{W_2(n)}{W_3(n)},$$

where

$$W_1(n) \equiv \sum_{x \in D_n} \frac{\mu(x)}{\mu(m^{(n)})} h(x),$$

$$W_2(n) \equiv \sum_{x \in D_n} \frac{\mu(x)}{\mu(m^{(n)})} h(x)^2,$$

$$W_3(n) \equiv \sum_{x \in A_n} \frac{\mu(x)}{\mu(m^{(n)})} h(x)^2.$$

Then, due to (7.29), $\lim_{n\uparrow\infty} \sup_{\Omega_n} |W_1(n)/W_2(n) - 1| = 0$. Since $W_3(n) \geq 1$, in order to prove that $\lim_{n\uparrow\infty} \sup_{\Omega_n} |W_2(n)/W_3(n) - 1| = 0$, it is enough to show

$$\lim_{n\uparrow\infty} \sup_{\Omega_n} \sum_{x \in A_n \setminus D_n} e^{-(V(x) - V(m^{(n)}))} = 0. \tag{7.30}$$

To this aim, it is more convenient to work on the rescaled lattice $\mathbb{Z}/N$, where $N = \ln^2 n$, and compare $V^{(N)}$ with $B$. Let us set here $m \equiv m(V^{(N)})$, $m_1 \equiv m_1(V^{(N)})$ and $m_2 \equiv m_2(V^{(N)})$. We can write

$$\sum_{x \in A_n \setminus D_n} e^{-(V(x) - V(m^{(n)}))} = \sum_{x \in A_n/N \setminus D_n/N} e^{-\sqrt{N}(V^{(N)}(x) - V^{(N)}(m))}. \tag{7.31}$$

Due to Lemma 5 applied with $\varepsilon = 4\varepsilon_n$ and with $[-1, 1]$ substituted with $[-1/\delta', 1/\delta']$, and due to the definition of $\Omega_n$, there exists a 1-minimum $m^*$ of

$$\bar{B} = (B(x) : x \in [-1/\delta', 1/\delta'])$$

such that

$$|B(m^*) - V(m)| \leq \varepsilon_n, \qquad |m - m^*| \leq \delta_n/2.$$

Let us denote $m_1^*$, $m_2^*$ the first 1-maximum of $\bar{B}$ respectively on the left and on the right of $m^*$. Due to Lemma 5 [see, in particular, (2.17)] and the definition of $\Omega_n$,

$$|m_1 - m_1^*| \leq \delta_n/2, \qquad |m_2 - m_2^*| \leq \delta_n/2.$$

In particular,

$$\text{r.h.s. of (7.31)} \leq \sum_{x \in \Delta_1 \cup \Delta_2} e^{2C_1 \ln N} e^{-\sqrt{N}(B(x) - B(m^*))}, \tag{7.32}$$

where

$$\Delta_1 \equiv (m_1^* - \delta_n/2, m^* - \delta_n/2] \cap \mathbb{Z}/N, \qquad \Delta_2 \equiv [m^* + \delta_n/2, m_2^* + \delta_n/2) \cap \mathbb{Z}/N.$$

Let us estimate the contribution in (7.32) of the addenda $x \in \Delta_1$ (the case $x \in \Delta_2$ can be treated similarly).

Due to (7.2) and Lemma 5, $m_1^*$, $m_2^*$ are 1-maxima of $B$. Hence, there exists $k \in \mathbb{Z}$ such that

$$m_1^* = S_{k-1}^{(1)}, \qquad m^* = S_k^{(1)}, \qquad m_2^* = S_{k+1}^{(1)}.$$

We can write

$$\sum_{x \in \Delta_1} e^{2C_1 \ln N} e^{-\sqrt{N}(B(x) - B(m^*))} = I_1 + I_2 + I_3,$$

where $I_i$ is given by the sum over $x \in R_i$ and

$$R_1 = (m_1^* - \delta_n, m_1^*] \cap \mathbb{Z}/N,$$
$$R_2 = (m_1^*, m^* - T_{k,-}^{(1)}],$$
$$R_3 = (m^* - T_{k,-}^{(1)}, m^* - \delta_n).$$

Consider the case $x \in R_1$. Due to (7.8) and (7.11),

$$\begin{aligned} B(x) - B(m^*) &= (B(x) - B(m_1^*)) + (B(m_1^*) - B(m^*)) \\ &\geq -1/2 + 1 + C_4\varepsilon_n \geq 1/2. \end{aligned}$$

Since, due to (7.1), $|A_n| \leq cN$, we get

$$I_1 \leq |A_n| e^{2C_1 \ln N - \sqrt{N}/2} \leq cN^{1+2C_1} e^{-\sqrt{N}/2}. \tag{7.33}$$

Consider the case $x \in R_2$. Due to (7.8), $B(x) - B(m^*) \geq C_4\varepsilon_n$ [otherwise between $m_1^*$ and $m^*$ there would be a $(1 - C_4\varepsilon_n)$-minimum in contradiction with (7.8)]. Hence,

$$I_2 \leq |A_n| N^{2C_1} e^{-C_4\varepsilon_n\sqrt{N}} \leq cN^{1+2C_1-C_1C_4}. \tag{7.34}$$

Consider the case $x \in R_3$. Due to (7.10), $B(x) - B(m^*) \geq C_4\varepsilon_n$. Hence,

$$I_3 \leq |A_n| N^{2C_1} e^{-C_4\varepsilon_n\sqrt{N}} \leq cN^{1+2C_1-C_1C_4}. \tag{7.35}$$

By collecting together (7.33), (7.34) and (7.35), we get

$$\sum_{x \in \Delta_1} e^{2C_1 \ln N} e^{-\sqrt{N}(B(x) - B(m^*))} \leq cN^{1+2C_1} e^{-\sqrt{N}/2} + cN^{1+2C_1-C_1C_4}.$$

Since, by assumption, $1 + 2C_1 - C_1C_4$, the r.h.s. goes to 0 as $n$ goes to $\infty$. □

Let us conclude this section with some remarks, and the proof of Theorem 4. We will throughout the remainder of this discussion assume without further mentioning that the random environment is such that the hypothesis of our main statements are verified for all Dirichlet operators we will consider. The reader can check that this holds with high probability (see also the proof of Theorem 4 where the technical steps are discussed in more detail).

First, we note that the choice of the set $A_n$ in the lower bound (7.14), although probabilistically justified by the fact that the process will not have left $A_n$ by time $n$ and will not have remained in a much smaller set, either, with high probability, seems awkward from a spectral point of view. In fact, we should obtain the same localization result if we choose instead of $A_n$ a much larger set. To see this in some detail, let us consider any interval $A \supset A_n$. Obviously, we have that

$$\mathrm{P}_0^\omega(X_n \in D) \geq \sum_{j=1}^{|A|} (1-\lambda_j^{(A)})^n (\psi_j^{(A)}, 1_D)\psi_j^{(A)}(0) \qquad \forall D \subset A, \tag{7.36}$$

where $\lambda_j^{(A)}$, $\psi_j^{(A)}$ are the eigenvalues and eigenfunctions of $\mathbb{L}(A)$, with $\lambda_j^{(A)}$ increasingly ordered.

Let us understand what we can say about the spectrum of $\mathbb{L}(A)$. Let us write $\lambda_1^{(n)}$ for the principal eigenvalue of $\mathbb{L}(A_n)$. We know that $\lambda_1^{(A)} \leq \lambda_1^{(n)}$. Let $k$ be the number of eigenvalues of $\mathbb{L}(A)$ which are smaller or equal than $\lambda_1^{(n)}$. If $k = 1$, then the analysis above remains essentially unchanged. In what follows we suppose $k \geq 2$.

From our analysis of eigenvalues, this means that the potential $V^{(1)}$ [recall the definition (1.9)] restricted to $A$ has $k$ $O(\ln n)$-minima (we always assume $n$ large). Let us denote these minima by $x_1, \dots, x_k$, labeled as in Section 5 to correspond to increasing eigenvalues of $\mathbb{L}(A)$. Clearly, one of these minima is $m^{(n)}$, defined as in Theorem 3, say, $x_l = m^{(n)}$. Let us denote by $B_i$ small neighborhoods of the minima $x_i$.

Using the same arguments as before, we see that we get, up to terms tending to zero with $n$,

$$\mathrm{P}_0^\omega(X_n \in B_i) \geq \sum_{j=1}^{k} (\psi_j^{(A)}, 1_{B_i})\psi_j^{(A)}(0). \tag{7.37}$$

Now we know that the left-hand side of the equation equals one, if $i = l$, and zero, otherwise. On the other hand, we also know, from our estimate on the eigenfunctions, that

$$(\psi_j^{(A)}, 1_{B_j})\psi_j^{(A)}(0) \sim \frac{\mu(B_j)}{\|h_{x_j,S_{j-1}^*}\|_2^2} h_{x_j,S_{j-1}^*}(0) \sim h_{x_j,S_{j-1}^*}(0), \tag{7.38}$$

where $S^*_{j-1} = \{x_1, x_2, \ldots, x_{j-1}\} \cup (\mathbb{Z} \setminus A)$. Note that the right-hand side is essentially one, if "0 is in the valley of $x_j$," where we call the valley of $x_j$ the interval between the two highest maxima to the right and to the left of $x_j$ one needs to cross to reach $S^*_{j-1}$ from $x_j$.

Let us now look at the probability to be in $B_l$. Up to terms tending to zero with $n$, we can write this as

$$\begin{aligned} \mathrm{P}_0^\omega(X_n \in B_l) \geq{}& (\psi_l^{(A)}, 1_{B_l})\psi_l^{(A)}(0) \\ &+ \sum_{j:\mu(x_j)>\mu(x_l)} (\psi_j^{(A)}, 1_{B_l})\psi_j^{(A)}(0) \\ &+ \sum_{j:\mu(x_j)<\mu(x_l)} (\psi_j^{(A)}, 1_{B_l})\psi_j^{(A)}(0). \end{aligned} \tag{7.39}$$

We already know that the first term equals one, as does the left-hand side. Now for the first sum we get an easy asymptotic bound using again our estimates for the eigenfunctions, namely (up to terms tending to zero with $n$),

$$\sum_{j:\mu(x_j)>\mu(x_l)} (\psi_j^{(A)}, 1_{B_l})\psi_j^{(A)}(0) \leq \sum_{j:\mu(x_j)>\mu(x_l)} \frac{\mu(B_l)}{\mu(B_j)}, \tag{7.40}$$

which will tend to zero with $n$ (in the good subspace of environments). To deal with the second sum, we need to be more careful. First, note that $x_l$ is the minimum of the $\ln n$-valley that contains 0; thus, it is not possible that any of the valleys of the $x_j$ with $V(x_j) > V(x_l)$ contains the origin. Using these facts, and the precise representation of the eigenfunction (6.17), pointwise estimates on the $h^\lambda$ (see Lemma 4.3 of [3]), and the usual estimates on the equilibrium potential, one may show that indeed all terms in this sum also tend to zero with $n$. We leave the details for the interested reader.

A more interesting observation ensues when regarding a neighborhood $B_i$ with $\mu(x_i) > \mu(x_l)$ and such that 0 is contained in the valley of $x_i$. Then we know that, up to terms tending to zero with $n$,

$$\begin{aligned} o(1) &= \mathrm{P}_0^\omega(X_n \in B_i) \\ &\geq (\psi_i^{(A)}, 1_{B_i})\psi_i^{(A)}(0) \\ &\quad + \sum_{j:\mu(x_j)>\mu(x_i)} (\psi_j^{(A)}, 1_{B_i})\psi_j^{(A)}(0) \\ &\quad + \sum_{j:\mu(x_j)<\mu(x_i)} (\psi_j^{(A)}, 1_{B_i})\psi_j^{(A)}(0). \end{aligned} \tag{7.41}$$

Now the first term in the r.h.s. is close to one, while the first sum, by the same estimates as before, tends to zero. Thus, we can conclude that

$$\sum_{j:\mu(x_j)<\mu(x_i)} (\psi_j^{(A)}, 1_{B_i})\psi_j^{(A)}(0) \sim -1. \tag{7.42}$$

We see that the small negative parts of the eigenfunctions play a crucial role here and cannot be neglected! Deriving (7.42) directly from our estimates on the eigenfunctions is not possible.

Proof of Theorem 4. Let us now exploit these observations to prove Theorem 4. To do this, we recall the construction of the sequence of boxes $A_{n_k}$ [recall the definition (7.13) of $A_n$ and $D_n$]: start with $A_{n_0}$, $n_0$ large. Then increase $n$ to $n_1$ such that, for the first time, $m^{(n_1)} \neq m^{(n_0)}$, and so on.

Let

$$\lambda_1^{(n_k)} < \lambda_2^{(n_k)} < \cdots < \lambda_{|A_{n_k}|-1}^{(n_k)} < \lambda_{|A_{n_k}|}^{(n_k)}$$

be the eigenvalues of the generator $\mathbb{L}(A_{n_k})$ and call $\psi_j^{(n_k)}$ the eigenvector associated to $\lambda_j^{(n_k)}$. Due to Lemma 1,

$$1 - \lambda_i^{(n_k)} = -(1 - \lambda_{|A_{n_k}|-i+1}^{(n_k)}) \qquad \forall 1 \le i \le |A_{n_k}|, \tag{7.43}$$

and we can assume that $\psi_j^{(n_k)}$ and $\psi_{|A_{n_k}|-i+1}^{(n_k)}$ coincide on even sites and are opposite on odd sites.

Given positive constants $\delta, \delta', \beta$, we define $\Omega_n = \Omega_n(\delta, \delta', \beta) \subset \Omega$ as in the proof of Theorem 3 [see (7.12)] with the following additional assumption: let $h_1 > h_2 > h_3$ be the minimal values such that

$$|M_{h_i}^-(V^{(\ln^2 n_k)}) \cap (m_1, m_2)| = i, \qquad 1 \le i \le 3$$

(it is understood that we exclude degenerate cases), then we require that

$$h_1 \ge h_2 + \delta, \qquad h_2 \ge h_3 + \delta. \tag{7.44}$$

Given $\alpha > 0$, we fix $\delta, \delta', \beta > 0$ such that $\mathbf{P}(\Omega_n) \ge 1 - \alpha$. Note that, due to our choice of $\Omega_n$, all results obtained in the proof of Theorem 3 remain valid.

For each $T_k \in \mathbb{N}$, we can write

$$P_0^\omega(X_{T_k} \in D_{n_k}) = I_k + E_k, \tag{7.45}$$

where

$$I_k = \sum_{j=1}^{|A_{n_k}|} I_{k,j}, \qquad I_{k,j} = (1 - \lambda_j^{(n_k)})^{T_k} (\psi_j^{(n_k)}, 1_{D_{n_k}})\psi_j^{(n_k)}(0), \tag{7.46}$$

and the error term $E_k$ is bounded by

$$0 \le E_k \le P_0^\omega(\tau_{A_{n_k}^c} < T_k). \tag{7.47}$$

LEMMA 10. *Assume that* $T_k \lambda_1^{(n_k)} = o(1)$. *Then*

$$P_0^\omega(\tau_{A_{n_k}^c} < T_k) = o(1), \tag{7.48}$$

*where* $o(1)$ *is uniform for all* $\omega \in \Omega_{n_k}$.

PROOF. The proof relies on the representation that is completely analogous to the second line in (7.14):

$$P_0^\omega(\tau_{A_{n_k}^c} \ge T_k) = \frac{1}{\mu(0)}(1_0, \mathbb{P}(A_{n_k})^{T_k - 1} 1_{A_{n_k}}). \tag{7.49}$$

To show that this quantity goes to one [and, hence, the expression in equation (7.48) goes to zero], we proceed exactly as in the proof of Theorem 3. It is clear that we have succeeded if we establish the analogues of relations (7.20), (7.21) and (7.22), with $D_n$ replaced by $A_{n_k}$ and $n$ replaced by $n_k$, with exception of the exponent of $(1 - \lambda_j^{(n_k)})$, which becomes $T_k - 1$.

There is, however, nothing to be done: the first two relations follow without change, and for the last it is enough to notice that the bound on $\|1_{A_{n_k}}\|$ is essentially equal to that of $\|1_{D_{n_k}}\|$, since $D_{n_k}$ is a neighborhood of the deepest minimum on $A_{n_k}$. The assertion of the proof is thus obvious. □

Due to (7.23) and (7.43), the limits (7.20) and (7.22) imply

$$\lim_{k\uparrow\infty} \sup_{\omega\in\Omega_{n_k}} |(\psi_1^{(n_k)}, 1_{D_{n_k}})\psi_1^{(n_k)}(0) - 1| = 0, \tag{7.50}$$

$$\lim_{k\uparrow\infty} \sup_{\omega\in\Omega_{n_k}} |(\psi_{|A_{n_k}|}^{(n_k)}, 1_{D_{n_k}})\psi_{|A_{n_k}|}^{(n_k)}(0)| = 0. \tag{7.51}$$

Moreover, as done in (7.19) and (7.24), if $\omega \in \Omega_{n_k}$, we can always bound

$$\sum_{j=3}^{|A_{n_k}|-2} I_{k,j} \le |1 - \lambda_3^{(n_k)}|^{T_{n_k}} \|1_{D_{n_k}}\| \le c|1 - \lambda_3^{(n_k)}|^{T_k} \ln n_k \cdot n_k^{1/(2\delta)}. \tag{7.52}$$

We consider below two cases: $n_k^{h_3} \ll T_k \ll n_k^{h_2}$ and $T_k = t/\lambda_2^{(n_k)}$. Note that, since $\lambda_1^{(n_k)} \sim n_k^{-h_1}$, in both cases $T_k \lambda_1^{(n_k)} = o(1)$.

All our estimates below have to be considered uniform on $\omega \in \Omega_{n_k}$.

*Case* $n_k^{h_3} \ll T_k \ll n_k^{h_2}$. Take, for example, $T_k = n_k^{(h_2+h_3)/2}$.

Then [see (1.28)]

$$(1 - \lambda_i^{(n_k)})^{T_k} \sim e^{-T_k \lambda_i^{(n_k)}} \sim e^{-T_k / n_k^{h_i}}, \qquad i = 1, 2, 3. \tag{7.53}$$

In particular, for $i=1,2$, $(1-\lambda_i^{(n_k)})^{T_k}=1+o(1)$. Due to (7.50) and (7.51), $I_{k,1}=1+o(1)$ and $I_{k,|A_{n_k}|}=o(1)$. Moreover, (7.53) and (7.52) imply that $\sum_{j=3}^{|A_{n_k}|-2} I_{k,j}=o(1)$. Hence,

$$I_k = 1 + I_{k,2} + I_{k,|A_{n_k}|-1} + o(1).$$

By Lemma 10, also the error term $E_k$ goes to 0 uniformly in $\omega\in\Omega_{n_k}$. On the other hand, due to Theorem 3 and the geometry of $V$ over $A_{n_k}$,

$$P_0^\omega(X_{T_k}\in D_{n_k}) = o(1).$$

(7.45) and the above observations imply that

$$I_{k,2} + I_{k,|A_{n_k}|-1} = -1 + o(1).$$

Reasoning as in the proof of (7.21), from the above expression we derive that $I_{k,2}=-1+o(1)$, $I_{k,|A_{n_k}|-1}=o(1)$. Since $1-\lambda_2^{(n_k)}=1+o(1)$ and due to (7.43), the previous identities imply that

$$(\psi_2^{(n_k)}, 1_{D_{n_k}})\psi_2^{(n_k)}(0) = 1+o(1), \tag{7.54}$$

$$(\psi_{|A_{n_k}|-1}^{(n_k)}, 1_{D_{n_k}})\psi_{|A_{n_k}|-1}^{(n_k)}(0) = o(1). \tag{7.55}$$

*Case* $T_k = t/\lambda_2^{(n_k)}$. In this case

$$(1-\lambda_i^{(n_k)})^{T_k} \begin{cases} =1+o(1), & \text{if } i=1,\\ =e^{-t}+o(1), & \text{if } i=2,\\ \le \exp(-ctn_k^{h_2}/n_k^{h_3}), & \text{if } i=3. \end{cases} \tag{7.56}$$

By means of (7.56) with $i=1$, (7.50) and (7.51), we get that $I_{k,1}=1+o(1)$, $I_{k,|A_{n_k}|}=o(1)$. By means of (7.56) with $i=2$, (7.54) and (7.55), we get that $I_{k,2}=e^{-t}+o(1)$, $I_{k,|A_{n_k}|-1}=o(1)$. Moreover, (7.56) with $i=3$ and (7.52) imply that $\sum_{j=3}^{|A_{n_k}|-2} I_{k,j}=o(1)$. Hence,

$$I_k = 1 - e^{-t} + o(1).$$

The assertion now follows from (7.45) and the fact that $I_k=1-e^{-t}+o(1)$, $E_k=o(1)$. $\square$

Note that this observation suggests the following trap model caricature of Sinai's random walk: Take the sequence of values $\lambda_2^{(n_k)}\equiv\Lambda_k$; this sequence is fully determined by the random potential. Now consider the continuous time Markov chain on the positive integers that jumps from site $k$ to site $k+1$ with rate $\Lambda_k$.

## APPENDIX A: RG-ALGORITHM LABELING THE $H$-MINIMA

To compare our spectral results with [15], we show in this Appendix that the renormalization group algorithm of [10], Section II, leads to the labeling $M_h^-(\gamma) = \{x_1, x_2, \dots, x_q\}$, fulfilling (1.16), whenever $\gamma \in C([-1,1])$ satisfies $|M_h^-(\gamma)| = q \geq 1$ and

$$||\gamma(y) - \gamma(x)| - |\gamma(y') - \gamma(x')|| \geq \delta \qquad \forall (x,y) \neq (x', y') \in M_h^-(\gamma) \times M_h^+(\gamma). \tag{A.1}$$

It is simple to check that (A.1) implies (1.15). Let us first describe the RG-algorithm in terms of $h$-extrema. To this aim, we label the points of $M_h^-(\gamma) \cup M_h^+(\gamma)$ as

$$z_1^{(1)} < z_2^{(1)} < \cdots < z_{2q+1}^{(1)}.$$

As discussed in Lemma 2, $z_j^{(1)}$ is a $h$-maximum if $j$ is odd, otherwise it is a $h$-minimum. We introduce a coarse-grained potential $\mathcal{V}^{(1)}$ on $[-1,1]$ by setting

$$\mathcal{V}^{(1)}(x) = \begin{cases} -\infty, & \text{if } x \in \{-1,1\} \setminus \{z_1^{(1)}, z_{2q+1}^{(1)}\}, \\ V(x), & \text{if } x = z_i^{(1)},\ 1 \leq i \leq 2q+1, \end{cases}$$

and by extending $\mathcal{V}^{(1)}$ to all $[-1,1]$ by linear interpolation (see Figure 5). Note that $\mathcal{V}^{(1)} \equiv -\infty$ on $[-1,1] \setminus [z_1^{(1)}, z_{2q+1}^{(1)}]$.

We now define inductively by decimation of the less deep valley new potentials $\mathcal{V}^{(2)}, \mathcal{V}^{(3)}, \dots, \mathcal{V}^{(q)}$ on $[-1,1]$ satisfying the following property: For each $2 \leq i \leq q$, there exist $-1 \leq a_i < b_i \leq 1$ such that $\mathcal{V}^{(i)} \equiv -\infty$ on $[-1,1] \setminus [a_i, b_i]$ and $\mathcal{V}^{(i)}$ is piecewise-linear on $[a_i, b_i]$, with $a_i, b_i$ local maxima and having $q - i + 1$ local minima in $[a_i, b_i]$. To this aim, suppose $\mathcal{V}^{(i)}$ to be defined for some $1 \leq i \leq q-1$, fulfilling the above properties, and write

$$a_i = z_1^{(i)} < z_2^{(i)} < \cdots < z_{2(q-i+1)+1}^{(i)} = b_i$$

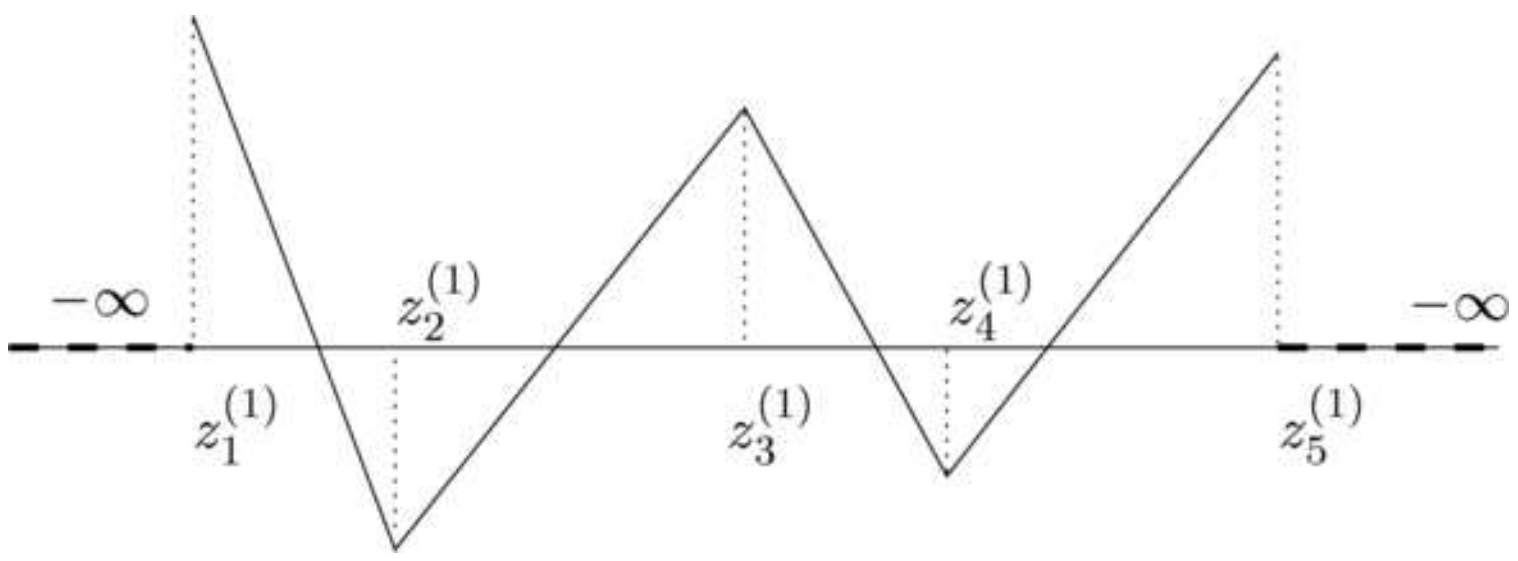

FIG. 5. *Potential* $\mathcal{V}^{(1)}$, $q = 2$.

for its $h$-extrema on $[a_i, b_i]$. Let us consider the bond $[z_k^{(i)}, z_{k+1}^{(i)}]$, $k = k(i)$, with the smallest variation of $\mathcal{V}^{(i)}$:

$$\begin{aligned} &|\mathcal{V}^{(i)}(z_k^{(i)}) - \mathcal{V}^{(i)}(z_{k+1}^{(i)})| \\ &\qquad = \min\{|\mathcal{V}^{(i)}(z_s^{(i)}) - \mathcal{V}^{(i)}(z_{s+1}^{(i)})| : 1 \le s \le 2(q-i+1)\}. \end{aligned} \tag{A.2}$$

Note that the index $k$ is uniquely defined due to (A.1).

Let us define $D_i \equiv \{z_1^{(i)}, z_2^{(i)}, \ldots, z_{2(q-i+1)+1}^{(i)}\}$ and $D_{i+1} \equiv D_i \setminus \{z_k^{(i)}, z_{k+1}^{(i)}\}$. Then $\mathcal{V}^{(i+1)}$ is defined by setting

$$\mathcal{V}^{(i+1)}(x) = \begin{cases} -\infty, & \text{if } x \in \{-1, 1\} \setminus D_{i+1}, \\ V(x), & \text{if } x \in D_{i+1}, \end{cases}$$

and by extending $\mathcal{V}^{(i+1)}$ to all $[-1, 1]$ by linear interpolation. In Figure 6 we consider the case $1 < k < 2(q - i + 1) + 1$.

Finally, we denote by $T_j$ the r.h.s. of (A.2) and by $y_i$ the local minimum of $V^{(i)}$ in $\{z_k^{(i)}, z_{k+1}^{(i)}\}$. Since for a given curve $\gamma$ they depend on $h$, we will sometimes write $y_i(h)$, $T_i(h)$ in order to underline this dependence. We can now state the relation between the above RG-construction and the labeling satisfying (1.16 ):

Proposition 7. *Let $h > 0$ and $\gamma \in C([-1, 1])$ satisfying $|M_h^-(\gamma)| = q \ge 1$ and (A.1). Moreover, let $\{x_1, x_2, \ldots, x_q\}$ be the labeling of $M_h^-(\gamma)$ satisfying (1.16) and let $y_1, y_2, \ldots, y_q$, $T_1, T_2, \ldots, T_q$ defined as in the above RG-construction. Then*

$$x_k = y_{q-k+1}, \qquad \gamma(z^*(x_k, S_{h,k-1})) - \gamma(x_k) = T_{q-k+1} \qquad \forall 1 \le k \le q,$$

*where $S_{h,k-1}(\gamma) = \{x_1, x_2, \ldots, x_{k-1}\} \cup \{-1, 1\}$.*

Proof. We prove the proposition by induction on $q$. It is simple to check that the assertion holds for all $h > 0$ if $q = 1$. Assume that it is valid for all $h > 0$ if $q = \bar{q} - 1$, for some $\bar{q} \ge 2$. We fix $\gamma \in C([-1, 1])$ and $h > 0$ such that $M_h^-(\gamma) = \bar q$. Let $y_1(h), \dots, y_{\bar q}(h)$, $T_1(h), \dots, T_{\bar q}(h)$ be defined by the RG-procedure described above. We observe that $M_{h'}^-(\gamma) = M_h^-(\gamma) \setminus \{y_1\}$, where $h' \equiv T_1(h) + \delta$ and $y_k(h') = y_{k+1}(h)$ for all $1 \le k \le \bar q - 1$. Setting

$$\begin{aligned} X_j &\equiv y_{\bar q - j}(h') = y_{\bar q - j + 1}(h), \\ S_{h',k} &= \{X_1, X_2, \dots, X_k\} \cup \{-1, 1\} \qquad \forall 1 \le j \le \bar q - 1, \end{aligned}$$

by the inductive hypothesis, we obtain that

$$\begin{aligned} &\gamma(z^*(X_k, S_{h',k-1})) - \gamma(X_k) \\ &\qquad \ge \max_{\bar q - 1 \ge j > k} \{\gamma(z^*(X_j, S_{h',k-1})) - \gamma(X_j)\} + \delta \\ &\qquad\qquad \forall 1 \le k \le \bar q - 2 \end{aligned} \tag{A.3}$$

and

$$\gamma(z^*(X_k, S_{h',k-1})) - \gamma(X_k) = T_{\bar q - k}(h') \qquad \forall 1 \le k \le \bar q - 1. \tag{A.4}$$

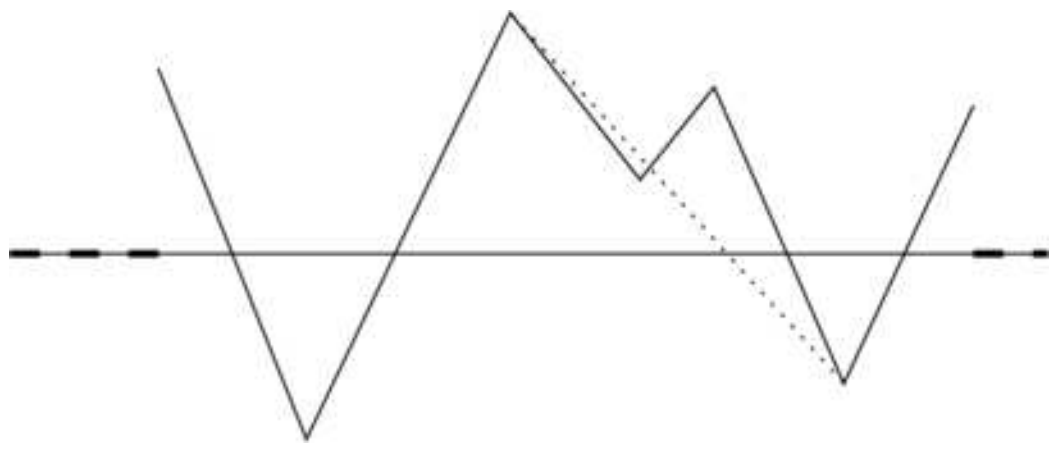

Fig. 6. *Decimation of the less deep valley.*

Let us now define $x_1 \equiv X_1, \dots, x_{\bar q - 1} \equiv X_{\bar q - 1}, x_{\bar q} \equiv y_1$. We claim that $\{x_1, \dots, x_{\bar q}\}$ satisfies (1.16). In fact, by observing that $S_{h,k} = S_{h',k}$ for $1 \le k \le \bar q - 1$, due to (A.3), we only need to prove

$$\gamma(z^*(x_k, S_{h,k-1})) - \gamma(x_k) \ge \gamma(z^*(y_1, S_{h,k-1})) - \gamma(y_1) + \delta \qquad \forall 1 \le k \le \bar q - 1.$$

The above inequalities follow easily from (A.1) and the fact that (1.15) implies

$$\begin{aligned} \gamma(z^*(x_k, S_{h,k-1})) - \gamma(x_k) &= \max_{x \in [-1,1] \setminus S_{h,k-1}} \gamma(z^*(x, S_{h,k-1})) - \gamma(x) \\ &\qquad \forall k : 1 \le k \le q - 1. \end{aligned} \tag{A.5}$$

To conclude the proof, we need to show that

$$\gamma(z^*(x_k, S_{h,k-1})) - \gamma(x_k) = T_{\bar q - k + 1}(h) \qquad \forall 1 \le k \le \bar q.$$

Since $T_{\bar q - k + 1}(h) = T_{\bar q - k}(h')$ for all $1 \le k < \bar q$ and due to (A.4), one only needs to check the trivial identity

$$\gamma(z^*(y_1, \{x_1, \dots, x_{\bar q - 1}\})) - \gamma(y_1) = T_1(h). \qquad \square$$

## APPENDIX B: PROOF OF LEMMA 4

Recall the definition of the random variable $X^{(h)}$ given in (2.8).

Let us first prove (2.9). Due to Proposition 1 and (2.7), we obtain

$$\begin{aligned} \mathbf{P}_B(|\mathcal{E}_h(\gamma)| \ge 4) &\ge \mathbf{P}_B(|\mathcal{E}_1(\gamma) \cap [-h^{-2}, h^{-2}]| \ge 4) \\ &\ge \mathbf{P}_B(S_1^{(1)} \le 1/h^2) P(X^{(1)} \le 1/(3h^2))^3. \end{aligned} \tag{B.1}$$

By (2.7), the Schwarz inequality and since $E(X^{(1)}) = 1/\sigma^2$, $E((X^{(1)})^2) = 1/2\sigma^2$, we obtain that, for all $t \geq 0$,

$$\begin{aligned}\mathbf{P}_B(S_1^{(1)} \leq t) &\geq 1 - E(X^{(1)}; X^{(1)} > t)/E(X^{(1)}) \\ &\geq 1 - E((X^{(1)})^2)^{1/2} P(X^{(1)} > t)^{1/2}/E(X^{(1)}) \\ &= 1 - \frac{\sigma}{\sqrt{2}} P(X^{(1)} > t)^{1/2}.\end{aligned} \tag{B.2}$$

Since $X^{(1)}$ has density $\sigma^2 f(\sigma^2 x)\,dx$ with $f(x)$ as in (2.6), for each $\alpha > 0$, there exists $c(\sigma, \alpha) > 0$ such that $P(X^{(1)} \geq t) \leq c(\sigma, \alpha) t^{-\alpha}$ for all $t > 0$. This allows to bound from below $P_B(S_1^{(1)} \leq 1/h^2)$ [due to (B.2)] and $P_B(X^{(1)} \leq 1/(3h^2))$. These lower bounds together with (B.1) imply (2.9).

In order to prove (2.10), we observe that Proposition 1 implies

$$\mathbf{P}_B(\gamma : |\mathcal{E}_h(\gamma) \cap [-1, 1]| \geq n) \leq P(Z_n^{(h)} \leq 2), \tag{B.3}$$

where

$$Z_n^{(h)} \equiv X_1^{(h)} + X_2^{(h)} + \cdots + X_n^{(h)}$$

and $X_1^{(h)}, X_2^{(h)}, \ldots, X_n^{(h)}$ are i.i.d. random variables having Laplace transform given by the r.h.s. of (2.8). In particular, for all $t > 0$,

$$P(Z_n^{(h)} \leq t) \leq eE(\exp\{-Z_n^{(h)}/t\}) = e/\cosh^n\left(\frac{\sqrt{2}h}{\sqrt{t}\sigma}\right) \leq e\left(1 + \frac{h^2}{t\sigma^2}\right)^{-n}, \tag{B.4}$$

where in the last inequality we have used the bound $\cosh x \geq 1 + x^2/2$. By taking $t = 2$, we get (2.10).

To prove (2.11), we define the increasing sequence $\tilde{S}_1^{(h)} < \tilde{S}_2^{(h)} < \cdots$ as the sequence of $h$-extrema of $\gamma$ not larger than $-1$ (note that such a sequence is well defined $\mathbf{P}_B$-almost surely). Then, the l.h.s. of (2.11) can be bounded by

$$\begin{aligned}\sum_{n=2}^{\infty} \mathbf{P}_B(|\mathcal{E}_h(\gamma) \cap [-1, 1]| = n, \exists j : 1 \leq j \leq n - 1, \\ \text{s.t. } |\gamma(\tilde{S}_j^{(h)}) - \gamma(\tilde{S}_{j+1}^{(h)})| < h + \delta).\end{aligned} \tag{B.5}$$

By Proposition 1, for all $n \in \mathbb{Z}$, $|\gamma(S_n^{(h)}) - \gamma(S_{n+1}^{(h)})| - h$ is an exponential variable with mean $h$ and, therefore,

$$\mathbf{P}_B(|\gamma(S_n^{(h)}) - \gamma(S_{n+1}^{(h)})| < h + \delta) \leq 1 - e^{-\delta/h}.$$

Therefore, due to the bound (B.5) and (2.10),

$$\text{l.h.s. of (2.11)} \leq (1 - e^{-\delta/h}) \sum_{n=2}^{\infty} e\left(1 + \frac{h^2}{2\sigma^2}\right)^{-n} n.$$

Since for all $a > 1$ $\sum_{n=1}^{\infty} na^{-n} = a(a-1)^{-2}$, we get (2.11).

To prove (2.12), given $\gamma \in C(\mathbb{R})$ and $a_1 < a_2 < \cdots < a_n$, we say that condition $C((a_1, a_2, \ldots, a_n), \gamma, \delta)$ is fulfilled if

$$||\gamma(a_i) - \gamma(a_j)| - |\gamma(a_{i'}) - \gamma(a_{j'})|| \geq \delta$$

for all $(i,j) \neq (i',j')$ with $i, i'$ odd, $j, j'$ even and $i < j$, $i' < j'$. Due to Proposition 1 and since $(B_{th^2}/h, t \in \mathbb{R}) \stackrel{\text{law}}{=} (B_t, t \in \mathbb{R})$,

$$\begin{aligned}
\text{l.h.s. of (2.12)} &\leq \sum_{n=1}^{\infty} \mathbf{P}_B(|\mathcal{E}_h(\gamma) \cap [-1,1]| = n, C((\tilde{S}_1^{(h)}, \tilde{S}_2^{(h)}, \ldots, \tilde{S}_n^{(h)}), \gamma, \delta)) \\
&\leq \sum_{n=1}^{\infty} \mathbf{P}_B(|\mathcal{E}_h(\gamma) \cap [-1,1]| = n) \\
&\qquad \times \sum_{a,b,a',b'}^{(n)} P(||\Sigma(a,b)| - |\Sigma(a',b')|| \leq \delta/h),
\end{aligned}$$

where the summation $\sum_{a,b,a',b'}^{(n)}$ is over all odd integers $a, b, a', b'$ with $1 \leq a \leq b \leq n$, $1 \leq a' \leq b' \leq n$ and, given $a \leq b$ odd with $1 \leq a \leq b \leq n$,

$$\Sigma(a,b) \equiv (Y_a - Y_{a+1}) + (Y_{a+2} - Y_{a+3}) + \cdots + (Y_{b-2} - Y_{b-1}) + Y_b + 1, \tag{B.6}$$

where $Y_z, z \in \mathbb{Z}$, are independent exponential variables with mean 1.

We claim that there exists a constant $c_0 > 0$, independent of all other parameters, such that

$$P(||\Sigma(a,b)| - |\Sigma(a',b')|| \leq \delta/h) \leq c_0 \delta/h \tag{B.7}$$

for all $a, b, a', b'$ as above. This, together with (B.3) and (B.4), implies

$$\text{l.h.s. of (2.12)} \leq \frac{ec_0\delta}{h} \sum_{n=1}^{\infty} n^4 \left(1 + \frac{h^2}{2\sigma^2}\right)^{-n}.$$

Since $\sum_{n=1}^{\infty} n^4 a^{-n} \leq ca^4/(a-1)^5$ for all $a > 1$, the above bound implies (2.12).

Let us prove (B.7). Since $E(\exp\{itX\}) = 1/(1-it)$ if $X$ is an exponential variable with mean 1, we obtain that the characteristic function $\phi_{a,b}(t) \equiv E(\exp\{it\Sigma(a,b)\})$ satisfies

$$|\phi_{a,b}(t)| \leq (1+t^2)^{-(b-a)/2} |1-it|^{-1}.$$

In particular, by the inverse formula of the Fourier transform, if $a < b$ or due the explicit expression if $a = b$, we get that $0 \leq f_{a,b}(x) \leq c'$ $\forall x \in \mathbb{R}$, where $f_{a,b}$ is the density function of $\Sigma(a,b)$ and the constant $c'$ is independent of all parameters. Let us first suppose that $a \le b < a' \le b'$ and bound

$$\begin{aligned} &P(||\Sigma(a,b)| - |\Sigma(a',b')|| \le \delta/h) \\ &\quad \le P(|\Sigma(a,b) - \Sigma(a',b')| \le \delta/h) \\ &\qquad + P(|\Sigma(a,b) + \Sigma(a',b')| \le \delta/h). \end{aligned} \tag{B.8}$$

Since $\Sigma(a,b), \Sigma(a'b')$ are independent,

$$\begin{aligned} &P(|\Sigma(a,b) - \Sigma(a',b')| \le \delta/h) \\ &\quad = \int_{\mathbb{R}} dx' \, f_{a',b'}(x') \int_{\mathbb{R}} dx \, f_{a,b}(x) \mathbb{I}_{|x-x'| \le \delta/h} \le 2c'\delta/h. \end{aligned}$$

Similarly, one can bound the last member in (B.8) by $2c'\delta/h$. It is easy to adapt the above argument when the sets $[a,b] \cap \mathbb{Z}$, $[a',b'] \cap \mathbb{Z}$ have nonempty intersection in order to get a similar bound for the l.h.s. of (B.8), completing the proof of (B.7).

Let us prove (2.13) by using Lemma 3. Note that this lemma gives the statistics of the $h$-slopes that are not crossing a given point, and therefore cannot be applied directly to the $h$-slope crossing $-1$ or 1. In order to avoid this problem, we look to the behavior of the $h$-slopes in a larger interval $[-L-1, L+1]$ requiring that the first $h$-extremum in such an interval is smaller than $-1$ and the last one is larger than 1 (in this way, all the $h$-slopes covering part of the interval $[-1,1]$ cannot cross the boundary $\{-L-1, L+1\}$). For any $\alpha > 0$, the probability that the previous condition is not satisfied can be bounded from above by

$$2\mathbf{P}_B(S_1^{(h)} > L) = 2\mathbf{P}_B(S_1^{(1)} > Lh^{-2}) \le c(\alpha,\sigma) L^{-\alpha} h^{2\alpha}$$

due to (B.2) and the subsequent discussion there.

Let us define $\mathcal{E}_{h,\beta,\varepsilon}$ as the event that $\exists n \in \mathbb{Z}$ with $S_n^{(h)}, S_{n+1}^{(h)} \in [-L-1, L+1]$ and

$$\Big( \inf_{t \in (\beta, T_{n,+}^{(h)}]} |B_{S_n^{(h)}+t} - B_{S_n^{(h)}}| \Big) \wedge \Big( \inf_{t \in (\beta, T_{n+1,-}^{(h)}]} |B_{S_{n+1}^{(h)}-t} - B_{S_{n+1}^{(h)}}| \Big) < \varepsilon.$$

Then, due to Lemma 3,

$$\begin{aligned} \mathbf{P}_B(\mathcal{D}_{h,\beta,\varepsilon}) &\le \mathbf{P}_B(\mathcal{E}_{h,\beta,\varepsilon}) + c(\alpha,\sigma) L^{-\alpha} h^{2\alpha} \\ &\le cn \sum_{n=2}^{\infty} \mathbf{P}_B(\gamma : |\mathcal{E}_h(\gamma) \cap [-L-1, L+1]| = n) \varepsilon/\sqrt{\beta} \\ &\quad + c(\alpha,\sigma) L^{-\alpha} h^{2\alpha}. \end{aligned} \tag{B.9}$$

By (2.3) and (2.10),

$$\begin{aligned}
&\mathbf{P}_B(\gamma : |\mathcal{E}_h(\gamma) \cap [-L-1, L+1]| = n)\\
&\qquad = \mathbf{P}_B(\gamma : |\mathcal{E}_{h/\sqrt{L+1}}(\gamma) \cap [-1,1]| = n)\\
&\qquad \le c\left(1 + \frac{h^2}{2\sigma^2(L+1)}\right)^{-n}.
\end{aligned}$$

Since for all $a > 1$, $\sum_{n=1}^{\infty} na^{-n} = a/(a-1)^2$, (2.13) follows from the above estimates by taking $\alpha \equiv 2$, $h^2/L = \varepsilon^{1/4}/\beta^{1/8}$.

To prove (2.14), we observe that $1 \le |M_h^-(\gamma^*)| \le Q$ if $4 \le |\mathcal{E}_h(\gamma) \cap [-1,1]| \le n$. Due to (2.9) and (2.10), by choosing $h$ small enough, the last event is verified with probability at least $1 - \alpha/5$. Let us assume that $M_h^-(\gamma^*) \ne \varnothing$. In order to verify conditions (1.15) and (1.16), we have to take in consideration that the smallest and the largest elements of $M_h^+(\gamma^*)$ could not be $h$-maxima of $\gamma$. By choosing $h'$ small enough, we have that $M_h^+(\gamma^*) \subset M_{h'}^+(\gamma)$ with probability at least $1 - \alpha/5$. In this case, condition (1.16) is implied by the event $(C_{h',\delta})^c$. Due to (2.12), $\mathbf{P}_B(C_{h',\delta}) < \alpha/5$ if $\delta$ is small enough. Similarly, due to (2.11), we can assume that the event $\mathcal{B}_{h,\delta}$ has probability less than $\alpha/5$ if $\delta$ is small enough. At this point, in order to verify condition (1.15), it remains to observe that if $\delta$ is small enough, then with probability at least $1 - \alpha/5$, one has $\gamma(w_1) - \gamma(u_1) > h + \delta$ and $\gamma(w_{q+1}) - \gamma(u_q) > h + \delta$, where $w_1, w_{q+1}, u_1, u_q$ are as in Lemma 2 with $\gamma$ replaced by $\gamma^*$.

## APPENDIX C: STURM OSCILLATION THEORY

As discussed in [19], the qualitative theory of second order Sturm–Liouville equations

$$\frac{d}{dt}\left(p(t)\frac{du}{dt}(t)\right) + q(t)u(t) = 0, \qquad p \in C^1, q \in C^0, p > 0$$

can be generalized to difference equations, that is, equations of the form $Hu = 0$ with $H$ a Jacobian matrix, namely, $H = (H_{i,j})_{i,j \in I}$ is a symmetric matrix indexed on a (possibly infinite) interval $I \subset \mathbb{Z}$ such that $H_{i,j} = 0$ whenever $|i - j| > 1$. In what follows we derive from [19] some results mainly related to Sturm oscillation theory for the Dirichlet operator $\mathbb{L}(D)$, $D \equiv \{a, a+1, \dots, b\} \subset \mathbb{Z}$. To this aim, we introduce the following notation: given $u \in \mathbb{R}^D$, the continuous function $\hat{u}$ is defined on $[a, b]$ by setting $\hat{u}(x) \equiv u(x)$ for all $x \in [a, b] \cap \mathbb{Z}$ and by extending $\hat{u}$ on $[a, b]$ by linear interpolation.

Let us first observe that, due to a simple iterative procedure, the system

$$((\mathbb{L}(D) - \lambda)u)(x) = 0 \qquad \forall x \in D \setminus \{b\},$$

uniquely determines $u \in \mathbb{R}^D$ when given the value $u(a)$ [in particular, the eigenvalues of $\mathbb{L}(D)$ are all simple] and each eigenvector cannot have two consecutive zeros and cannot vanish on $a$ or $b$. A deeper insight of the qualitative behavior of the eigenvectors is given by the following result:

Proposition 8 (Sturm oscillation theorem). *Let $\lambda_1 < \lambda_2 < \cdots < \lambda_r$ be the eigenvalues of $\mathbb{L}(D)$, where $D \equiv \{a, a+1, \dots, b\}$, $r = b - a + 1$. For each $1 \le i \le r$, let $f^{(i)}$ be an eigenvector of $\mathbb{L}(D)$ with eigenvalue $\lambda_i$. Then the function $\hat{f}^{(i)}$ has $i - 1$ zeros in $[a, b]$.*

Proof. Without loss of generality, we assume that $a = 1 < b = r$. Let us consider the matrix $H = (H_{i,j})_{i,j\in D}$ defined as $H_{i,j} \equiv (\mu(i)/\mu(j))^{1/2}\mathbb{L}_{i,j}$. Due to (1.12), $H$ is a Jacobian matrix. Moreover, since $H = A^{-1}\mathbb{L}(D)A$ where $A_{i,j} \equiv \delta_{i,j}\mu(j)^{-1/2}$, $f \in \mathbb{R}^D$ is an eigenvector of $H$ with eigenvalue $\lambda$ iff $Af$ is an eigenvector of $\mathbb{L}(D)$ with eigenvalue $\lambda$.

Given $\lambda \in \mathbb{R}$, let $\{u_j(\lambda)\}_{j\in D}$ be the unique solution of the system

$$\begin{cases} \sum_{j\in D} H_{i,j}u_j(\lambda) = \lambda u_i(\lambda) & \forall 1 \le i \le r-1, \\ u_1(\lambda) = 1. \end{cases}$$

By solving the above equations iteratively, one easily checks that $u_j(\lambda)$ is a polynomial of degree $j - 1$ with leading term $(a_1a_2\cdots a_{j-1})^{-1}\lambda^{j-1}$ for all $1 \le j \le r$, where $a_i \equiv H_{i,i+1} < 0$ for all $1 \le i \le r - 1$. Let us introduce the monomic polynomials

$$P_i(\lambda) \equiv \begin{cases} 1, & \text{if } i = 0, \\ (a_1a_2\cdots a_i)u_{i+1}(\lambda), & \text{if } 1 \le i < r, \\ \det(\lambda\mathbb{I} - H), & \text{if } i = r, \end{cases}$$

and define the function $\tilde{y}_\lambda(x)$ on $[0, r]$ by linear interpolation of the values $\tilde{y}_\lambda(i) \equiv (-1)^i P_i(\lambda)$, $i \in [0, r] \cap \mathbb{Z}$. Then, as stated after Proposition 2.4 in [19], the number of eigenvalues of $H$ below $\lambda$ equals the number of zeros of $\tilde{y}_\lambda$ on $[0, r)$.

If $\lambda = \lambda_k$ for some $1 \le k \le r$, then $\{u_j(\lambda)\}_{j\in D}$ is the unique eigenvector of $H$ with eigenvalue $\lambda$ such that $u_1(\lambda) = 1$. Moreover, $(-1)^i\mathrm{sgn}(P_i(\lambda)) = \mathrm{sgn}(u_{i+1}(\lambda))$ for all $0 \le i \le r-1$ since $a_1, \dots, a_{r-1}$ are negative, while $P_r(\lambda) = 0$ since $\lambda$ is an eigenvalue of $H$. In conclusion, the number of zeros of $\tilde{y}_\lambda$ on $[0, r)$ equals the number of zeros of the function $\hat{u}$ on $[1, r]$ defined by linear interpolation from the values $\hat{u}(i) \equiv u_i(\lambda)$, $i \in [1, r] \cap \mathbb{Z}$, which trivially equals the number of zeros of $\hat{f}^{(k)}$ on $[1, r]$. □

The above proposition and the observation that any eigenvector of $\mathbb{L}(D)$ cannot have two consecutive zeros easily imply the following result.

Corollary 2. *Let $\lambda_1 < \lambda_2 < \cdots < \lambda_r$ be the eigenvalues of $\mathbb{L}(D)$, where $D \equiv \{a, a+1, \dots, b\}$, $r = b - a + 1$. Given $1 \le i \le r$, let $f^{(i)}$ be an eigenvector of $\mathbb{L}(D)$ with eigenvalue $\lambda_i$. Then $f^{(1)}$ is of constant sign on $D$ while, for each index $i$ with $2 \le i \le r$, there exist integer numbers*

$$a \le y_1 < y_2 < \cdots < y_{i-1} < b$$

*such that $f^{(i)}$ is alternately nonnegative or negative on the $i$ intervals $[a, y_1] \cap \mathbb{Z}$, $[y_1 + 1, y_2] \cap \mathbb{Z}$, $[y_2 + 1, y_3] \cap \mathbb{Z}, \dots, [y_{i-1} + 1, b] \cap \mathbb{Z}$.*

A simple application of the above corollary is the following:

Corollary 3. *Let $A, B$ be finite subsets of $\mathbb{Z}$ with $A \subset B$ and $A \neq B$ and let $\lambda_A, \lambda_B$ be respectively the principal eigenvalue of $\mathbb{L}(A)$ and $\mathbb{L}(B)$. Then $\lambda_B < \lambda_A$.*

Proof. Let $f_A \in \mathbb{R}^A$ be a principal eigenvector of $\mathbb{L}(A)$ and let $\tilde{f}_A \in \mathbb{R}^B$ be defined as $\tilde{f}_A \equiv \mathbb{I}_A f_A$. Since $\mathbb{L}(B)\tilde{f}_A(x) = \lambda_A f_A(x)$ for all $x \in A$, we get

$$(\tilde{f}_A, \mathbb{L}(B)\tilde{f}_A)_{L^2(B,\mu)} = \lambda_A (\tilde{f}_A, \tilde{f}_A)_{L^2(B,\mu)}$$

and, consequently, $\lambda_B \leq \lambda_A$. Note that if $\lambda_B = \lambda_A$, then the above identity would imply that $\tilde{f}_A$ is proportional to $f_B$, in contradiction with Corollary 2. □

We can finally apply the Sturm oscillation theorem in order to show a spectral interlacing property for couples of Dirichlet operators.

Proposition 9. *Given points $a < z_1 < z_2 < \cdots < z_k < b$ in $\mathbb{Z}$, we define $D \equiv [a, b] \cap \mathbb{Z}$ and $D_k \equiv D \setminus \{z_1, \dots, z_k\}$. If $\gamma$ denotes the principle eigenvalue of $\mathbb{L}(D_k)$, then*

$$|\sigma(\mathbb{L}(D)) \cap [0, \gamma)| \leq k.$$

Proof. Let $\lambda_1 < \lambda_2 < \cdots < \lambda_r$ be the eigenvalues of $\mathbb{L}(D)$, where $r = b - a + 1 > k$, and let $f$ be an eigenvector of $\mathbb{L}(D)$ with eigenvalue $\lambda_{k+1}$. By the above corollary, there exist integers $a \leq y_1 < y_2 < \cdots < y_k < b$ such that $f$ is alternately nonnegative or negative on the intervals $[a, y_1] \cap \mathbb{Z}$, $[y_1 + 1, y_2] \cap \mathbb{Z}$, $[y_2 + 1, y_3] \cap \mathbb{Z}, \dots, [y_k + 1, b] \cap \mathbb{Z}$. Since these intervals are $k + 1$, at least one of them has empty intersection with $\{z_1, z_2, \dots, z_k\}$. Let us write such an interval as $[v, w] \cap \mathbb{Z}$, with $v, w \in \mathbb{Z}$, and let $j \in \{0, 1, \dots, k\}$ be such that $z_j < v \leq w < z_{j+1}$, where $z_0 \equiv a - 1$, $z_{k+1} = b + 1$. Finally, let us consider the Dirichlet operator $\mathbb{L}(I)$, $I \equiv (z_j, z_{j+1}) \cap \mathbb{Z}$, and denote by $\beta$ its principal eigenvalue and by $g$ a related eigenvector. Since $\mathbb{L}(D_k)\tilde{g} = \beta\tilde{g}$ where $\tilde{g} \in \mathbb{R}^{D_k}$ is defined as $\tilde{g} \equiv g\mathbb{I}_I$, it must be $\gamma \leq \beta$. In particular, the assertion follows if we prove that $\beta < \lambda_{k+1}$. Due to the variational characterization of $\beta$, in order to prove that $\beta \leq \lambda_{k+1}$, it is enough to show that

$$(h, \mathbb{L}(I)h)_{L^2(I,\mu)} \leq \lambda_{k+1} \sum_{v \leq x \leq w} \mu(x) f^2(x) = \lambda_{k+1}(h, h)_{L^2(I,\mu)}, \tag{C.1}$$

where $h \in \mathbb{R}^I$ is defined as $h \equiv f\mathbb{I}_{[v,w]}$. In fact, it is simple to check that $h$ is not the zero function, since in this case it should be $v = w$, $f(v) = 0$ and $f(v-1)$, $f(v+1)$ should be both negative or both positive. All this is in contradiction with the identity $\mathbb{L}(D)f(v) = \lambda_{k+1}f(v)$. In order to prove (C.1), we note that the identity there follows from the definition of $h$. To prove the inequality, we observe that $\mathbb{L}(I)h(x) = \lambda_{k+1}h(x)$ if $v < x < w$, while

$$\mathbb{L}(I)h(v) = \omega_v(f(v) - f(v+1)) = \mathbb{L}f(v) - (1-\omega_v)(f(v) - f(v-1)),$$

where we set $f \equiv 0$ outside $D$. Suppose, for example, that $f(v) \geq 0$, then $f(v-1) < 0$ because of the initial discussion. In particular, $f(v)(f(v) - f(v-1)) \geq 0$ and, therefore,

$$h(v)\mathbb{L}(I)h(v) \leq f(v)\mathbb{L}f(v) = \lambda_{k+1}f^2(v).$$

The same conclusion holds if $f(v) < 0$ and if we replace $v$ with $w$, thus proving (C.1). Finally, we note that if $\beta = \lambda_{k+1}$, then, due to (C.1) and the variational characterization of $\beta$, it should be

$$(h, \mathbb{L}(I)h)_{L^2(I,\mu)} = \beta(h,h)_{L^2(I,\mu)}.$$

It is simple to check that the above identity would imply that $g = ch$ on $I$ for some non zero constant $c$. Since by Corollary 2 $g$ cannot vanish, this would imply that $z_j + 1 = v$ and $z_{j+1} - 1 = w$. Moreover, the identities $\mathbb{L}(I)g(v) = \beta g(v)$, $\mathbb{L}(D)f(v) = \lambda_{k+1}f(v) = \beta f(v)$ and $g = cf$ on $I$ would imply that $f(v) = f(v-1)$ in contradiction with the property that $f(v)$ and $f(v-1)$ cannot have the same sign. This shows that $\beta < \lambda_{k+1}$, thus concluding the proof. □

**Acknowledgments.** We wish to thank J.-P. Bouchaud, C. Monthus and Z. Shi for very useful discussions. One of us (A.F.) thanks the Weierstrass Institute for Applied Analysis and Stochastics in Berlin for the kind hospitality while part of this work was being done.

Weierstrass Institut für Angewandte Analysis und Stochastik
Mohrenstrasse 39
10117 Berlin
Germany
and
Mathematisches Institut
Technische Universität Berlin
Strasse des 17. Juni 136
10623 Berlin
Germany
E-mail: bovier@wias-berlin.de

Dipartimento di Matematica "G. Castelnuovo"
Università "La Sapienza"
P.le Aldo Moro 2
00185 Roma
Italy
E-mail: faggiona@mat.uniroma1.it